\documentclass{elsart}

\usepackage{graphicx}
\usepackage{psfrag}
\usepackage{color}

\usepackage{amssymb}
\usepackage{amsmath}
\usepackage{dsfont}
\usepackage{rotating}

\newcommand{\K}{\mathbf{K}}
\newcommand{\M}{\mathbf{M}}

\newcommand{\U}{\mathcal{U}}

\renewcommand{\S}{\mathbf{S}} 
\newcommand{\Q}{\mathbf{Q}} 
\newcommand{\q}{\mathbf{q}} 
\newcommand{\f}{\mathbf{f}} 
\newcommand{\g}{\mathbf{g}}

\newcommand{\F}{\mathbf{F}} 
 
\newcommand{\x}{\mathbf{x}} 
\newcommand{\xxi}{\boldsymbol{\xi}} 
\newcommand{\X}{\mathbf{X}} 
\renewcommand{\P}{\mathbf{P}} 
\newcommand{\w}{\mathbf{w}}

\renewcommand{\v}{\mathbf{v}}

\begin{document}

\begin{frontmatter}



\title{Lagrangian ADER-WENO Finite Volume Schemes on Unstructured Triangular Meshes Based On Genuinely Multidimensional HLL Riemann Solvers}


\author[UNITN]{Walter Boscheri}
\ead{walter.boscheri@unitn.it}
\author[ND]{Dinshaw S. Balsara}
\ead{d.balsara@nd.edu}
\author[UNITN]{Michael Dumbser}
\ead{michael.dumbser@unitn.it}
\address[UNITN]{Laboratory of Applied Mathematics, Department of Civil, Environmental and Mechanical Engineering \\ University of Trento, Via Mesiano 77, I-38123 Trento, Italy}
\address[ND]{Physics Department, University of Notre Dame du Lac, 225 Nieuwland Science Hall, Notre Dame, IN 46556, USA}

\begin{abstract}
In this paper we use the genuinely multidimensional HLL Riemann solvers recently developed by Balsara et al. in \cite{BalsaraMultiDRS} to construct a new class 
of computationally efficient high order Lagrangian ADER-WENO one-step ALE finite volume schemes on unstructured triangular meshes. A nonlinear WENO reconstruction  
operator allows the algorithm to achieve high order of accuracy in space, while high order of 
accuracy in time is obtained by the use of an ADER time-stepping technique based on a local space-time  Galerkin  predictor. The multidimensional HLL and HLLC Riemann solvers operate at each vertex of the 
grid, considering the entire Voronoi neighborhood of each node and allows for larger time steps than  conventional one-dimensional Riemann solvers. 
The results produced by the multidimensional Riemann solver are then used twice in our one-step ALE algorithm: first, as a node solver that assigns a 
unique velocity vector to each vertex, in order to preserve the continuity of the computational mesh; second, as a building block for genuinely  multidimensional numerical flux evaluation that allows the scheme to run with \textit{larger time steps} compared to conventional finite volume schemes 
that use classical one-dimensional Riemann solvers in normal direction. 
The space-time flux integral computation is carried out at the boundaries of each triangular space-time control volume using the \textit{Simpson} 
quadrature rule in space and Gauss-Legendre quadrature in time. 
A rezoning step may be necessary in order to overcome element overlapping or crossing-over. Since our one-step ALE finite volume scheme is 
based directly on a space-time conservation formulation of the governing PDE system, the remapping stage is not needed, making our algorithm a so-called 
\textit{direct} ALE method. 

We apply the method presented in this article to two systems of hyperbolic conservation laws, namely the Euler equations of compressible gas dynamics and the 
equations of ideal classical magneto-hydrodynamics (MHD). Convergence studies up to fourth order of accuracy in space and time have been carried out.  
Several numerical test problems have been solved to validate the new approach. Furthermore, the new high order Lagrangian schemes based on genuinely  multidimensional Riemann solvers have been carefully compared with high order Lagrangian finite volume schemes based on conventional one-dimensional 
Riemann solvers. It has been clearly shown that due to the less restrictive CFL condition the new schemes based on multidimensional HLL and HLLC Riemann solvers 
are computationally more efficient than the ones based on a conventional one-dimensional Riemann solver technique. 
\end{abstract}
\begin{keyword}
Arbitrary-Lagrangian-Eulerian (ALE) \sep multidimensional HLL and HLLC Riemann solvers \sep large time steps  
\sep direct ALE \sep local rezoning \sep high order WENO finite volume schemes \sep moving unstructured meshes  \sep hyperbolic conservation laws \sep ADER schemes 
\sep Euler equations \sep MHD equations 
\end{keyword}
\end{frontmatter}


\section{Introduction}
\label{sec.introduction}
Hyperbolic systems of conservation laws describe mathematically many important phenomena, such as environmental flows, hydrodynamic and thermodynamic problems, 
as well as the dynamics of many industrial and mechanical processes. Therefore a lot of research has been carried out in the past decades in order to  solve those  conservation laws numerically, starting from the one-dimensional case. A very famous and widespread approach is given by \textit{Godunov}-type  finite volume methods \cite{GodunovRS,vanLeerRS}, where the discrete solution is stored as constant data within each control volume of the computational mesh and is evolved in  time by using the integral form of the conservation law. Since the discrete solution in general exhibits jumps at the element interfaces, the introduction of  numerical fluxes across the discontinuities of each cell is necessary. Godunov suggested to obtain these numerical fluxes by solving \textit{Riemann problems}
at each interface. Early work regarded the exact solution of the Riemann problem \cite{GodunovRS,ColellaRS}, that was followed by the development of approximate Riemann  solvers, such as the linearized Riemann solver of Roe \cite{Roe}, the HLL and HLLE Riemann solvers \cite{HartenHLL,EinfeldtHLL} and the local Lax-Friedrichs (LLF) solver proposed by  Rusanov \cite{RusanovRS}, which can be reinterpreted as an HLL-type flux with a particular choice of the signal speeds. While the above-mentioned HLL schemes are very robust, they smear out contact discontinuities. An  improvement was  made by Einfeldt and Munz in \cite{munz91} with the introduction of the HLLEM Riemann solver, where the intermediate state was assumed piecewise linear instead of piecewise constant. Another well-known improvement of the original HLL scheme is due to Toro et al. in \cite{ToroHLLC} with the design of the HLLC Riemann solvers  that use an enhanced wave model that is able to capture also the intermediate contact wave. In \cite{OsherRS} Osher et al. introduced a class of approximate Riemann solvers based on path integrals, where the paths were obtained by an approximation of the solution of the Riemann problem by rarefaction fans. A simpler and more general version of the Osher flux has recently been forwarded by Dumbser and Toro in \cite{OsherUniversal,OsherNC}. 
All those one-dimensional Riemann solvers can be used even in two- and  three-dimensional problems, where the discontinuities are resolved at each boundary of the control volume along the normal direction. As shown in detail in \cite{ToroBook} the stability of any unsplit Godunov-type finite volume scheme using one-dimensional Riemann solvers in $d$ space dimensions is guaranteed under a  CFL condition of the type $\textnormal{CFL} < 1/d$ that becomes the more severe as the dimensionality of the problem increases. 

The increasingly severe restriction on the timestep in multiple space dimensions caused a lot of effort in research with the aim to introduce multidimensional  effects into the Riemann solvers \cite{Rumsey1993,Colella1990,WAF}. 
Further advances have been made in \cite{AbgrallmultiD1,AbgrallmultiD2,AbgrallmultiD3,Fey1,Fey2}, where new multidimensional Riemann solvers were  developed in extension of the linearized Riemann solver of Roe \cite{Roe}. Other schemes that consider multi-dimensional effects are the well-known flux-corrected transport (FCT) algorithm \cite{FCT1,FCT2} and the so-called finite volume evolution Galerkin method proposed in 
\cite{EvolutionGalerkin1,EvolutionGalerkin2}. 
In \cite{balsarahlle2d,balsarahllc2d} Balsara designed the first genuinely multidimensional HLL and HLLC Riemann solvers in two space dimensions for hyperbolic conservation laws on Cartesian grids. There, several applications for both hydrodynamics and  magneto-hydrodynamics (MHD) were shown in order to validate the method. In \cite{balsarahllc2d} it was shown that multid RS could be built on the following four precepts: 1) a self-similar wave model, 2) entropy enforcement, 
3) consistency with the conservation law and 4) preservation of internal sub-structures, like the contact discontinuity. In \cite{BalsaraMultiDRS} these precepts 
have been carried over to unstructured triangular meshes for the Euler and MHD equations.  

In the Lagrangian framework the mesh is moving together with the fluid, in order to identify material interfaces and to track them precisely. The use of Godunov-type methods is widespread among Lagrangian schemes. In \cite{munz94} Munz proposed the first Godunov-type finite volume schemes for Lagrangian gas dynamics based  on Roe and HLL-type Riemann solvers. Multi-dimensional cell-centered Lagrangian finite volume schemes have been considered by Despr\'es et al. in 
\cite{DepresMazeran2003,Despres2005,Despres2009} while ALE schemes with remapping for single and multi-material flows have been considered by the group of 
Shashkov et al. at Los Alamos in a recent series of papers \cite{ShashkovMultiMat1,ShashkovMultiMat2,ShashkovMultiMat3,ShashkovMultiMat4,ShashkovCellCentered,ShashkovRemap1,ShashkovRemap2,ShashkovRemap3,ShashkovRemap4,ShashkovRemap5,LoubereMaireShashkov}. 
Maire et al. \cite{Maire2009,Maire2010,Maire2011} introduced first and second order accurate cell-centered Lagrangian schemes in two- and three- space dimensions on general polygonal grids, computing the time derivatives of the fluxes with a node-centered solver that can be interpreted as a multi-dimensional extension of the Generalized Riemann Problem (GRP) methodology introduced by Ben-Artzi and Falcovitz \cite{Artzi}, Le Floch et al. \cite{Raviart.GRP.1,Raviart.GRP.2} and Titarev and Toro \cite{toro3,TitarevToro,Toro:2006a}. 

In the finite element framework Scovazzi et al. \cite{scovazzi1,scovazzi2} constructed better than second order accurate Lagrangian schemes, while better than second order accurate Lagrangian finite volume schemes have been proposed for the first time by Cheng and Shu \cite{chengshu1,chengshu2,chengshu3,chengshu4}, who  used an essentially non-oscillatory (ENO) reconstruction on curved \textit{structured} meshes to achieve high order of accuracy in space. Through the use of a 
Runge-Kutta or a Lax-Wendroff-type time stepping procedure also high order of accuracy in time has been obtained. Very recently, Dumbser et al. 
\cite{Lagrange1D,BoscheriDumbserLag,DumbserBoscheriLagNC} developed the first high order Lagrangian one-step ADER-WENO finite volume schemes for conservative and non-conservative hyperbolic systems on \textit{unstructured} triangular meshes. 
In \cite{LagNS} the multidimensional HLL Riemann solver presented in \cite{BalsaraMultiDRS} has been used as a node solver for the computation of the mesh velocity while in this paper we introduce the multidimensional HLL Riemann solver \cite{BalsaraMultiDRS} for the first time also into the  
space-time flux computation of a high order Lagrangian finite volume scheme, leading to a more efficient algorithm that can deal with CFL numbers close to 
unity even in two space dimensions.
The method presented in this paper belongs to the class of Arbitrary Lagrangian Eulerian (ALE) schemes \cite{Hirt1974,Peery2000,Smith1999,Feistauer1,Feistauer2,Feistauer3,Feistauer4}, where the mesh velocity can be chosen arbitrarily and does not necessarily have to coincide with the local fluid velocity. 

The rest of this paper is structured as follows. The numerical scheme is presented in Section \ref{sec.numethod}, while in Section \ref{sec.validation} numerical convergence studies are shown and some typical benchmark problems for hydrodynamics and magnetohydrodynamics are solved. Concluding remarks and an outlook to  future research and developments are given in Section \ref{sec.concl}. 

\section{Numerical Method}
\label{sec.numethod} 
In this paper we consider nonlinear systems of hyperbolic balance laws, which can generally be cast into the following form: 
\begin{equation}
\label{PDE}
  \frac{\partial \Q}{\partial t} + \nabla \cdot \F(\Q) = \S(\Q), \qquad \mathbf{x} \in \Omega(t) \subset \mathbb{R}^2, \quad t \in \mathbb{R}_0^+,      
\end{equation} 
where $t$ is the time and $\mathbf{x}=(x,y)$ represents the spatial position vector. The vector of conserved variables is represented by $\Q=(q_1,q_2,...,q_\nu)
\in \Omega_{\Q}$ and is defined in the space of the admissible states $\Omega_{\Q} \subset \mathbb{R}^\nu$, the nonlinear flux tensor is given by $\F(\Q)=\left( \f(\Q),\g(\Q) \right)$ and $\S(\Q)$ is a nonlinear but non-stiff algebraic source term. The two-dimensional time-dependent computational domain is denoted by $\Omega(t) \subset \mathds{R}^2$ and it is discretized by a total number $N_E$ of conforming triangles $T^n_i$ at a general time $t^n$. The \textit{current triangulation} $\mathcal{T}^n_{\Omega}$ of the domain $\Omega(t^n)=\Omega^n$ is simply the union of all the elements of the domain at a given time and it can be written as
\begin{equation}
\mathcal{T}^n_{\Omega} = \bigcup \limits_{i=1}^{N_E}{T^n_i}.  
\label{trian}
\end{equation}

Since the mesh is unstructured and we are working in the Lagrangian framework, which implies mesh motion and element deformation, it is convenient to introduce a \textit{local} spatial reference coordinate system $\xi-\eta$, where the physical element $T_i^n$ is mapped to a unit reference element $T_e$. The vector of spatial coordinates in the physical system is denoted by $\mathbf{x}=(x,y)$, while $\boldsymbol{\xi} = (\xi, \eta)$ is the position vector in the reference system. The unit triangle is defined by the nodes $\boldsymbol{\xi}_{e,1}=(\xi_{e,1},\eta_{e,1})=(0,0)$, $\boldsymbol{\xi}_{e,2}=(\xi_{e,2},\eta_{e,2})=(1,0)$ and $\boldsymbol{\xi}_{e,3}=(\xi_{e,2},\eta_{e,2})=(0,1)$. The spatial mapping reads 
\begin{equation} 
 \x = \x(\xxi,t^n) = \X^n_{1,i} + \left( \X^n_{2,i} - \X^n_{1,i} \right) \xi + \left( \X^n_{3,i} - \X^n_{1,i} \right) \eta, 
 \label{xietaTransf} 
\end{equation} 
with $\mathbf{X}^n_{k,i} = (X^n_{k,i},Y^n_{k,i})$ representing the vector of physical coordinates of the $k$-th vertex of triangle $T^n_i$ at time $t^n$.

We adopt a finite volume scheme, where data are represented as usual by piecewise constant cell averages and they are stored and evolved in time within the control volumes. Hence, the solution vector is given for each element $T_i^n$ by
\begin{equation}
  \Q_i^n = \frac{1}{|T_i^n|} \int_{T^n_i} \Q(\x,t^n) d\x,     
 \label{eqn.cellaverage}
\end{equation}
where $|T_i^n|$ is the volume of element $T_i^n$ at the current time $t^n$, which reduces to the surface of $T_i^n$ in two space dimensions. If only the cell averages \eqref{eqn.cellaverage} are used inside the numerical scheme, the resulting algorithm will only be first order accurate in space and time. Therefore, 
we use a high order WENO reconstruction technique to achieve higher order of accuracy in space. Within this reconstruction or recovery step, piecewise high 
order polynomials $\mathbf{w}_h(\x,t^n)$ are recovered from the known cell averages, as described in the next section.

\subsection{Polynomial WENO Reconstruction} 
\label{sec.weno} 
The original WENO scheme \cite{JiangShu1996,HuShuVortex1999,ZhangShu3D,BalsaraShu2000} uses a \textit{pointwise} formulation, while here we adopt the \textit{polynomial} approach presented in \cite{StencilRec1990,friedrich,KaeserIske2005,Dumbser2007204,DumbserKaeser06b,MOOD,MixedWENO2D,MixedWENO3D}. Since all the details of the algorithm can be found in the above-mentioned references, we only give a brief overview and description of the reconstruction procedure in the following. 

The reconstructed solution $\mathbf{w}_h(\x,t^n)$ is represented by piecewise polynomials of degree $M$ and is computed at each time level $t^n$ for each control volume $T_i^n$ of the domain $\Omega^n$. We need a set of \textit{reconstruction stencils} $\mathcal{S}_i^s$, each of which is composed of a total number of $n_e$ elements belonging to some neighborhood of $T_i^n$, i.e. 
\begin{equation}
\mathcal{S}_i^s = \bigcup \limits_{j=1}^{n_e} T^n_{m(j)}, 
\label{stencil}
\end{equation}
where $s$ denotes the stencil number. As stated in \cite{StencilRec1990,Olliver2002,KaeserIske2005}, the number of elements inside the stencil must be greater than the smallest number $\mathcal{M} = (M+1)(M+2)/2$ needed to reach the formal order of accuracy $M+1$, hence we typically take $n_e = 2 \mathcal{M}$ in two space dimensions. In \eqref{stencil} those elements that belong to the stencil are counted by the local index $1\leq j \leq n_e$, which is mapped to the global element number $m(j)$ of the triangulation \eqref{trian}. In order to avoid ill-conditioned reconstruction matrices, the reconstruction is performed in the \textit{reference system} $(\xi,\eta)$ according to the mapping \eqref{xietaTransf}.

The reconstruction polynomials $\w^s_h(x,y,t^n)$ are then expressed using a modal basis in terms of the \textit{orthogonal} Dubiner-type basis functions  $\psi_l(\xi,\eta)$ described in \cite{Dubiner,orth-basis,CBS-book} and they read  
\begin{equation}
\label{eqn.recpolydef} 
\w^s_h(x,y,t^n) = \sum \limits_{l=1}^\mathcal{M} \psi_l(\xi,\eta) \hat \w^{n,s}_{l,i} := \psi_l(\xi,\eta) \hat \w^{n,s}_{l,i}.   
\end{equation}

For each element $T^n_j \in \mathcal{S}_i^s$ the reconstruction is based on integral conservation, hence 
\begin{equation}
\label{intConsRec}
\frac{1}{|T^n_j|} \int \limits_{T^n_j} \psi_l(\xxi) \hat \w^{n,s}_{l,i} d\x = \Q^n_j, \qquad \forall T^n_j \in \mathcal{S}_i^s,      
\end{equation}
where the multi-dimensional integrals are evaluated using Gaussian quadrature formulae of suitable order, see \cite{stroud} for details. To make notation easier we will use in the rest of the paper classical tensor index notation with the Einstein summation convention, which implies summation over two equal indices. $\mathcal{M}$ is the total number of unknown degrees of freedom which has to be determined for each element $T_i^n$.

We remind that the number of stencil elements $n_e$ is greater than the number of unknowns $\hat \w^{n,s}_{l,i}$, so that Eqn. \eqref{intConsRec} yields an \textit{over-determined} linear algebraic system that is solved using either a constrained least-squares (LSQ) technique, see 
\cite{StencilRec1990,KaeserIske2005,DumbserKaeser06b}, or a more robust  singular value decomposition algorithm (SVD). 
The linear constraint requires the integral conservation equation \eqref{intConsRec} to hold exactly at least for element $T_i^n$. 
We note that as an alternative to LSQ and SVD the more flexible and more elegant kernel reconstruction recently proposed by Aboiyar et al. 
\cite{AboiyarIske} can be used, which automatically satisfies the constraint by construction. 
Since we are dealing with a \textit{moving} mesh algorithm and the integrals of \eqref{intConsRec} depend on the geometry, we can not compute a reconstruction matrix and store it once and for all elements in a preprocessing stage, but the system must be assembled and solved again at the beginning of each time step. However, this can be efficiently done using optimized standard LAPACK and BLAS routines. What remains constant during the whole computation is the definition 
of the reconstruction stencils, since the stencil search algorithm is rather expensive and should be called only once in the preprocessor stage of the 
algorithm before starting the simulation.  

As stated by the Godunov theorem \cite{GodunovRS}, no linear monotone schemes of order greater than one can exist, therefore the reconstruction scheme must be nonlinear in order to circumvent the theorem and to achieve higher order of accuracy. Hence, more than just one reconstruction stencil is needed, i.e. $s>1$ in Eqn. (\ref{stencil}), and for each stencil a reconstruction polynomial $\w^s_h$ is computed. In 2D we use 
in total 7 reconstruction stencils for each element, hence $1 \leq s \leq 7$: specifically we take one central stencil ($s=1$), three forward stencils ($2 < s \leq 4$) and three backward stencils ($5 < s \leq 7$), as proposed in \cite{KaeserIske2005,DumbserKaeser06b}. The final nonlinear WENO reconstruction polynomial is given as a weighted combination of the 7 reconstruction polynomials defined for each stencil, where the nonlinearity is inserted into the WENO weights, which also depend on $\w^s_h(x,y,t^n)$. We use the oscillation indicators $\sigma_s$ defined in \cite{JiangShu1996} and the oscillation indicator matrix $\Sigma_{lm}$ proposed in \cite{Dumbser2007204,DumbserKaeser06b}, which read  
\begin{equation}
\sigma_s = \Sigma_{lm} \hat w^{n,s}_{l,i} \hat w^{n,s}_{m,i}, \qquad 
\Sigma_{lm} = \sum \limits_{ \alpha + \beta \leq M}  \, \, \int \limits_{T_e} \frac{\partial^{\alpha+\beta} \psi_l(\xi,\eta)}{\partial \xi^\alpha \partial \eta^\beta} \cdot 
                                                                         \frac{\partial^{\alpha+\beta} \psi_m(\xi,\eta)}{\partial \xi^\alpha \partial \eta^\beta} d\xi d\eta.   
\end{equation} 
The nonlinear weights $\omega_s$ are defined by
\begin{equation}
\tilde{\omega}_s = \frac{\lambda_s}{\left(\sigma_s + \epsilon \right)^r}, \qquad 
\omega_s = \frac{\tilde{\omega}_s}{\sum_q \tilde{\omega}_q},  
\end{equation} 
where we set $\epsilon=10^{-14}$, $r=8$, $\lambda_s=1$ for the one-sided stencils and $\lambda_0=10^5$ for the central stencil, as done in \cite{DumbserEnauxToro,Dumbser2007204}. The final nonlinear WENO reconstruction polynomial and its coefficients are then given by 
\begin{equation}
\label{eqn.weno} 
 \w_h(x,y,t^n) = \sum \limits_{l=1}^{\mathcal{M}} \psi_l(\xi,\eta) \hat \w^{n}_{l,i}, \qquad \textnormal{ with } \qquad  
 \hat \w^{n}_{l,i} = \sum_s \omega_s \hat \w^{n,s}_{l,i}.   
\end{equation}

\subsection{Local Space-Time Predictor on Moving Curved Triangular Meshes} 
\label{sec.localCG} 
The finite volume algorithm presented in this paper is required to be high order accurate also in time, therefore the reconstructed polynomials $\w_h$ obtained at  the current time $t^n$ are then \textit{evolved} locally within each element $T_i(t)$ during one time step $[t^{n};t^{n+1}]$. This local evolution step is performed for each element and it leads to piecewise space-time polynomials of degree $M$, denoted by $\q_h(\x,t)$, and it was first introduced by Dumbser et al. in \cite{DumbserEnauxToro,Dumbser20088209} for the Eulerian case and then extended to moving meshes in \cite{Lagrange1D,BoscheriDumbserLag,DumbserBoscheriLagNC}. Such a procedure is carried on without any information from neighbor elements, hence improving the efficiency of the algorithm. As shown in \cite{BalsaraRumpf}, the local space--time predictor technique adopted in this paper, also known as ADER scheme, is almost two times more efficient than the strong stability preserving Runge-Kutta time stepping schemes. Only later in the finite volume scheme, where flux computation occurs, we couple the information with neighbor data. 

A \textit{weak} formulation of the governing PDE \eqref{PDE} is used to obtain the space-time solution $\q_h(\x,t)$. First we rewrite the balance law \eqref{PDE} in the local reference system, i.e. 
\begin{equation}
\frac{\partial \Q}{\partial \tau}\tau_t + \frac{\partial \Q}{\partial \xi}\xi_t + \frac{\partial \Q}{\partial \eta}\eta_t + \frac{\partial \f}{\partial \tau}\tau_x + \frac{\partial \f}{\partial \xi}\xi_x + \frac{\partial \f}{\partial \eta}\eta_x + \frac{\partial \g}{\partial \tau}\tau_y + \frac{\partial \g}{\partial \xi}\xi_y + \frac{\partial \g}{\partial \eta}\eta_y = \mathbf{S}(\Q),  
\label{PDEweak}
\end{equation}
where $\mathbf{x}=(x,y)$ and $\boldsymbol{\xi}=(\xi,\eta)$ represent the spatial coordinate vectors in physical and reference coordinates, respectively, while $\mathbf{\tilde{x}}=(x,y,t)$ and $\boldsymbol{\tilde{\xi}}=(\xi,\eta,\tau)$ are the corresponding space-time coordinates. The mapping in time, which is simply given by  
\begin{equation}
t = t_n + \tau \, \Delta t, \qquad  \tau = \frac{t - t^n}{\Delta t}, \qquad \Rightarrow \qquad \widehat{t}_l = t_n + \tau_l \, \Delta t, 
\label{timeTransf}
\end{equation} 
together with the local space transformation \eqref{xietaTransf}, are used to define the following Jacobian matrix and its inverse that are needed to formulate Eqn. \eqref{PDEweak}:
\begin{equation}
J_{st} = \frac{\partial \mathbf{\tilde{x}}}{\partial \boldsymbol{\tilde{\xi}}} = \left( \begin{array}{ccc} x_{\xi} & x_{\eta} & x_{\tau} \\ y_{\xi} & y_{\eta} & y_{\tau} \\ 0 & 0 & \Delta_t \\ \end{array} \right), \quad J_{st}^{-1} = \frac{\partial \boldsymbol{\tilde{\xi}}}{\partial \mathbf{\tilde{x}}} = \left( \begin{array}{ccc} \xi_{x} & \xi_{y} & \xi_{t} \\ \eta_{x} & \eta_{y} & \eta_{t} \\ 0 & 0 & \frac{1}{\Delta t} \\ \end{array} \right). 
\label{Jac}
\end{equation}
Here, we introduced the simplifications $\tau_x = \tau_y = 0$ and $\tau_t = \frac{1}{\Delta t}$, according to the definition (\ref{timeTransf}). With the inverse of the Jacobian matrix \eqref{Jac} the weak form given by \eqref{PDEweak} reduces to 
\begin{equation}
\frac{\partial \Q}{\partial \tau} + \Delta t \left[ \frac{\partial \Q}{\partial \xi}\xi_t + \frac{\partial \Q}{\partial \eta}\eta_t + \frac{\partial \f}{\partial \xi}\xi_x + \frac{\partial \f}{\partial \eta}\eta_x + \frac{\partial \g}{\partial \xi}\xi_y + \frac{\partial \g}{\partial \eta}\eta_y  \right] = \Delta t \mathbf{S}(\Q),
\label{PDECG}
\end{equation}
which can be further abbreviated by 
\begin{equation}
\Q_\tau =  \Delta t \P,
\label{PCG}
\end{equation}
using the term $\P$, which reads
\begin{equation}
 \P := \S(\Q) - \left( \Q_{\xi}\xi_t + \Q_{\eta}\eta_t + \f_{\xi}\xi_x + \f_{\eta}\eta_x + \g_{\xi}\xi_y + \g_{\eta}\eta_y \right).
 \label{eqn.pdef}
\end{equation}

In order to discretize Eqn. \eqref{PDECG}, we rely on a set of space-time basis functions $\theta_l=\theta_l(\boldsymbol{\tilde{\xi}})=\theta_l(\xi,\eta,\tau)$ which are defined by the Lagrange interpolation polynomials passing through a set of space-time nodes $\boldsymbol{\tilde{\xi}}_m=(\xi_m,\eta_m,\tau_m)$, see \cite{Dumbser20088209} for details. A nodal approach is then adopted to express the space-time solution $\q_h$ and the term $P_h$ as
\begin{equation}
\q_h=\q_h(\xi,\eta,\tau) = \theta_{l}(\xi,\eta,\tau) \widehat{\q}_{l,i}, \qquad \P_h=\P_h(\xi,\eta,\tau) = \theta_{l}(\xi,\eta,\tau) \hat{\P}_{l,i}. 
\label{thetaSol}
\end{equation}
In this article an \textit{isoparametric} approach is adopted, hence mapping the physical space-time coordinate vector $\mathbf{\tilde{x}}$ to the reference space-time coordinate vector $\boldsymbol{\tilde{\xi}}$ using the \textit{same} basis functions $\theta_l$ which represent also the solution $\q_h$. Therefore one obtains
\begin{equation}
 \x(\xi,\eta,\tau) = \theta_l(\xi,\eta,\tau) \widehat{\x}_{l,i}, \qquad t(\xi,\eta,\tau) = \theta_l(\xi,\eta,\tau) \widehat{t}_l, 
 \label{eqn.isoparametric} 
\end{equation} 
where the $\widehat{\mathbf{x}}_{l,i} = (\widehat{x}_{l,i},\widehat{y}_{l,i})$ denote the degrees of freedom of the vector of physical coordinates in space, while  the degrees of freedom $\widehat{t}_l$ denote the physical time at each space-time node $\tilde{\x}_{l,i} = (\widehat{x}_{l,i}, \widehat{y}_{l,i}, \widehat{t}_l)$. The spatial degrees of freedom are partially unknown, whereas the temporal degrees of freedom are \textit{known}. 

To make notation easier, let us introduce the following two integral operators: 
\begin{equation}
\left[f,g\right]^{\tau} = \int \limits_{T_e} f(\xi,\eta,\tau) g(\xi,\eta,\tau) d\xi d\eta, \quad \left\langle f,g \right\rangle = \int \limits_{0}^{1} \int \limits_{T_e} f(\xi,\eta,\tau)g(\xi,\eta,\tau) d\xi d\eta d\tau,  
\label{intOperators}
\end{equation}
which denote the scalar products of two functions $f$ and $g$ over the spatial reference element $T_e$ at time $\tau$ and over the space-time reference element $T_e\times \left[0,1\right]$, respectively. Let us furthermore introduce the following matrices,  
\begin{equation}
\K_{\tau} = \left\langle \theta_k,\frac{\partial \theta_l}{\partial \tau} \right\rangle, \qquad \M = \left\langle \theta_k,\theta_l \right\rangle. 
\label{Ktau}
\end{equation}
which do not depend on the current geometry configuration and therefore can be precomputed and stored once and for all in a preprocessing step. 

Next, using the definitions given in \eqref{thetaSol}, the weak formulation of the governing PDE \eqref{PDE} is obtained by multiplying \eqref{PDECG} with a test function which is given by the same space-time basis functions $\theta_k(\xi,\eta,\tau)$ and then integrating it over the unit reference space-time element $T_e \times [0,1]$. In a compact matrix notation it reads
\begin{equation}
\K_{\tau}\widehat{\q}_{l,i} = \Delta t \M \widehat{\P}_{l,i}. 
\label{LagrSTPDECGmatrix}
\end{equation}
The vector $\widehat{\q}_{l,i}$ can be split into two parts, yielding 
\begin{equation}
\widehat{\q}_{l,i}=(\widehat{\q}_{l,i}^{0},\widehat{\q}_{l,i}^{1}),
\label{eq.qhDOF}
\end{equation}
where $\widehat{\q}_{l,i}^{0}$ represent the degrees of freedom that are known from the initial condition $\w_h$ by setting the corresponding degrees of freedom to the known values, see \cite{Dumbser20088209} for details, while $\widehat{\q}_{l,i}^{1}$ are the unknown degrees of freedom for $\tau>0$. The known degrees of freedom $\widehat{\q}_{l,i}^{0}$ are moved onto the right-hand side of \eqref{LagrSTPDECGmatrix}, hence obtaining the following nonlinear algebraic equation system \eqref{PCG}, which can be solved by an iterative procedure, i.e.
\begin{equation}
  \K_{\tau} \widehat{\q}_{l,i}^{r+1} = \Delta t \M \widehat{\P}_{l,i}^r,  
\label{CGfinal}
\end{equation}
with the superscript $r$ denoting the iteration number. For an efficient initial guess ($r=0$) based on a second order MUSCL-type scheme see \cite{HidalgoDumbser}, otherwise one can simply take the reconstruction polynomial $\w_h$ at the initial time level.

In the ALE formulation the mesh is moving in time, hence the space-time control volume of each element $T_i^n$ is also changing within each timestep, i.e. the vertex coordinates of the local space-time element are evolving in time. Therefore we also have to solve the following ODE system 
\begin{equation}
\frac{d \mathbf{x}}{dt} = \mathbf{V}(\Q,\x,t),
\label{ODEmesh}
\end{equation}
which governs the element motion. Here, $\mathbf{V}=\mathbf{V}(\Q,\x,t)$ is the local mesh velocity, which does not necessarily have to coincide with the local fluid velocity, since our numerical algorithm is designed to be an \textit{arbitrary Lagrangian-Eulerian} scheme, where the mesh velocity can be chosen independently from the physical flow motion. That allows us to obtain either purely Lagrangian schemes, if the local mesh velocity is equal to the local fluid velocity, or purely   Eulerian schemes in the case of $\mathbf{V}=0$. The velocity $\mathbf{V}$ is also approximated using a nodal approach, hence yielding 
\begin{equation}
\mathbf{V}_h=\mathbf{V}_h(\xi,\eta,\tau) = \theta_{l}(\xi,\eta,\tau) \widehat{\mathbf{V}}_{l,i}, \quad  \widehat{\mathbf{V}}_{l,i} = \mathbf{V}(\widehat{\q}_{l,i}, \widehat{\x}_{l,i},\hat t_l).
\label{Vdof}
\end{equation}
As suggested in \cite{Lagrange1D,BoscheriDumbserLag}, the local space-time Galerkin method is also used to solve the system \eqref{ODEmesh} for the unknown  coordinate vector $\widehat{\mathbf{x}}_{l,i}$: 
\begin{equation}
\left\langle \theta_k,\frac{\partial \theta_l}{\partial \tau} \right\rangle \widehat{\mathbf{x}}_{l,i} = \Delta t \left\langle \theta_k,\theta_l \right\rangle \widehat{\mathbf{V}}_{l,i},
\label{VCG}
\end{equation}
hence obtaining the following iteration scheme for the vertex coordinates of the local element $T_i^n$:
\begin{equation}
\K_{\tau} \widehat{\mathbf{x}}^{r+1}_{l,i} = \Delta t \M \widehat{\mathbf{V}}^r_{l,i}.
\label{newVertPos}
\end{equation}
The initial condition of the ODE system is given by the nodal degrees of freedom $\widehat{\x}_l$ at relative time $\tau=0$, which are known from the current configuration of triangle $T_i^n$ at time $t^n$ and the mapping \eqref{xietaTransf}. Eqn. \eqref{newVertPos} is iterated \textit{together} with Eqn. \eqref{CGfinal}. The iteration stops when the residuals of both systems are less than a prescribed tolerance. 
Note that Eqn. \eqref{newVertPos} gives a high order accurate predictor of the mesh motion since all integrals present in \eqref{VCG} are evaluated with 
high order of accuracy $M$. 

Once we have carried out the above iterative procedure for all elements of the computational domain, we end up with an \textit{element-local predictor} for the numerical solution $\q_h$, for the fluxes $\mathbf{F}_h=(\f_h,\g_h)$, for the source term $\S_h$ and also for the mesh velocity $\mathbf{V}_h$. Since the space-time predictor procedure has been performed \textit{locally}, the predicted geometry at the new time level $t^{n+1}$ may be discontinuous. Therefore in the next Section \ref{sec.meshmotion} we will show how to resolve this discontinuity by using a multi-dimensional HLL Riemann solver to obtain a \textit{unique} velocity vector for each vertex of the computational mesh.

\subsection{Mesh motion}
\label{sec.meshmotion}
At the end of the local space-time predictor step the coordinate vector of node $k$ has been computed separately for each surrounding element $T_j^n \in \mathcal{V}_k$, where $\mathcal{V}_k$ denotes the Voronoi neighborhood of vertex $k$. Hence, node $k$ is in principle assigned with several node velocity vectors that would lead to different positions of the node at the new time level $t^{n+1}$. 
We therefore need a \textit{local} node-based strategy that defines a \textit{unique} time-averaged velocity vector $\overline{\mathbf{V}}_k^n$ at each mesh 
node $k$, so that the new vertex position can be simply computed as 
\begin{equation} 
	\mathbf{X}^{n+1}_{k}	= \mathbf{X}^{n}_{k}	+ \Delta t \, \overline{\mathbf{V}}_k^n. 
	\label{eqn.vertex.update}
\end{equation}
Once each vertex is given a unique new position $\mathbf{X}^{n+1}_{k}$, we can update all the other geometric quantities needed for the computation, e.g. normal vectors, volumes, side lengths, barycenter position, \textit{etc.} 

The procedure adopted to move the mesh will be described in the following and can be summarized in three main steps:
\begin{itemize}
	\item \textit{Lagrangian step}: a \textit{node solver} algorithm allows each node of the computational mesh to be assigned with a \textit{unique} velocity starting from the predicted solution $\q_h$ and the new node position is computed according to \eqref{eqn.vertex.update};
	\item \textit{rezoning step}: since the Lagrangian motion may lead to very distorted and stretched elements, in some cases a rezoning strategy is needed in order to achieve or recover a better mesh quality, i.e. without tangled elements;
	\item \textit{relaxation algorithm}: the aim of this step is to define the final node position by performing a linear convex combination between its Lagrangian position and its rezoned position, attempting to preserve the excellent properties in the resolution of contact waves, typically achieved by Lagrangian algorithms,  together with a good mesh quality without invalid elements. 
\end{itemize}

\subsubsection{Lagrangian step: the node solver}
\begin{figure}[htbp]
	\centering
		\includegraphics[width=0.70\textwidth]{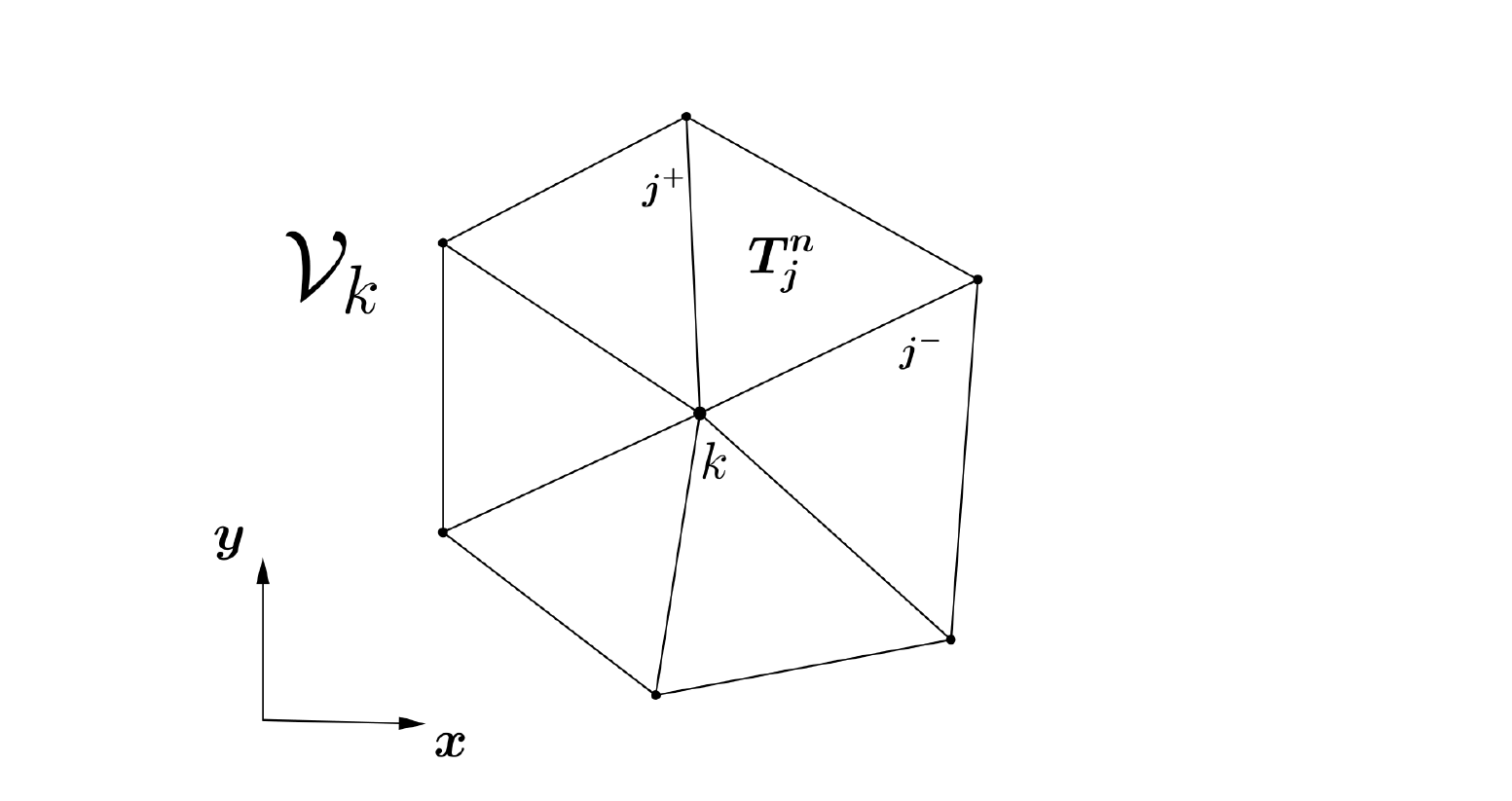}
	\caption{Geometrical notation for the node solver algorithm: $k$ is the local node, $T_j^n$ denotes one element of the neighborhood $\mathcal{V}_k$ and $(j^-,j^+)$ are the counterclockwise ordered sides of $T_j^n$ which share vertex $k$.}
	\label{fig.NS}
\end{figure}

Node solver algorithms are needed in all cell-centered Lagrangian schemes to fix a unique node velocity $\overline{\mathbf{V}}_k^n$ at each node $k$ of the mesh, starting from different contributions that come from each control volume $T_j^n \in \mathcal{V}_k$ attached to node $k$, as depicted in Figure \ref{fig.NS}. In \cite{LagNS} Boscheri et al. apply three different node solver methods to both hydrodynamics and magnetohydrodynamics and compare the numerical results among the 
various approaches. 

The simplest but most general node solver was proposed by Cheng and Shu \cite{chengshu1,chengshu3,chengshu4}, which computes the final node velocity as the arithmetic average velocity among all the contributions coming from the neighbor elements of node $k$. As done in \cite{LagNS}, this method may be improved by inserting a weighted average that depends on the geometry of the neighbors, i.e.
\begin{equation}
\overline{\mathbf{V}}_k^n = \frac{1}{\mu_k}\sum \limits_{T_j^n \in \mathcal{V}_k}{\mu_{k,j}\mathbf{V}_{k,j}^n}, 
\label{eqnNScs}
\end{equation}
where the weights $\mu_{k,j}=\rho_j |T_j^n|$ are the masses of the elements given in terms of density $\rho_j$ and volume $|T_j^n|$, while $\mu_k$ denotes the sum of all weights and the $\mathbf{V}_{k,j}^n$ are the time-averaged vertex-extrapolated velocities from element $T_j^n$ at vertex $k$ .

In \cite{Maire2010,Maire2011} Maire proposed a more sophisticated node solver algorithm for compressible hydrodynamics, which is based on the conservation of total energy. Forces on node $k$ that depend on pressure and velocity are computed for each neighbor element $T_j^n$ solving multiple approximate half-Riemann problems around a vertex on a series of sides $(j^+,{j^-})$ with the use of the acoustic Riemann solver \cite{DukowiczRS}. Finally, the node velocity is obtained as the solution of a linear algebraic equation system. 

\begin{figure}[htbp]
	\centering
		\includegraphics[width=0.50\textwidth]{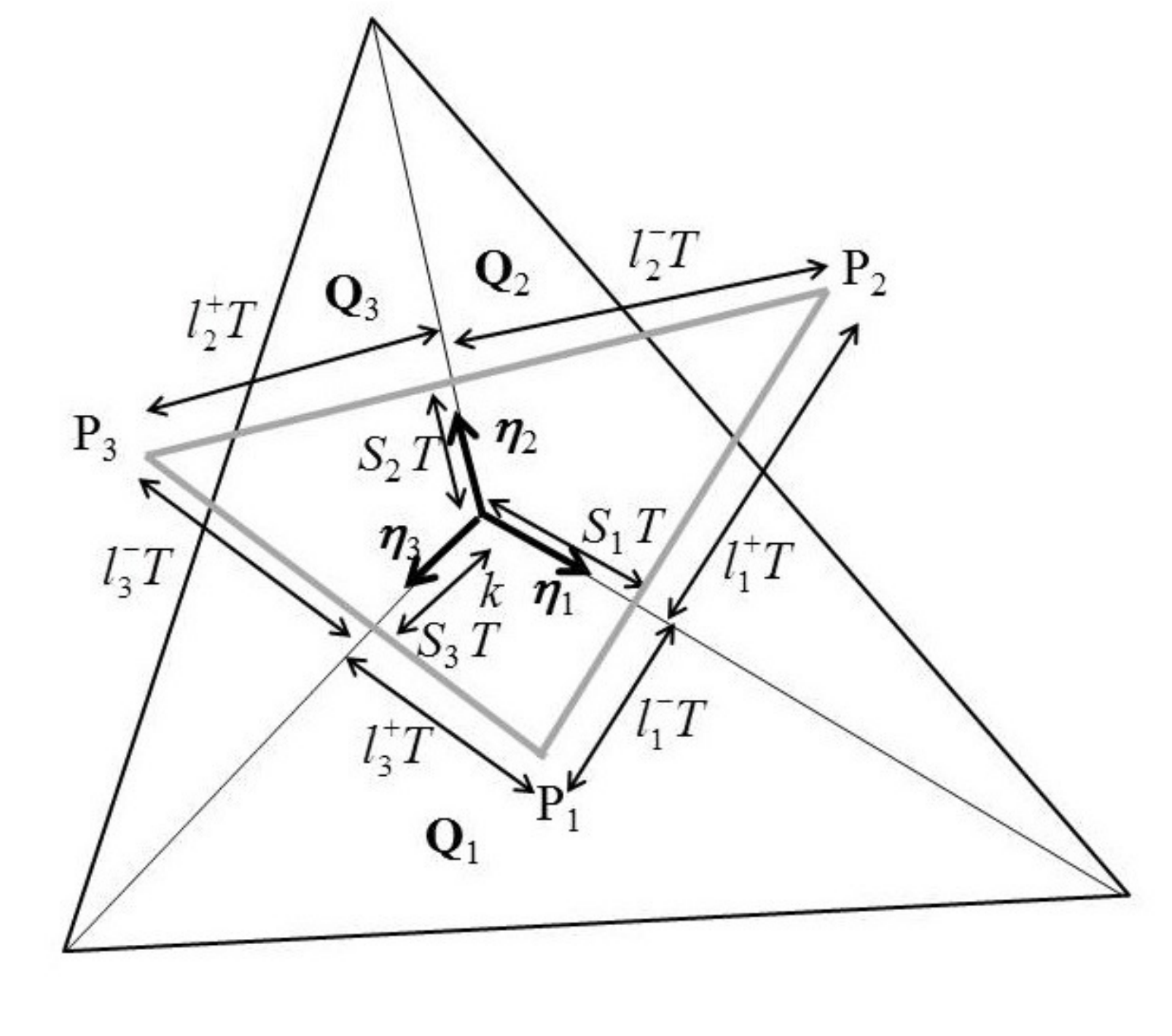}
	\caption{From \cite{BalsaraMultiDRS}: multidimensional problem at vertex $k$, where three different states $(\Q_1,\Q_2,\Q_3)$ come together. The gray lines highlights the control volume generated by the propagation of the wavespeeds $(\S_1,\S_2,\S_3)$ within a time step $\Delta t$.}
	\label{fig.NSb1}
\end{figure}

\begin{figure}[htbp]
	\centering
		\includegraphics[width=0.80\textwidth]{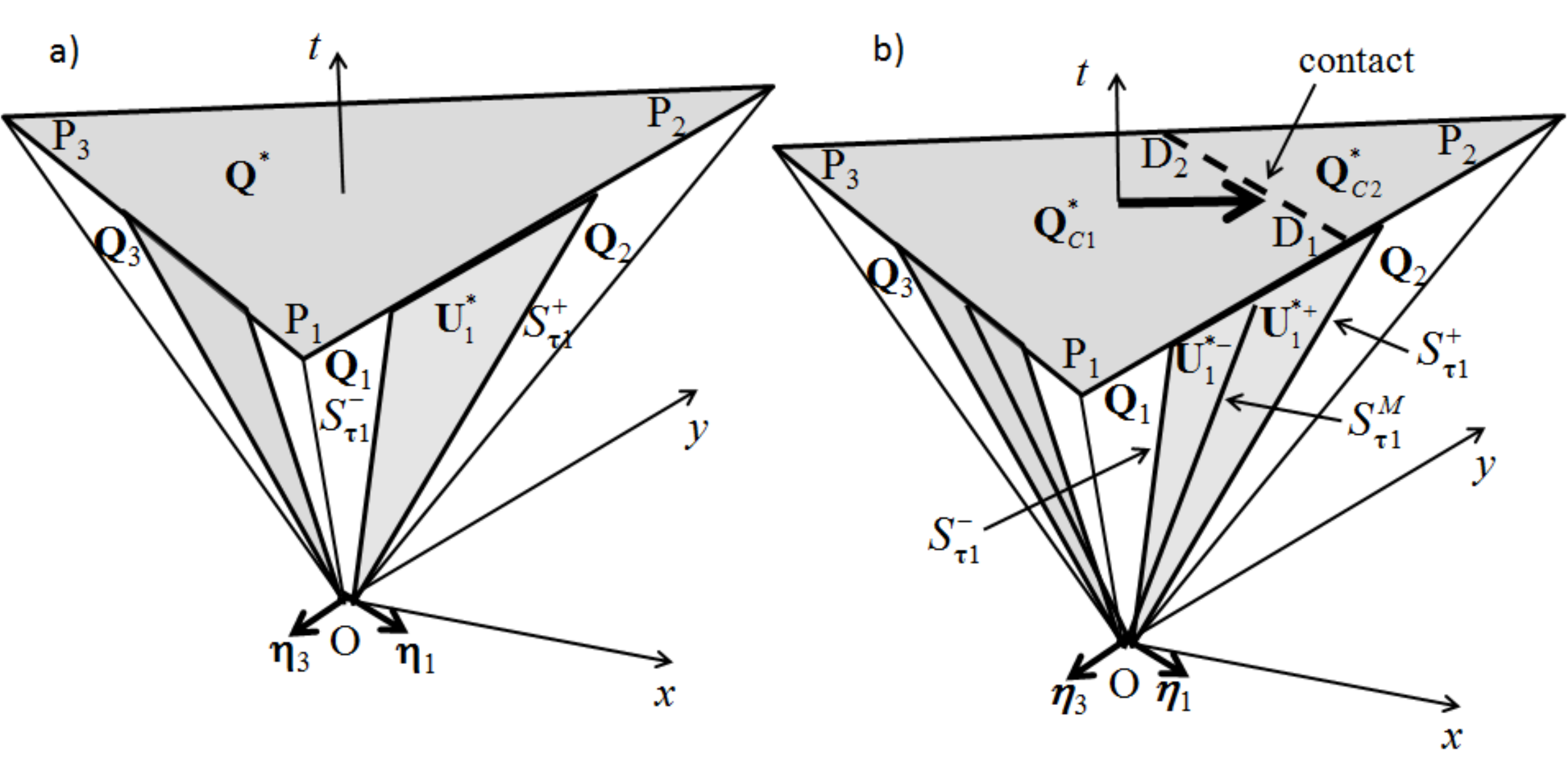}
	\caption{Figures 3a and 3b from Balsara, Dumbser and Abgrall \cite{BalsaraMultiDRS} show the space-time diagram when three states $\Q_1$, $\Q_2$ and $\Q_3$ come together at a node. Two space and one time dimensions are shown. Fig. 3a is useful for node motion. The strongly interacting state $\Q^*$ occupies a self-similar region in space-time that looks like an inverted triangular pyramid (because we have 3 incoming states). The side panels of Fig. 3a depict the one-dimensional HLL Riemann problems. Fig. 3b is useful for the corrector step. The contact discontinuity $D_1D_2$  in Fig. 3b splits the HLL state $\Q^*$ from Fig. 3a into two HLLC states $\Q^*_{C1}$ and $\Q^*_{C2}$. The side panels of Fig. 3b depict the one-dimensional HLLC Riemann problems. }
	\label{fig.NSb2}
\end{figure}

In this paper we will use the very recent approach introduced in \cite{LagNS}, where the final node velocity vector is extracted from the  multidimensional state $\Q^*$ that has been obtained as the strongly interacting state of a genuinely multi-dimensional Riemann solver. In order to obtain the  strongly interacting state $\Q^*$ we rely on the genuinely multidimensional formulation of the HLL Riemann solver for hyperbolic conservation laws on unstructured meshes proposed by Balsara et al. in \cite{BalsaraMultiDRS}. There, instead of the classical edge-based one-dimensional Riemann solvers in normal direction typically used in unsplit Godunov-type finite volume methods, the authors adopt a node-based numerical flux that takes into account the multidimensional nature of the physical flow structure. Both HLL and HLLC Riemann solvers have been considered in \cite{BalsaraMultiDRS}. We use the same notation adopted by the authors in \cite{BalsaraMultiDRS}, depicted in  Figure \ref{fig.NSb1}, where three different states $(\Q_1,\Q_2,\Q_3)$ are coming together at a vertex $k$. Let $\Q_j$ be the generic state of the neighbor element $T_j^n$ and let $\pmb{\eta}_j$ be the unit outward-pointing edge vector that separates the counterclockwise ordered states $\Q_j$ and $\Q_{j+1}$. Together with $\pmb{\eta}_j$ we define the orthogonal unit vector $\pmb{\tau}_j$ in such a way that the normal vectors $(\pmb{\eta}_j,\pmb{\tau}_j)$ form a local edge-aligned reference system. The waves propagate towards the outside of the Voronoi neighborhood along the edge direction $\pmb{\eta}_j$ with speeds $\mathbf{S}_j$ and after one timestep $\mathbf{T}=\Delta t = t^{n+1}-t^n$ they are located in the polygon bounded by vertices $P_j$. This irregular polygonal surface is denoted by $\Omega_{HLL}$ and is uniquely defined by the intersection between the lines orthogonal to $\pmb{\eta}_j$ and located at a distance $d_j=\mathbf{S}_j \mathbf{T}$ from vertex $k$ along direction $\pmb{\eta}_j$. Figure \ref{fig.NSb2} shows the time evolution of the polygonal area $\Omega_{HLL}$ that circumscribes the strongly interacting state and becomes a prism in the space-time reference system.

In order to evaluate the multidimensional state $\Q^*$, one has to solve first the one-dimensional Riemann problems perpendicular to each edge $j$ defined by $\pmb{\eta}_j$, i.e. along the $\pmb{\tau}_j$ directions, so that the resolved one-dimensional states $\Q_j^*$ are known, as represented by the darkly shaded areas on the side panels of Figure \ref{fig.NSb2}. Next, from $\Q_j^*$ we compute the wave speeds $\mathbf{S}_j$ which propagate along the edge direction within one timestep $\Delta t$ and using the \textit{multidimensional wave model} we are able to define the multidimensional area $\Omega_{HLL}$. Finally, the strongly interacting state $\Q^*$ is computed by integrating the conservation law \eqref{PDE} over the three-dimensional space-time control volume, as shown in Figure \ref{fig.NSb2}. The details for the computation of the multidimensional state $\Q^*$ can be found in \cite{BalsaraMultiDRS}, where an explicit formula for getting the multidimensional  state $\Q^*$ has been derived. The final value of the velocity vector for node $k$ is then easily  extracted from the multidimensional state $\Q^*$. Since $\overline{\mathbf{V}}_k^n$ is a time-averaged velocity, we use a standard Gauss-Legendre quadrature formula in time, hence the above procedure needs to be done for each temporal Gaussian quadrature point. 

\subsubsection{Rezoning step}
Lagrangian schemes have been developed so that the mesh follows the fluid motion as far as possible, hence allowing material interfaces or contact waves to be precisely located throughout the whole computation. When the flow motion becomes very complex, involving strong shock waves or vortex motion, the Lagrangian mesh quality drastically decreases producing highly distorted and twisted elements, that leads sometimes to an invalid computational grid, i.e. with some control volumes with negative Jacobians. Therefore a lot of effort has been made to develop proper and robust rezoning algorithms in order to maintain or to recover mesh quality during the simulation. Any Lagrangian scheme should be able to maintain the excellent resolution of the waves together with a reasonably well shaped computational mesh and this task may become very challenging in some cases. In the following we describe a suitable rezoning algorithm first presented in \cite{KnuppRezoning}, that is very efficient due to its node-centered formulation which can be carried out \textit{locally} considering each vertex and its surrounding neighbor elements.

Each rezoning algorithm starts with computing the Lagrangian coordinate vector $\mathbf{x}_k^{n+1,Lag}$ of each node $k$ of the mesh, given by \eqref{eqn.vertex.update}. For the sake of simplicity the generic element $T_j^{n+1}$ of the neighborhood $\mathcal{V}_k$ will be addressed with $j$ and let $\mathbf{x}_{j,l}=(x_{j,l},y_{j,l},z_{j,l})$ be the three counterclockwise ordered nodes $l=1,2,3$ associated with the neighbor triangle $T_j^{n+1}$. They are reordered keeping the counterclockwise convention in such a way that node $k$ corresponds to local node $l=1$ for element $T_j^{n+1}$. The rezoning algorithm is based on the optimization of a local objective function $\mathcal{K}_k$ that is defined for each node and expressed in terms of the Jacobian matrix $\mathbf{J}_{j}$ of the mapping \eqref{xietaTransf} of each neighbor element $T_j^{n+1}$, which reads 
\begin{equation}
\mathbf{J}_{j} = \left( \begin{array}{cc} x_{j,2}-x_k & y_{j,2}-y_k \\ x_{j,3}-x_k & y_{j,3}-y_k \end{array} \right).
\label{eq:locJac}
\end{equation}
$\kappa(\mathbf{J}_{j})$ represents the condition number of $\mathbf{J}_{j}$ and the objective function is evaluated considering all the elements surrounding node $k$, as done in \cite{KnuppRezoning}:
\begin{equation}
\mathcal{K}_k = \sum\limits_{T_j^{n+1} \in \mathcal{V}_k}{\kappa(\mathbf{J}_{j})}.
\end{equation}
The minimization of the above-defined function yields the optimal location of the free vertex $k$. The optimization procedure is simply chosen to be the first step of the Newton method, as proposed in \cite{MaireRezoning}, hence we need to compute the Hessian $\mathbf{H}_k$ and the gradient $\nabla \mathcal{K}_k$ of the function $\mathcal{K}_k$:
\begin{equation}
\mathbf{H}_k = \sum\limits_{T_j^{n+1} \in \mathcal{V}_k}{\left( \begin{array}{cc} 
\frac{\partial^2 \kappa(\mathbf{J}_{j})}{\partial x^2}          & \frac{\partial^2 \kappa(\mathbf{J}_{j})}{\partial x \partial y} \\ 
\frac{\partial^2 \kappa(\mathbf{J}_{j})}{\partial y \partial x} & \frac{\partial^2 \kappa(\mathbf{J}_{j})}{\partial y^2} \end{array} \right)}, 
\quad 
\nabla \mathcal{K}_k = \sum\limits_{T_j^{n+1} \in \mathcal{V}_k} \nabla \kappa(\mathbf{J}_{j}).
\label{eqn.HessGrad}
\end{equation}
The rezoned coordinate vector $\mathbf{x}_k^{Rez}$ is then given by one Newton step as follows: 
\begin{equation}
 \mathbf{x}_k^{Rez} = \mathbf{x}_k^{n+1,Lag} - \mathbf{H}_k^{-1}\left(\mathcal{K}_k\right) \cdot \nabla \mathcal{K}_k.
\label{eqn.vertex.rez}
\end{equation}

\subsubsection{Relaxation algorithm}
The rezoned coordinates may lead to an excessive change of the mesh configuration between the current time level $t^n$ and the next one $t^{n+1}$, hence causing a  loss of the excellent resolution capabilities of the Lagrangian framework. On the other side taking the pure Lagrangian coordinates could yield a very bad quality mesh in which even tangled elements might occur. Therefore the final node positions will be determined by a linear convex combination between the Lagrangian $\mathbf{x}_k^{Lag}$ and the rezoned $\mathbf{x}_k^{Rez}$ coordinates, i.e.
\begin{equation}
\mathbf{x}_k^{n+1} = \mathbf{x}_k^{n+1,Lag} + \omega_k \left( \mathbf{x}_k^{Rez} - \mathbf{x}_k^{n+1,Lag} \right),
\label{eqn.relaxation}
\end{equation}
where $\omega_k$ is a coefficient bounded in the interval $[0,1]$. The computation of $\omega_k$ is based on the Lagrangian grid deformation within a timestep and it results in such a way that for rigid body motion, i.e. pure translation and pure rotation, one has $\omega_k=0$, so that the fully Lagrangian motion of the mesh is guaranteed. If some elements are highly compressed or twisted, then the coefficient $\omega_k$ will be closer to its upper limit value of $1$. All the details of the computation of $\omega_k$ can be found in \cite{MaireRezoning}.

Since we want our scheme to be as Lagrangian as possible, we apply the rezoning and the relaxation algorithm only if really needed to carry out the computation. 
In Section \ref{sec.validation} we write explicitly whether the rezoning strategy has been used, or not.

\subsection{Finite Volume Scheme}
\label{sec.SolAlg}
In order to evolve the cell-averaged vector of conserved variables $\Q_i^n$ to the new time level $t^{n+1}$, according to \cite{BoscheriDumbserLag} we adopt a compact space-time divergence formulation of the governing PDE \eqref{PDE}, which reads
\begin{equation}
\tilde \nabla \cdot \tilde{\F} = \mathbf{S}(\Q), 
\label{PDEdiv3D}
\end{equation} 
where the space-time nabla operator is given by
\begin{equation}
\tilde \nabla  = \left( \frac{\partial}{\partial x}, \, \frac{\partial}{\partial y}, \, \frac{\partial}{\partial t} \right)^T,  \qquad 
\tilde{\F}  = \left( \mathbf{F}, \, \Q \right) = \left( \mathbf{f}, \, \mathbf{g}, \, \Q \right).
\end{equation}

\begin{figure}[htbp]
	\centering
		\includegraphics[width=0.8\textwidth]{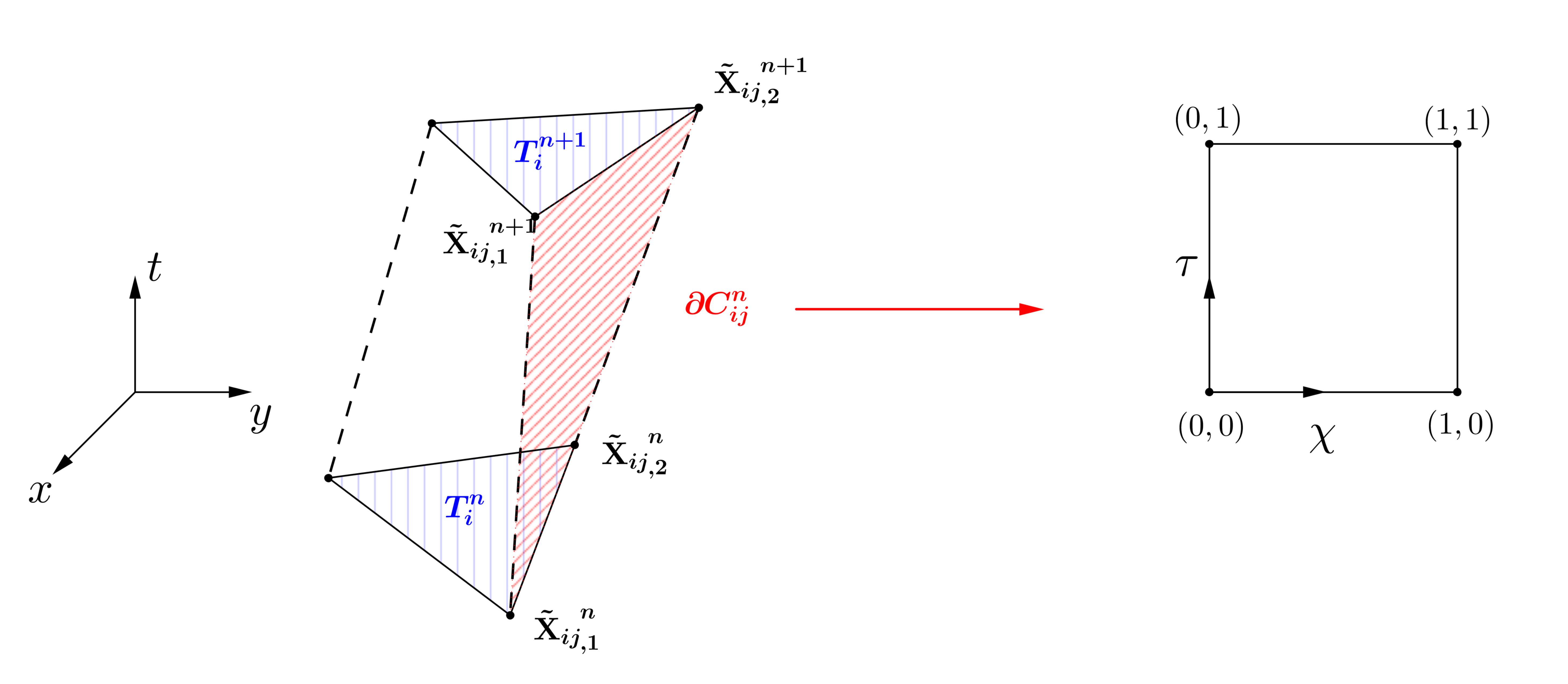}
	\caption{Physical space-time control volume $C^n_i$ and reference system $(\chi,\tau)$ adopted for the bilinear parametrization of the lateral sub-surfaces $\partial C^n_{ij}$.}
	\label{fig.FVsurf}
\end{figure}

The ALE framework involves a moving space-time control volume for each element $T_i^n$, that is obtained by connecting with \textit{straight} lines the vertices of element $T_i$ at time level $t^n$ with those at the new time level $t^{n+1}$, which are known from the predictor step and the node solver and rezoning algorithm, as fully  explained in Sections \ref{sec.localCG} and \ref{sec.meshmotion}, respectively. Thus, Eqn. \eqref{PDEdiv3D} is integrated over the space-time control volume $C^n_i = T_i(t) \times \left[t^{n}; t^{n+1}\right]$ shown in Figure \ref{fig.FVsurf}, yielding 
\begin{equation}
\int\limits_{t^{n}}^{t^{n+1}} \int \limits_{T_i(t)} \tilde \nabla \cdot \tilde{\F} \, d\mathbf{x} dt = \int\limits_{t^{n}}^{t^{n+1}} \int \limits_{T_i(t)} \S \, d\mathbf{x} dt,   
\label{STPDE}
\end{equation} 
that simplifies to
\begin{equation}
\int \limits_{\partial C^n_i} \tilde{\F} \cdot \ \mathbf{\tilde n} \, \, dS = 
\int\limits_{t^{n}}^{t^{n+1}} \int \limits_{T_i(t)} \S \, d\mathbf{x} dt,   
\label{I1}
\end{equation}   
where the left space-time volume integral has been rewritten using Gauss' theorem as the sum of the flux integrals computed over the space-time surface $\partial C^n_i$. Here, the outward pointing space-time unit normal vector is denoted by $\mathbf{\tilde n} = (\tilde n_x,\tilde n_y,\tilde n_t)$ and it is defined on the space-time surface $\partial C^n_i$. 

The surface $\partial C^n_i$ is bounded in time between the triangle at the current time level $T_i^{n}$ and the triangle at the new time level $T_i^{n+1}$. It is then closed laterally by a total number $\mathcal{N}_i$ of lateral sub-surfaces $\partial C^n_{ij} = \partial T_{ij}(t) \times [t^n;t^{n+1}]$, with $\mathcal{N}_i=3$ equal to the number of direct neighbors $T_j$ of element $T_i$, i.e. the so-called \textit{Neumann neighborhood} of $T_i$. Each of the three lateral sub-surfaces    is first mapped onto a side-aligned local reference system $(\chi,\tau)$ and then parametrized using a set of bilinear basis functions $\beta_k(\chi,\tau)$ \cite{BoscheriDumbserLag}, as depicted in Figure \ref{fig.FVsurf}. Furthermore the space-time unit normal vector $\mathbf{\tilde n}$ is computed from the parametrization of each sub-surface, see \cite{BoscheriDumbserLag,DumbserBoscheriLagNC,LagNS} for details.

Let $|T_i^{n}|$ denote the surface of triangle $T_i$ at the current time level $t^n$ and let $|\partial C_{ij}^n|$ be the determinant of the coordinate transformation of each sub-surface $\partial C^n_{ij}$. Discretization of Eqn. \eqref{I1} yields the following two-dimensional ALE one-step finite volume scheme on moving meshes:
\begin{equation}
|T_i^{n+1}| \, \Q_i^{n+1} = |T_i^n| \, \Q_i^n - \sum \limits_{T_j \in \mathcal{N}_i} \,\, {\int \limits_0^1 
\int \limits_0^1 | \partial C_{ij}^n| \tilde{\mathbf{G}}_{ij} \, d\tau d\chi}
+ \int\limits_{t^{n}}^{t^{n+1}} \int \limits_{T_i(t)} \S(\mathbf{q}_h) \, d\mathbf{x} dt, 
\label{PDEfinal}
\end{equation}
where $\tilde{\mathbf{G}}_{ij}$ denotes the numerical flux $\tilde{\mathbf{G}}_{ij} = \tilde{\F}_{ij} \cdot \mathbf{\tilde n}_{ij}$ used to resolve the discontinuity of the predictor solution $\mathbf{q}_h$ at the space-time  sub-face $\partial C_{ij}^n$. 

\begin{figure}[htbp]
	\centering
		\includegraphics[width=0.4\textwidth]{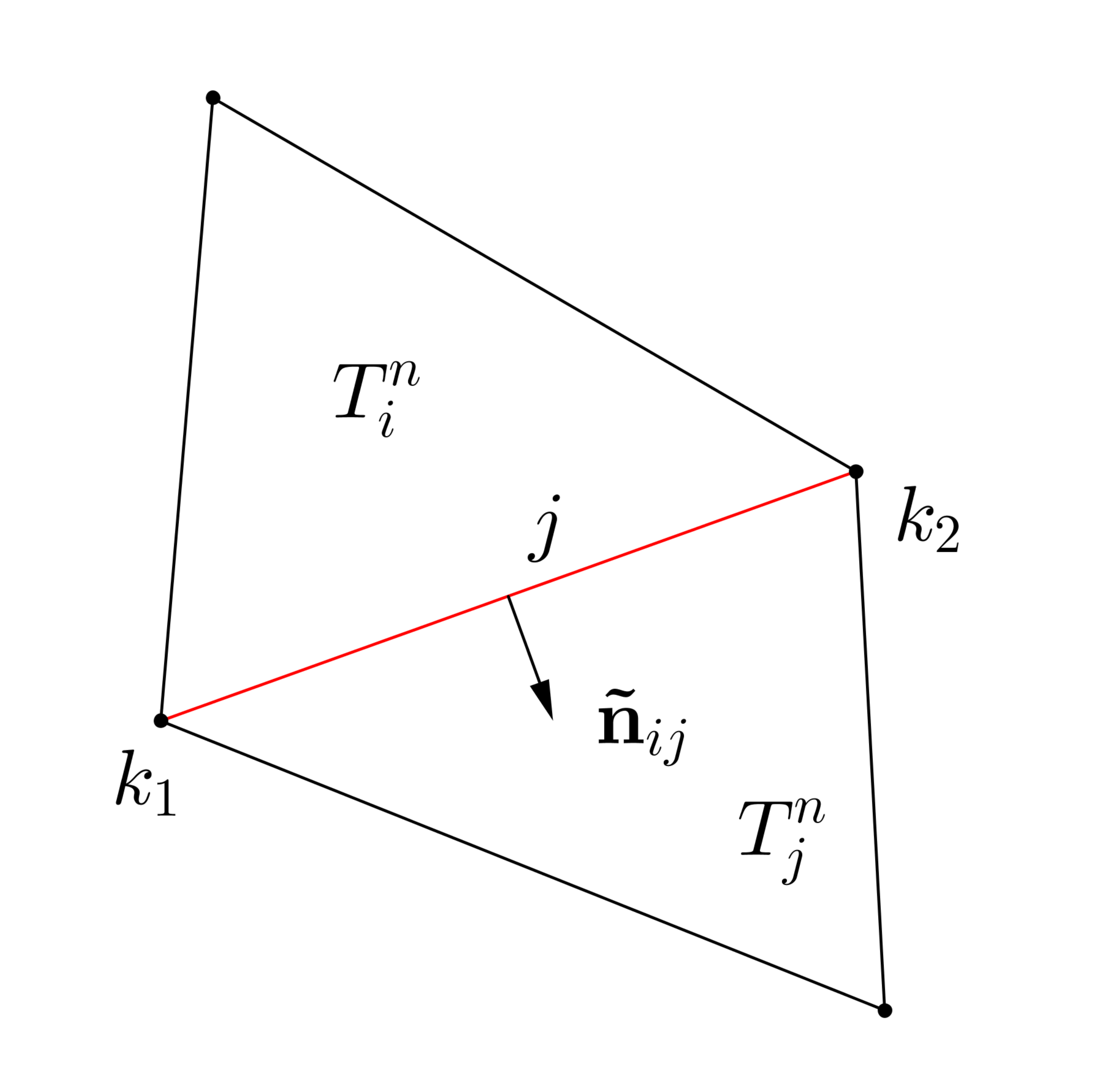}
	\caption{Notation used for the finite volume scheme. Element $T_i^n$ and its direct neighbor $T_j^n$ share edge $j$, which is bounded by vertices $k_1$ and $k_2$.}
	\label{fig.ALEflux}
\end{figure}

The numerical flux is evaluated at each sub-face by taking into account a multidimensional vertex-based flux and a one-dimensional edge-based flux. According to Figure \ref{fig.ALEflux}, let $(k_1,k_2)$ be the two vertices that bound edge $j$ and let $(\q_h^-,\q_h^+)$ be the numerical solution inside element $T_i(t)$ and inside the neighbor element $T_j(t)$, respectively. The term $\tilde{\mathbf{G}}_{ij}$ is computed as follows:
\begin{itemize}
	\item first we solve the Riemann problem around the two vertices $(k_1,k_2)$ of face $\partial C_{ij}$ using the multidimensional HLL Riemann solver \cite{BalsaraMultiDRS}, hence obtaining the multidimensional state $(\Q^*_1,\Q^*_2)$ and the multidimensional numerical fluxes $(\tilde{\F}^*_1,\tilde{\F}^*_2)$. The  multidimensional HLL formulation adopted here is the same algorithm used as node solver and explained in Section \ref{sec.meshmotion}. Now we do not limit to evaluate the interacting state $\Q^*$, but we also compute the multidimensional fluxes $\tilde{\F}^*$ for each vertex $(k_1,k_2)$ of edge $j$. Positivity of density and pressure is guaranteed by using the self-adjusting positivity preserving scheme of \cite{BalsaraPositivity2012}, extended to unstructured meshes. According to \cite{BalsaraMultiDRS}, the final expression for the multidimensional vertex--based fluxes $\tilde{\F}^*$ are computed by a blending between the multidimensional HLL and HLLC fluxes, with the blending factor given by the flattener variable introduced and presented in details in \cite{BalsaraPositivity2012};  
	\item then a classical Godunov-type one-dimensional edge flux $\tilde{\F}_{edge}$ has to be determined, that is projected orthogonally w.r.t. the edge $j$, as usually done on unstructured meshes. The one-dimensional ALE-type HLL flux can be formulated as 
\begin{equation}
  \tilde{\F}_{edge} \cdot \mathbf{\tilde n}_{ij} =  
  \frac{1}{s_R - s_L} \left[ \left( s_R \tilde{\F}(\q_h^-) - s_L \tilde{\F}(\q_h^+)  \right) \cdot \mathbf{\tilde n}_{ij}  + s_L s_R \left( \q_h^+ - \q_h^- \right) \right],  
  \label{eqn.rusanov} 
\end{equation} 
where $s_L$ and $s_R$ are the usual HLL estimates of the left and right signal speeds, associated with the ALE Jacobian matrix in spatial normal direction, given by 
\begin{equation} 
\mathbf{A}^{\!\! \mathbf{V}}_{\mathbf{n}}(\Q)=\left(\sqrt{\tilde n_x^2 + \tilde n_y^2}\right)\left[\frac{\partial \mathbf{F}}{\partial \Q} \cdot \mathbf{n}  - 
(\mathbf{V} \cdot \mathbf{n}) \,  \mathbf{I}\right], \qquad 
\textnormal{ with } \qquad   
\mathbf{n} = \frac{(\tilde n_x, \tilde n_y)^T}{\sqrt{\tilde n_x^2 + \tilde n_y^2}}.  
\end{equation} 
The local normal mesh velocity is denoted by $\mathbf{V} \cdot \mathbf{n}$ and $\mathbf{I}$ represents the identity matrix. 


\item the spatial part of the space-time surface integral at the space-time sub-face $\partial C_{ij}^n$ is  computed using the Simpson rule, which achieves up to fourth order of accuracy. In time, classical Gauss-Legendre  quadrature with two quadrature points is used. The final approximation of the lateral space-time surface integrals reads  
\begin{eqnarray}
\int \limits_0^1 \int \limits_0^1 | \partial C_{ij}^n| \tilde{\mathbf{G}}_{ij} d\tau d\chi & \approx & 
\sum_j \omega_j \left( \frac{1}{6} |\partial C_{ij}^n|(0,\tau_j) \tilde{\F}^*_1(\tau_j) \cdot \mathbf{\tilde n}_{ij}(0,\tau_j)  \right. 
\nonumber \\
&& \hspace{10mm} + \frac{4}{6} | \partial C_{ij}^n|(\frac{1}{2},\tau_j) \tilde{\F}_{edge}(\tau_j) \cdot \mathbf{\tilde n}_{ij}(\frac{1}{2},\tau_j)   \nonumber \\ 
&& \hspace{10mm} \left. + \frac{1}{6} | \partial C_{ij}^n|(1,\tau_j) \tilde{\F}^*_2(\tau_j) \cdot \mathbf{\tilde n}_{ij}(1,\tau_j)  \right),
\label{eqn.SimpsonFlux}
\end{eqnarray}
\end{itemize}
where $\tau_j$ and $\omega_j$ are the temporal quadrature points and weights, respectively. 

The timestep $\Delta t$ is evaluated as 
\begin{equation}
\Delta t = \textnormal{CFL} \, \min \limits_{T_i^n} \frac{D_i^n}{|\lambda^n_{\max,i}|}, \qquad \forall T_i^n \in \Omega^n, 
\label{eq:timestep}
\end{equation}
where $D_i^n$ is the incircle diameter of element $T_i^n$ and $|\lambda^n_{\max,i}|$ is the maximum absolute value of the eigenvalues computed from the solution 
$\Q_i^n$ in $T_i^n$. For unsplit Godunov-type schemes in in two space dimensions based on \textit{one-dimensional Riemann solvers} the Courant number CFL must satisfy $\textnormal{CFL} < 0.5$ for linear stability, as mentioned in \cite{ToroBook}. However, numerical evidence indicates that our finite volume schemes based on \textit{multidimensional Riemann solvers} are able to run in a stable manner also with a much less restrictive CFL condition of $\textnormal{CFL}<1$, 
because of the multidimensionality introduced in the numerical flux evaluation, see \cite{BalsaraMultiDRS}. 
Hence, for most of the test problems presented in the  next Section \ref{sec.validation} the CFL number has been actually set very close to this experimentally observed limit by choosing $\textnormal{CFL}=0.95$. Thus, the multidimensional finite volume scheme can run the same test case much more efficiently than a classical edge-based finite volume algorithm. This leads to a significant improvement in terms of computational efforts, especially 
in the Lagrangian framework, which is typically characterized by very small timesteps caused by strongly deformed and distorted elements. Furthermore, the high order WENO reconstruction on moving unstructured meshes is very expensive, since the reconstruction equations can no longer be solved once and for all in a preprocessor stage, as it was the case for Eulerian schemes on fixed meshes in \cite{DumbserKaeser06b,Dumbser2007204}. Hence, the possibility to use larger time 
steps leads to less reconstructions to be done when running a simulation to a given final time reducing thus the total computational effort. 
In the next section we will assess the possible gains by using multidimensional Riemann solvers also quantitatively in terms of achievable accuracy as a 
function of CPU time. 

\section{Test problems}
\label{sec.validation} 
In this section we show some classical numerical test problems in order to validate the numerical method presented in this article. We will consider the Euler equations for compressible gas dynamics as well as the ideal classical equations of magnetohydrodynamics (MHD). For each of those hyperbolic conservation laws 
we carry out numerical convergence studies and run the scheme for some well-established benchmark problems. 

The Euler equations of compressible gas dynamics can be written in terms of the vector of conserved variables $\Q$ and the flux tensor $\F=(\f,\g)$ as 
\begin{equation}
\label{eulerTerms}
\Q = \left( \begin{array}{c} \rho \\ \rho u \\ \rho v \\ \rho E \end{array} \right), \quad \f = \left( \begin{array}{c} \rho u \\ \rho u^2 + p \\ \rho uv \\ u(\rho E + p) \end{array} \right), \quad \g = \left( \begin{array}{c} \rho v \\ \rho uv \\ \rho v^2 + p  \\ v(\rho E + p) \end{array} \right),  
\end{equation}
where $\rho$ is the mass density, $\mathbf{v}=(u,v)$ denotes the velocity vector and $\rho E$ is the total energy density, while $p$ is the fluid pressure. The  system is closed by the equation of state (EOS) of an ideal gas: 
\begin{equation}
\label{eqn.eos} 
p = (\gamma-1)\left(\rho E - \frac{1}{2} \rho (u^2+v^2) \right).  
\end{equation}

The equations of ideal magneto hydrodynamics (MHD) constitute a more complicated system that takes into account also the magnetic field $\mathbf{B}=\left(B_x,B_y\right)$, the divergence of which must remain zero in time. This adds the following involution constraint to the PDE system 
\begin{equation}
\frac{\partial B_x}{\partial x} + \frac{\partial B_y}{\partial y} = 0,
\label{eq:divB}
\end{equation}
which always holds in the continuous case if the initial data for $\mathbf{B}$ are divergence-free. When we discretize the system, the divergence constraint may be violated, hence producing non-physical perturbations in the magnetic field. A possible solution to overcome this problem has been presented in \cite{Balsara2004,BalsaraRMHD,BalsaraRumpf}, where ADER-WENO MHD schemes that were divergence free were presented. In \cite{BalsaraAMRdivfree2001} Balsara showed a full-fledged scheme for second-order accurate, divergence-free evolution of vector fields on an adaptive mesh refinement (AMR), demonstrating that high order divergence-free reconstruction can be done at all orders of accuracy. Nevertheless the schemes cited so far have not been developed for unstructured meshes, therefore in this paper we adopt the hyperbolic version of the generalized Lagrangian multiplier (GLM) divergence cleaning approach proposed by Dedner et al. \cite{Dedneretal}, where one more variable $\Psi$ as well as one linear scalar PDE are added to the ideal classical MHD system in order to carry divergence errors out of the computational domain with an artificial divergence cleaning speed $c_h$. This leads to the so-called \textit{augmented} MHD system, which reads
\begin{equation}
\Q^{T} = \left( \rho , \rho u , \rho v , \rho E , B_x , B_y \right), \nonumber
\end{equation}
\begin{equation}
\label{MHDTerms}
\f = \left( \begin{array}{c} \rho u \\ \rho u^2 + \left(p+\frac{1}{8\pi}\mathbf{B}^2 \right) - \frac{B_x B_x}{4\pi} \\ \rho uv - \frac{B_y B_x}{4\pi} \\ u\left(\rho E + p+\frac{1}{8\pi}\mathbf{B}^2\right) - \frac{B_x\left(\mathbf{v} \cdot \mathbf{B} \right)}{4\pi} \\ B_x u - u B_x + \Psi \\ B_y u - v B_x \\ c_{h}^{2} B_x \end{array} \right), \qquad \g = \left( \begin{array}{c} \rho v \\ \rho uv - \frac{B_x B_y}{4\pi} \\ \rho v^2 + \left(p+\frac{1}{8\pi}\mathbf{B}^2 \right) - \frac{B_y B_y}{4\pi} \\ v\left(\rho E + p+\frac{1}{8\pi}\mathbf{B}^2\right) - \frac{B_y\left(\mathbf{v} \cdot \mathbf{B} \right)}{4\pi} \\ B_x v - u B_y \\ B_y v - v B_y + \Psi \\ c_{h}^{2} B_y \end{array} \right),  
\end{equation}
where the state vector is represented by $\Q$ and the fluxes are $\F=(\f,\g)$. The equation of state for closing the system is
\begin{equation}
p = (\gamma-1)\left(\rho E - \frac{1}{2} (u^2+v^2) - \frac{(B_x^2+B_y^2)}{8\pi} \right).  
\label{eqn.MHDeos}
\end{equation}

In the following we apply the new Lagrangian finite volume schemes based on genuinely multidimensional HLL Riemann solvers to the above-presented hyperbolic conservation laws for some classical test problems. Each test case is chosen to be run with the local mesh velocity being equal to the local fluid velocity,
hence $\mathbf{V}=\mathbf{v}$.  

\subsection{Numerical Convergence Study for the Euler equations}
\label{sec.conv.Rates-Eul}
We first perform a numerical convergence study for the Euler equations of compressible gas dynamics considering the smooth isentropic vortex, see e.g.  
\cite{HuShuVortex1999}. Let the initial computational domain be the square $\Omega(0)=[0;10]\times[0;10]$ and let the vector of primitive variables at the 
initial time be defined as 
\begin{equation}
\label{ShuVortIC}
(\rho, u, v, p) = (1+\delta \rho, 1+\delta u, 1+\delta v, 1+\delta p),
\end{equation}  
where $\delta \rho,\delta u, \delta v, \delta p$ are some perturbations defined in the following. We define a radial coordinate as $r^2=(x-5)^2+(y-5)^2$, set the vortex strength to $\epsilon=5$ and the ratio of specific heats is chosen as $\gamma=1.4$. The perturbations for density and pressure are given by 
\begin{equation}
\label{rhopressDelta}
\delta \rho = (1+\delta T)^{\frac{1}{\gamma-1}}-1, \quad \delta p = (1+\delta T)^{\frac{\gamma}{\gamma-1}}-1, 
\end{equation} 
where $\delta T$ denotes the perturbation of temperature, $\delta T = -\frac{(\gamma-1)\epsilon^2}{8\gamma\pi^2}e^{1-r^2}$. 
The perturbation of velocity $\mathbf{v}=(u,v)$ is taken to be 
\begin{equation}
\label{ShuVortDeltaV}
\left(\begin{array}{c} \delta u \\ \delta v \end{array}\right) = \frac{\epsilon}{2\pi}e^{\frac{1-r^2}{2}} \left(\begin{array}{c} -(y-5) \\ \phantom{-}(x-5) \end{array}\right). 
\end{equation} 

We set periodic boundary conditions on each side of the domain and the final time of the simulation is $t_f=1.0$. The vortex is also convected with velocity $\v_c=(1,1)$, hence the analytical solution $\Q_e$ is the time-shifted initial condition given by $\Q_e(\x,t_f)=\Q(\x-\v_c t_f,0)$. The HLLC flux has been used to run this  test problem on a series of successive refined meshes and the corresponding error has been measured in $L_2$  norm as 
\begin{equation}
  \epsilon_{L_2} = \sqrt{ \int \limits_{\Omega(t_f)} \left( \Q_e(x,y,t_f) - \w_h(x,y,t_f) \right)^2 dxdy }, 
  \label{eqnL2error}
\end{equation} 
where $h(\Omega(t_f))$ represents the mesh size which is the maximum diameter of the circumcircles of the triangles in the domain $\Omega(t_f)$ at the final time $t_f$. Table \ref{tab.convEul} shows the convergence results up to fourth order of accuracy run with $\textnormal{CFL}=0.95$.

\begin{table}  
\caption{Numerical convergence results for the compressible Euler equations using the Lagrangian one-step WENO finite volume schemes with the genuinely multidimensional HLL Riemann solvers presented in this article. The error norms refer to the variable $\rho$ (density) at time $t=1.0$ for first up to fourth order version of the scheme.} 
\begin{center} 
\begin{small}
\renewcommand{\arraystretch}{1.0}
\begin{tabular}{cccccc} 
\hline
  $h(\Omega(t_f))$ & $\epsilon_{L_2}$ & $\mathcal{O}(L_2)$ & $h(\Omega,t_f)$ & $\epsilon_{L_2}$ & $\mathcal{O}(L_2)$ \\
\hline
  \multicolumn{3}{c}{$\mathcal{O}1$} & \multicolumn{3}{c}{$\mathcal{O}2$} \\
\hline
3.60E-01 & 2.0390E-01 & -   & 3.41E-01  & 2.8130E-02 & -     \\ 
2.45E-01 & 1.5446E-01 & 0.7 & 2.49E-01  & 1.6398E-02 & 1.7   \\ 
1.71E-01 & 1.0728E-01 & 1.0 & 1.67E-01  & 7.8020E-03 & 1.9   \\ 
1.33E-01 & 8.1766E-02 & 1.1 & 1.28E-01  & 4.1857E-03 & 2.3   \\ 
\hline 
  \multicolumn{3}{c}{$\mathcal{O}3$} & \multicolumn{3}{c}{$\mathcal{O}4$}   \\
\hline
3.29E-01 & 2.1551E-02 & -   & 3.29E-01  & 5.9233E-03 & -      \\ 
2.52E-01 & 1.0161E-02 & 2.8 & 2.51E-01  & 2.1675E-03 & 3.7   \\ 
1.67E-01 & 3.7967E-03 & 2.4 & 1.67E-01  & 4.8533E-04 & 3.7   \\ 
1.28E-01 & 1.7601E-03 & 2.8 & 1.28E-01  & 1.5816E-05 & 4.1   \\ 
\hline 
\end{tabular}
\end{small}
\end{center}
\label{tab.convEul}
\end{table}

\begin{figure}[!htbp]
\begin{center}
\begin{tabular}{cc} 
\includegraphics[width=0.47\textwidth]{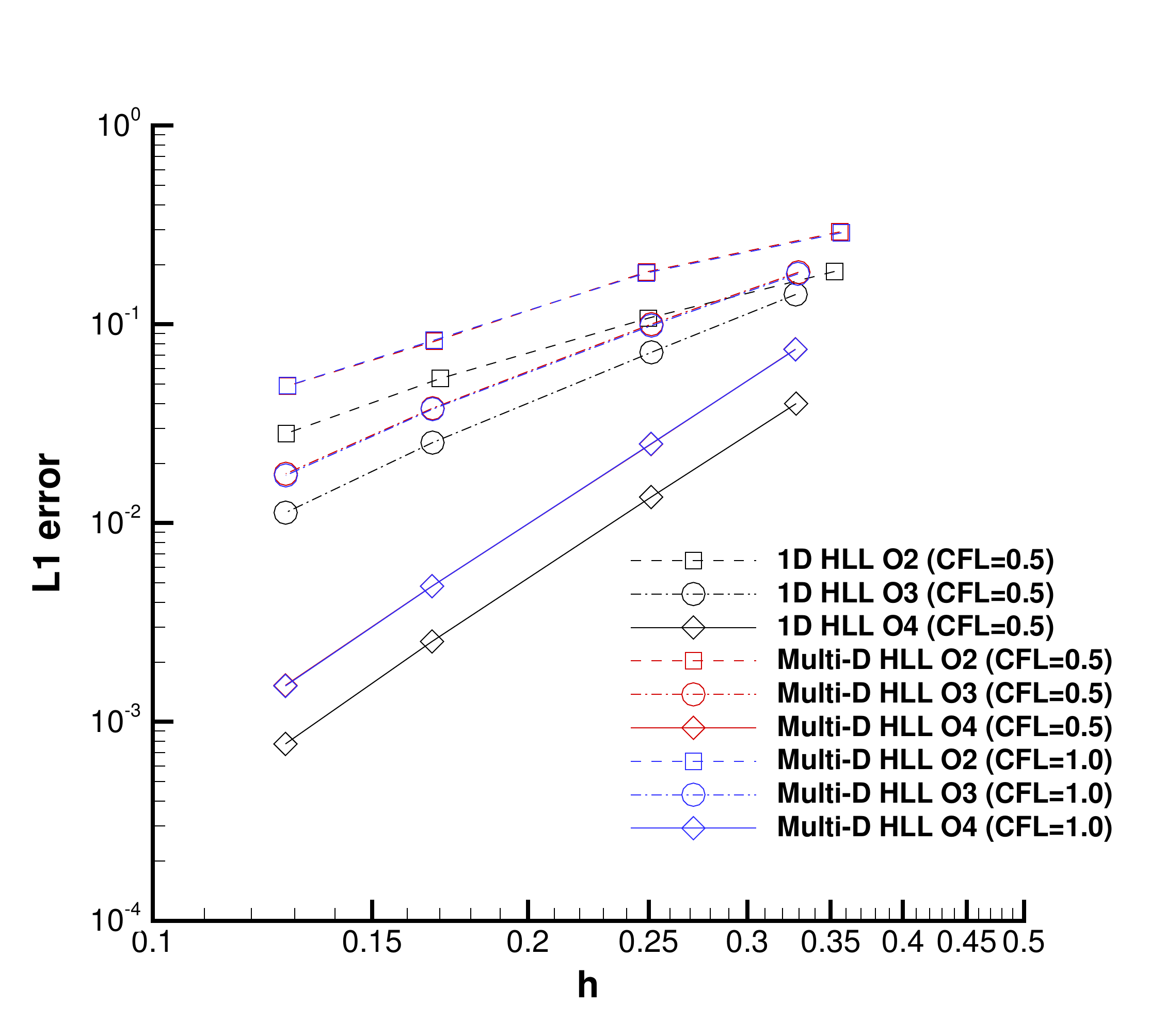}  &           
\includegraphics[width=0.47\textwidth]{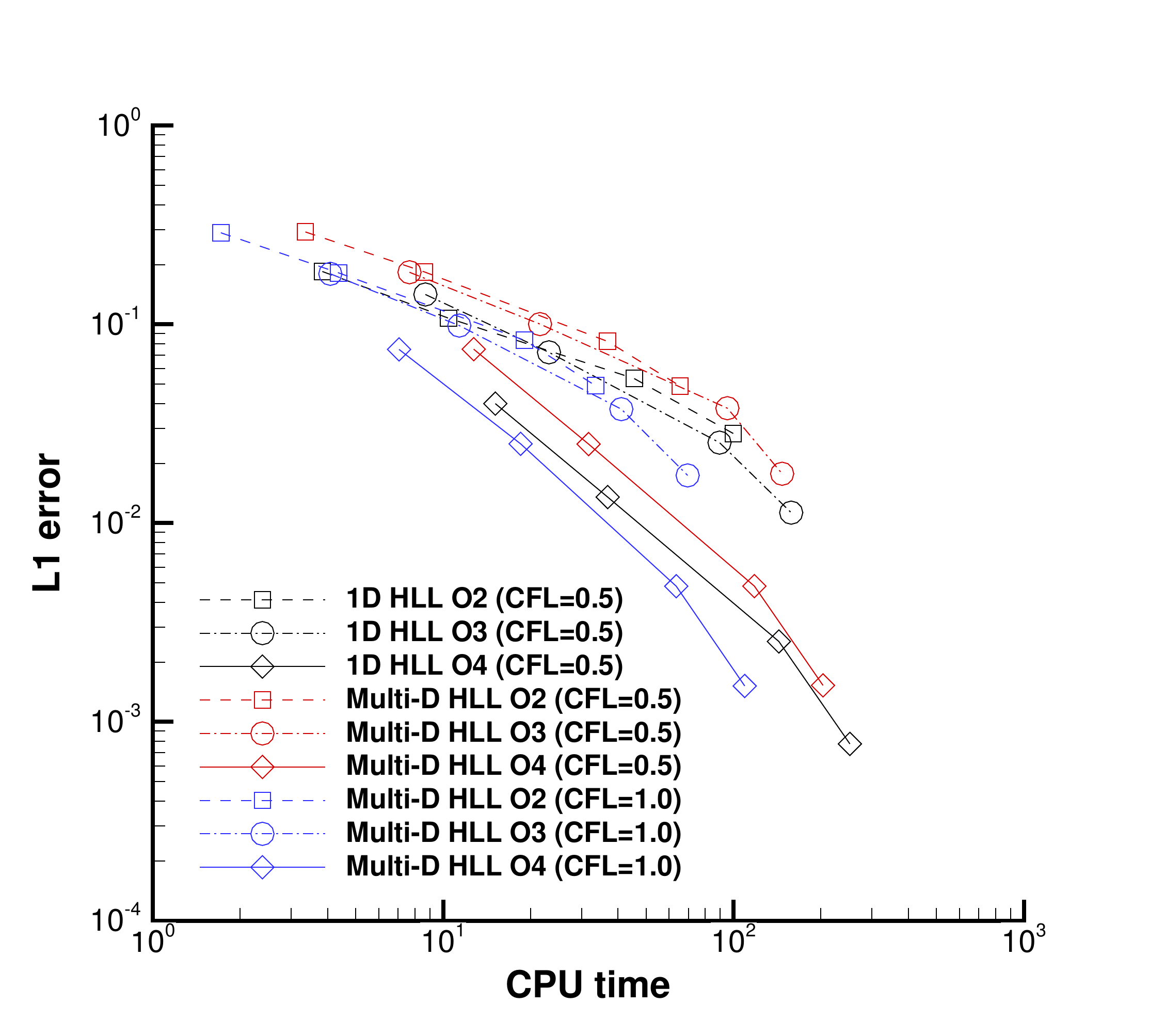} \\   
\end{tabular} 
\caption{Comparison between 1D HLL and Multi-D HLL Riemann solvers from first up to fourth order of accuracy with different CFL number. Left: dependency of the error norm on the mesh size. Right: dependency of the  error norm on the CPU time.} 
\label{fig.1DvsMultiD}
\end{center}
\end{figure}

In the following a fair comparison between the high order ADER-WENO ALE schemes based on one-dimensional Riemann solvers presented in 
\cite{BoscheriDumbserLag,LagNS} and the new ADER-WENO ALE algorithm based on multi-dimensional Riemann solvers illustrated in this paper is 
carried out. 
For this purpose we run the smooth vortex test problem again using both, a classical one-dimensional HLL Riemann solver and the genuinely 
multidimensional HLL solver, comparing in detail CPU time and accuracy. The behavior of the different solvers is depicted in Figure \ref{fig.1DvsMultiD}: 
the one-dimensional HLL Riemann solver with $\textnormal{CFL}=0.5$ is drawn by the black lines, while red and blue lines refer to the multidimensional HLL solver 
with $\textnormal{CFL}=0.5$ and $\textnormal{CFL}=1.0$, respectively. The error has been evaluated in $L_1$ norm as 
\begin{equation}
  \epsilon_{L_1} = \int \limits_{\Omega(t_f)} \left( \Q_e(x,y,t_f) - \w_h(x,y,t_f) \right) dxdy, 
\end{equation}
and the CPU time has been measured as the accumulated time obtained running the simulation in parallel on four Intel Core i7-2600 CPUs with a clock-speed of 3.40GHz. The multidimensional Riemann solver allows the scheme to be run with CFL condition of unity, hence representing clearly the most efficient algorithm 
in terms of computational efficiency (blue lines in the right panel of Fig. \ref{fig.1DvsMultiD}), although the most accurate one on a given mesh remains 
the classical one-dimensional Riemann solver (black lines in the left panel of Fig. \ref{fig.1DvsMultiD}). Numerical values of the mesh size $h$, the L1 
error norm  and the corresponding CPU time are reported in Table \ref{tab.1DvsMultiD} for each of the simulations contained in Figure \ref{fig.1DvsMultiD}.

\begin{table}  
\caption{Error norms and CPU times related to the comparison between 1D HLL and Multi-D HLL Riemann solvers from first up to fourth order of accuracy with different CFL number, shown in Figure \ref{fig.1DvsMultiD}}. 
\begin{center} 
\begin{scriptsize}
\begin{tabular}{|c||c|c|c|} 
\hline
\textbf{Numerical scheme} &  \textbf{mesh size} $h$    &    \textbf{L1 error}  & \textbf{CPU time}    \\  
\hline
1D HLL O2        &  3.5243E-01	& 1.8466E-01	& 3.8376E+00  \\
								 &	2.4976E-01	& 1.0750E-01	& 1.0436E+01  \\
								 &	1.7007E-01	& 5.3419E-02	& 4.5505E+01  \\
								 &	1.2795E-01	& 2.8263E-02	& 9.9841E+01  \\
\hline
1D HLL O3        &	3.2803E-01	& 1.4144E-01	& 8.6893E+00  \\
								 &	2.5128E-01	& 7.2359E-02	& 2.3151E+01  \\
								 &	1.6768E-01	& 2.5417E-02	& 8.9513E+01  \\
								 &	1.2783E-01	& 1.1291E-02	& 1.5808E+02  \\
\hline
1D HLL O4        &  3.2830E-01	& 3.9931E-02	& 1.5132E+01  \\
								 &	2.5113E-01	& 1.3526E-02	& 3.6879E+01  \\
								 &	1.6753E-01	& 2.5384E-03	& 1.4335E+02  \\
								 & 	1.2778E-01	& 7.7492E-04	& 2.5174E+02  \\															
\hline
\hline
MultiD HLL O2 (CFL=0.5) &	 3.5576E-01	& 2.9206E-01	& 3.3540E+00  \\
												&  2.4910E-01	& 1.8383E-01	& 8.6113E+00  \\
												&  1.6823E-01	& 8.2143E-02	& 3.6847E+01  \\
												&  1.2826E-01	& 4.8994E-02	& 6.5458E+01  \\
\hline
MultiD HLL O3 (CFL=0.5) & 3.2967E-01	& 1.8296E-01	& 7.6596E+00  \\
												&	2.5128E-01	& 1.0036E-01	& 2.1559E+01  \\
												&	1.6771E-01	& 3.7877E-02	& 9.5285E+01  \\
												&	1.2787E-01	& 1.7719E-02	& 1.4711E+02  \\
\hline
MultiD HLL O4 (CFL=0.5) &	3.2789E-01	& 7.5060E-02	& 1.2745E+01  \\
												&	2.5112E-01	& 2.4996E-02	& 3.1699E+01  \\
												&	1.6754E-01	& 4.8209E-03	& 1.1806E+02  \\
                        & 1.2777E-01	& 1.5259E-03	& 2.0344E+02  \\
\hline
\hline
MultiD HLL O2 (CFL=1.0) & 3.5699E-01	& 2.8910E-01	& 1.7160E+00  \\
												&	2.4897E-01	& 1.8169E-01	& 4.3680E+00  \\
												&	1.6815E-01	& 8.3235E-02	& 1.9063E+01  \\
												&	1.2827E-01	& 4.9115E-02	& 3.3540E+01  \\
\hline
MultiD HLL O3 (CFL=1.0) & 3.2937E-01	& 1.7966E-01	& 4.0872E+00  \\
												&	2.5128E-01	& 9.8654E-02	& 1.1372E+01  \\
												&	1.6770E-01	& 3.7463E-02	& 4.1091E+01  \\
												&	1.2787E-01	& 1.7384E-02	& 6.9545E+01  \\
\hline
MultiD HLL O4 (CFL=1.0) &	3.2794E-01	& 7.4847E-02	& 7.0512E+00  \\
												&	2.5107E-01	& 2.5067E-02	& 1.8486E+01  \\
												&	1.6754E-01	& 4.8157E-03	& 6.3539E+01  \\
												&	1.2777E-01	& 1.5162E-03	& 1.0936E+02  \\
\hline																																											
\end{tabular}
\end{scriptsize}
\end{center}
\label{tab.1DvsMultiD}
\end{table}

\subsection{Two-Dimensional Explosion Problem} 
\label{sec.EP2D}
\begin{figure}[!htbp]
\begin{center}
\begin{tabular}{cc} 
\includegraphics[width=0.47\textwidth]{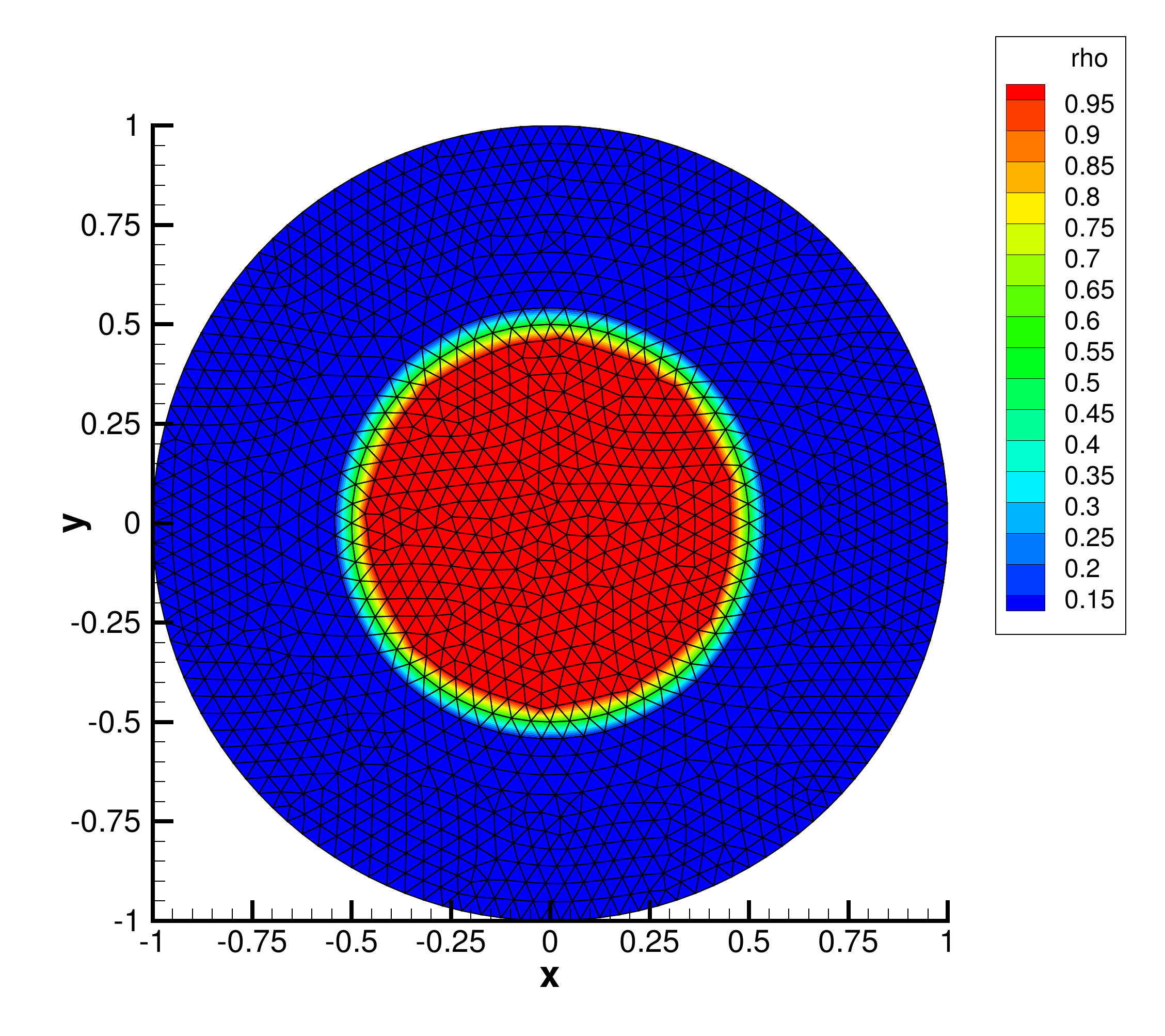}  &           
\includegraphics[width=0.47\textwidth]{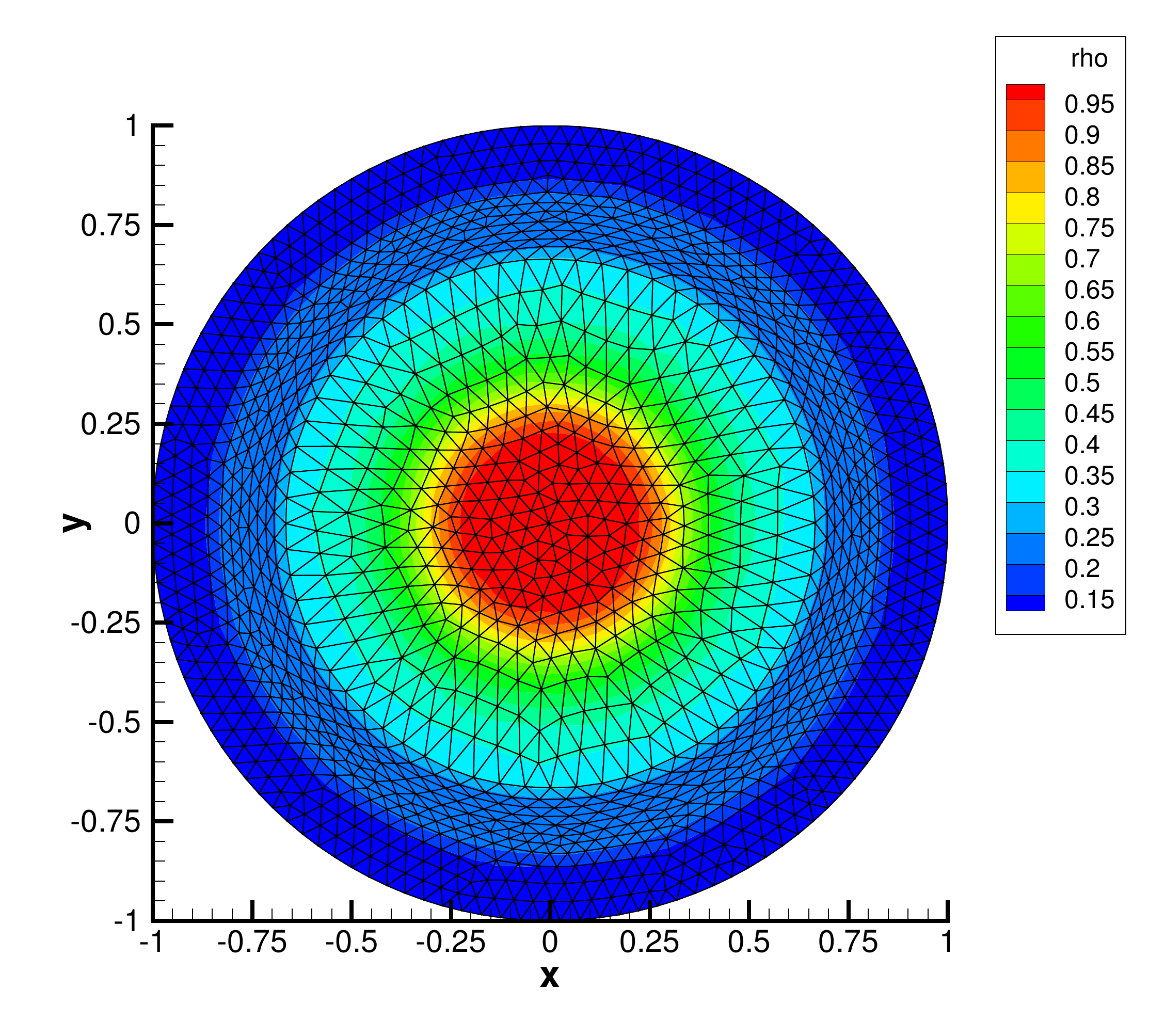} \\   
\end{tabular} 
\caption{Initial ($t=0$) and final ($t=t_f$) density distribution and mesh configuration for the two-dimensional explosion problem on a coarse grid.} 
\label{fig.EP2Dc}
\end{center}
\end{figure}

We use the third order version of our ALE finite volume scheme to solve a two-dimensional explosion problem.  
The initial domain $\Omega(0)$ is a circle of radius $R_o=1$ meshed with a total number of $N_E=68324$ elements. Let $r$ be the general radial position defined as $r=\sqrt{x^2+y^2}$. The initial condition is  given in terms of primitive variables $\U=(\rho,u,v,p)$ by two different states, separated by the circle 
of radius $R=0.5$: 
\begin{equation}
  \U(\x,0) = \left\{ \begin{array}{clcc} \U_i = & (1.0, 0.0, 0.0, 1.0)    & \textnormal{ if } & r \leq R, \\ 
                                         \U_o = & (0.125, 0.0, 0.0, 0.1)  & \textnormal{ if } & r > R,        
                      \end{array}  \right. 
\end{equation}
where $\U_i$ represents the \textit{inner state} and $\U_o$ the \textit{outer state} in primitive variables. The ratio of specific heats is assumed to be $\gamma = 1.4$ and the final simulation time is $t_f=0.25$. In Figure \ref{fig.EP2Dc} a coarse grid has been used to show the initial and the final density distribution and the corresponding mesh deformation for the two-dimensional explosion problem.

\begin{figure}[!htbp]
\begin{center}
\begin{tabular}{cc} 
\includegraphics[width=0.47\textwidth]{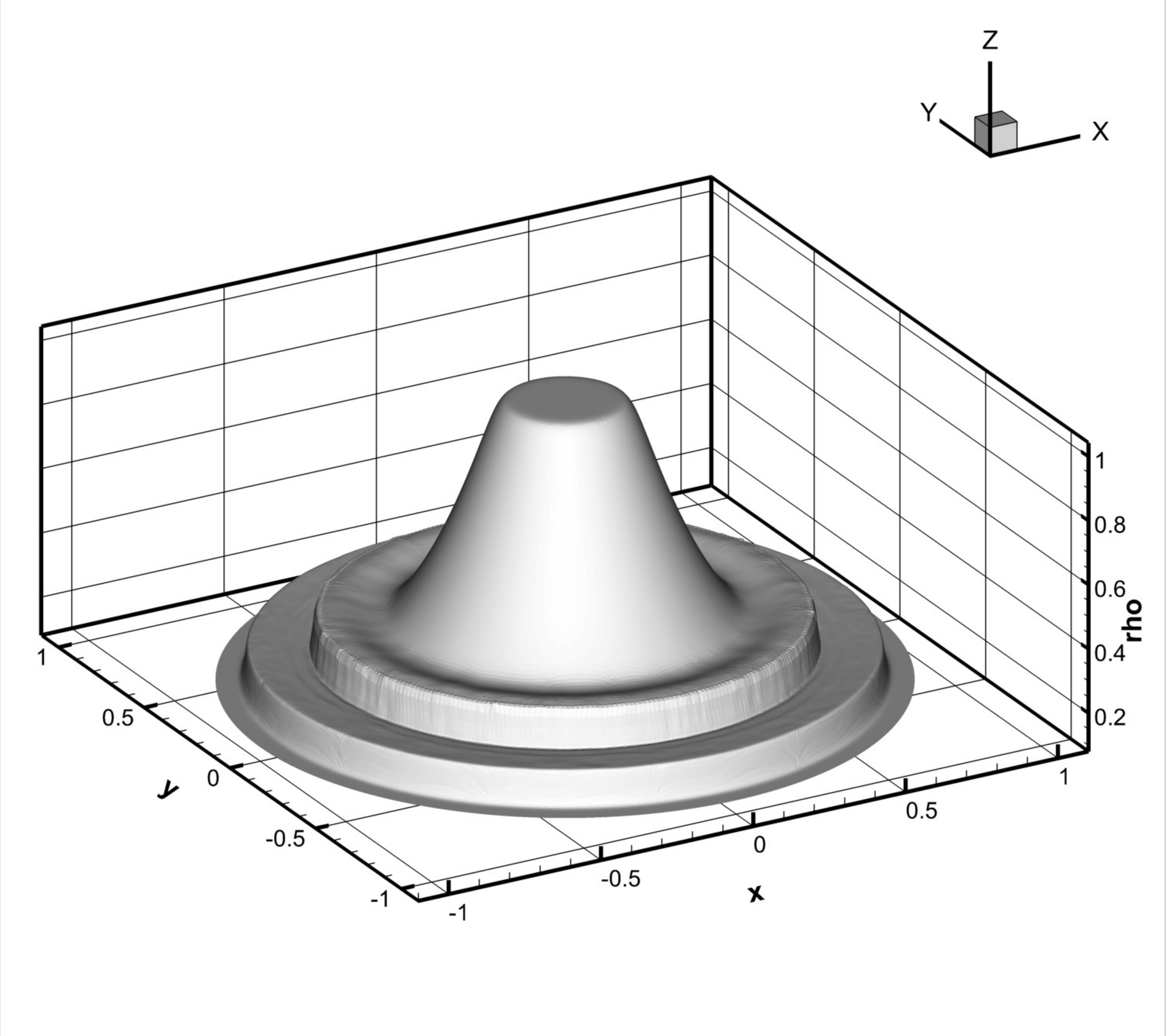}  &           
\includegraphics[width=0.47\textwidth]{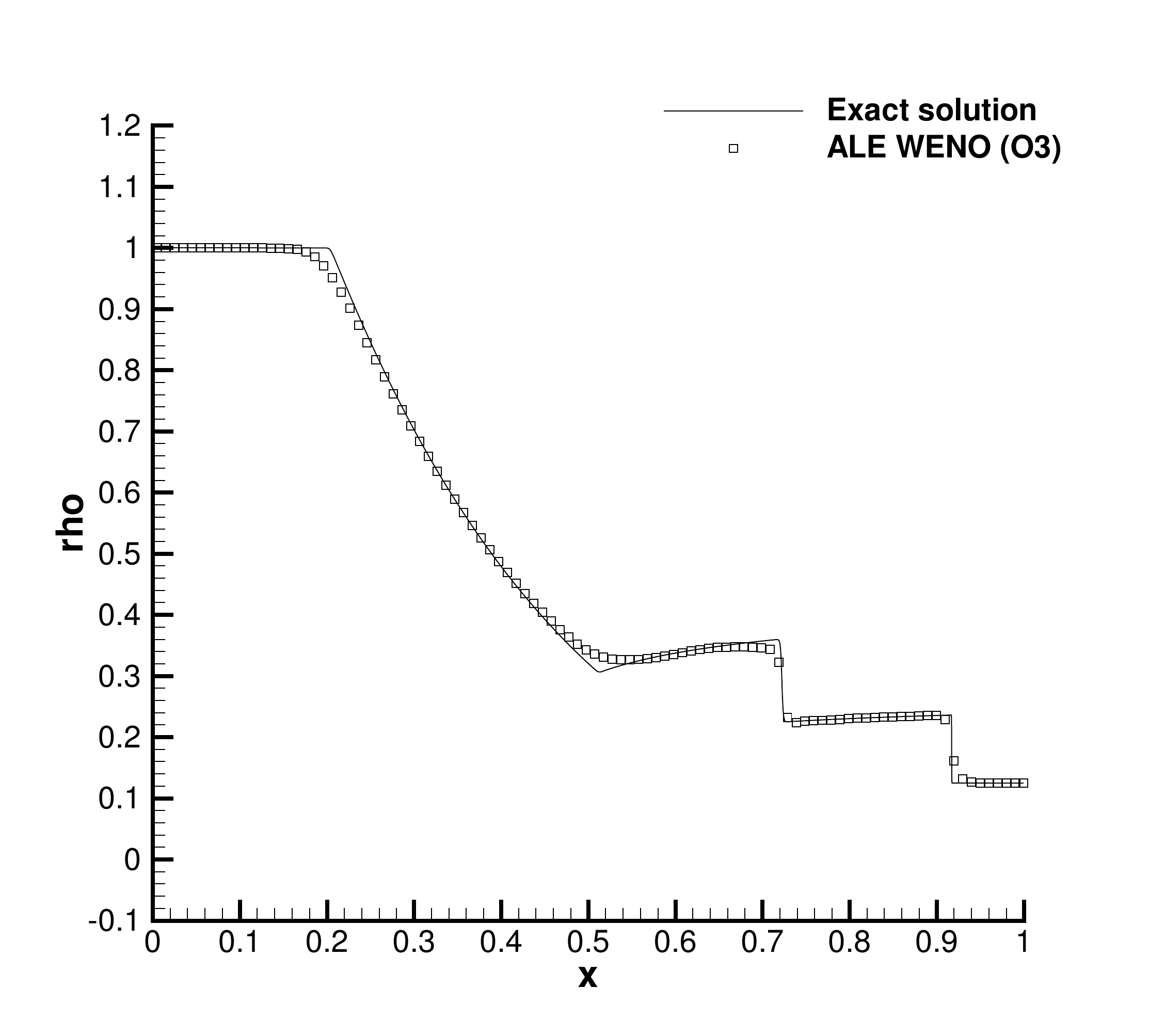} \\
\includegraphics[width=0.47\textwidth]{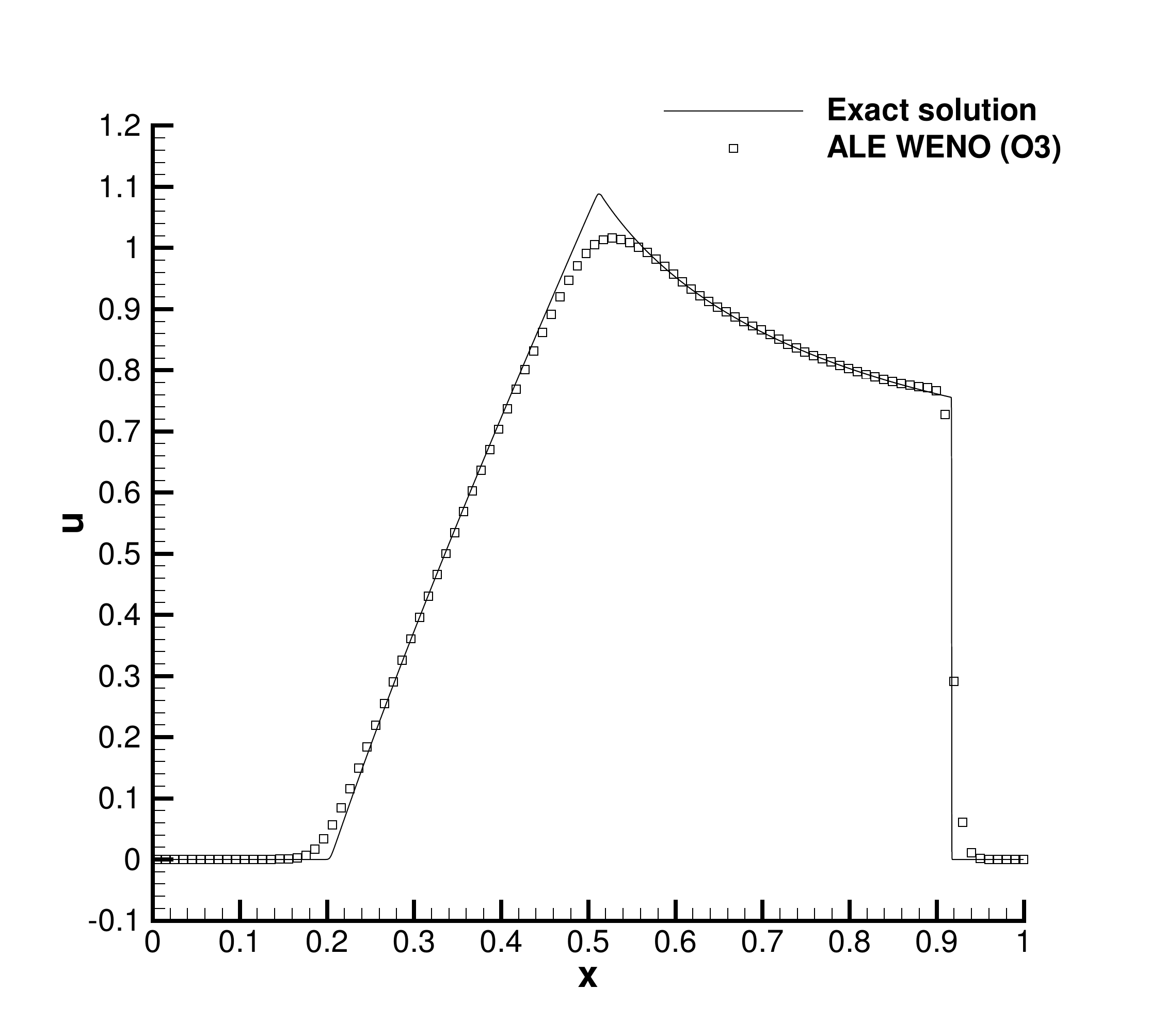}  &           
\includegraphics[width=0.47\textwidth]{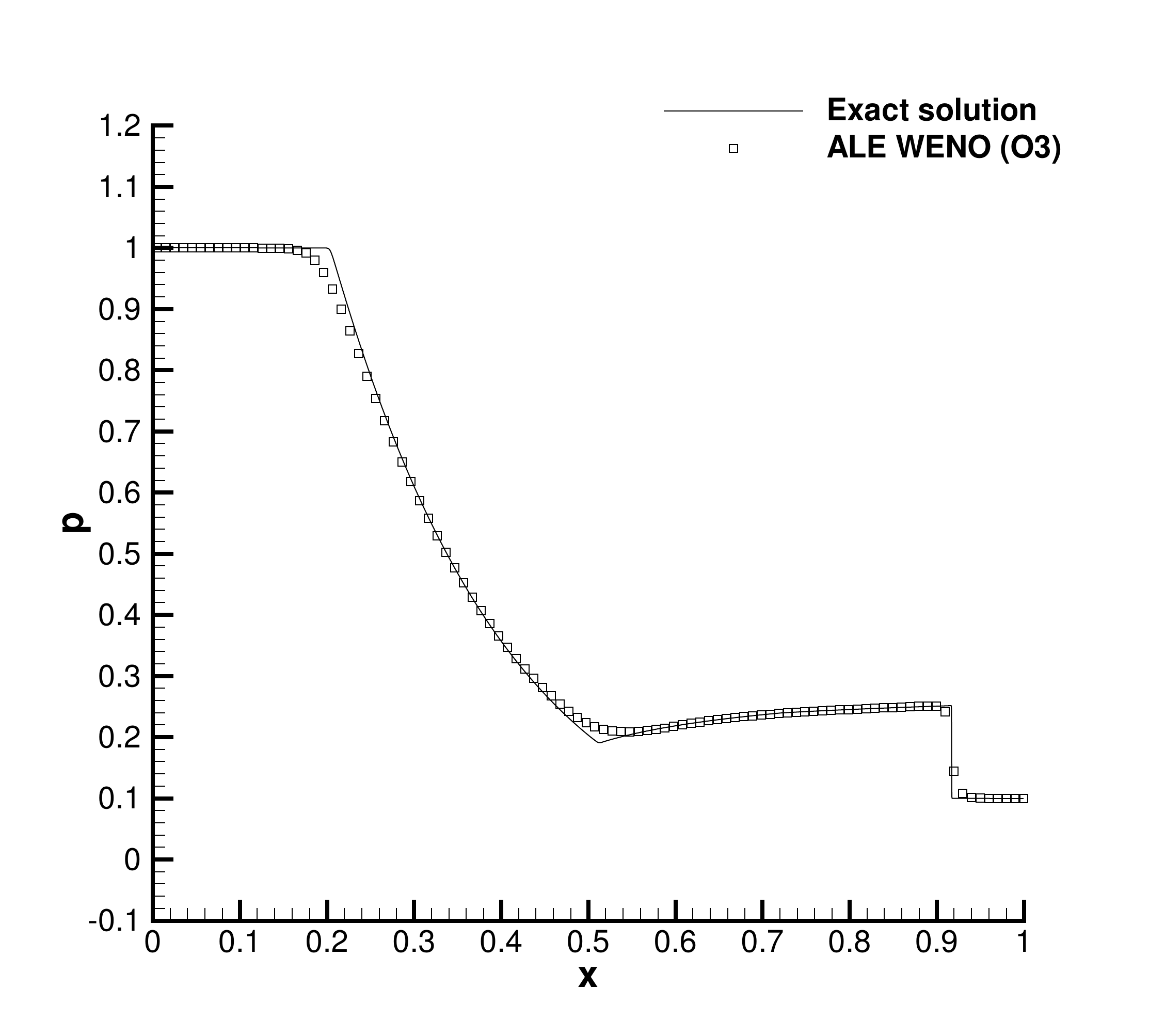} \\    
\end{tabular} 
\caption{Third order numerical results and comparison with exact solution for the two-dimensional explosion problem at time $t=0.25$.} 
\label{fig.EP2D}
\end{center}
\end{figure}

The reference solution is computed by solving a one-dimensional system with geometric source terms, as explained in \cite{BoscheriDumbserLag,ToroBook}. The Courant number is taken to be $\textnormal{CFL}=0.95$ and we use 
the HLLC flux to compute the numerical solution depicted in Figure \ref{fig.EP2D}: the solution consists in a circular rarefaction wave moving towards the center of the domain, a shock wave traveling outward and a contact wave in between them, which is very well resolved due to the use of a Lagrangian formalism.

\subsection{The Kidder Problem} 
\label{sec.Kidder}
\begin{figure}[!htbp]
\begin{center}
\begin{tabular}{cc} 
\includegraphics[width=0.47\textwidth]{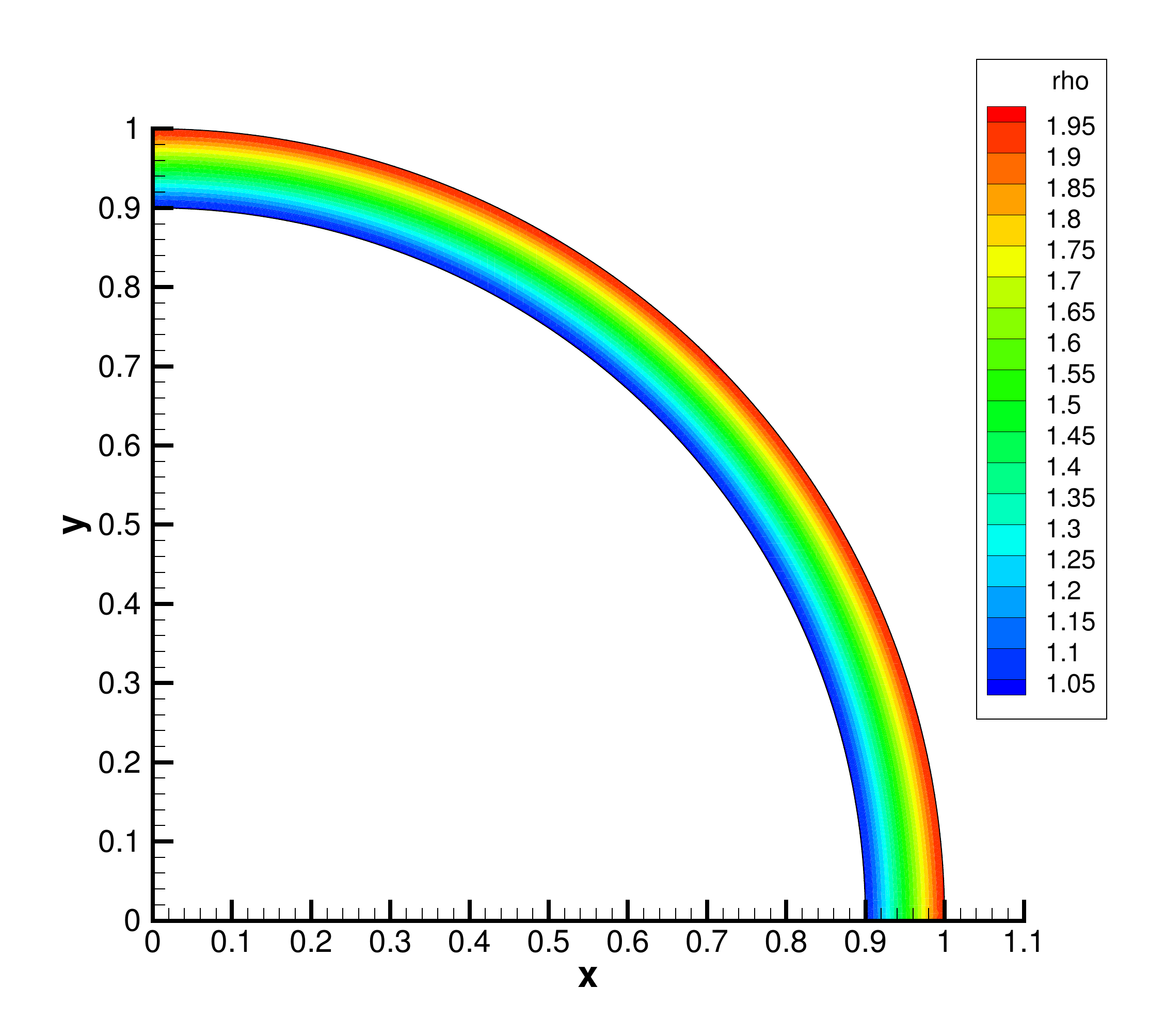}  &           
\includegraphics[width=0.47\textwidth]{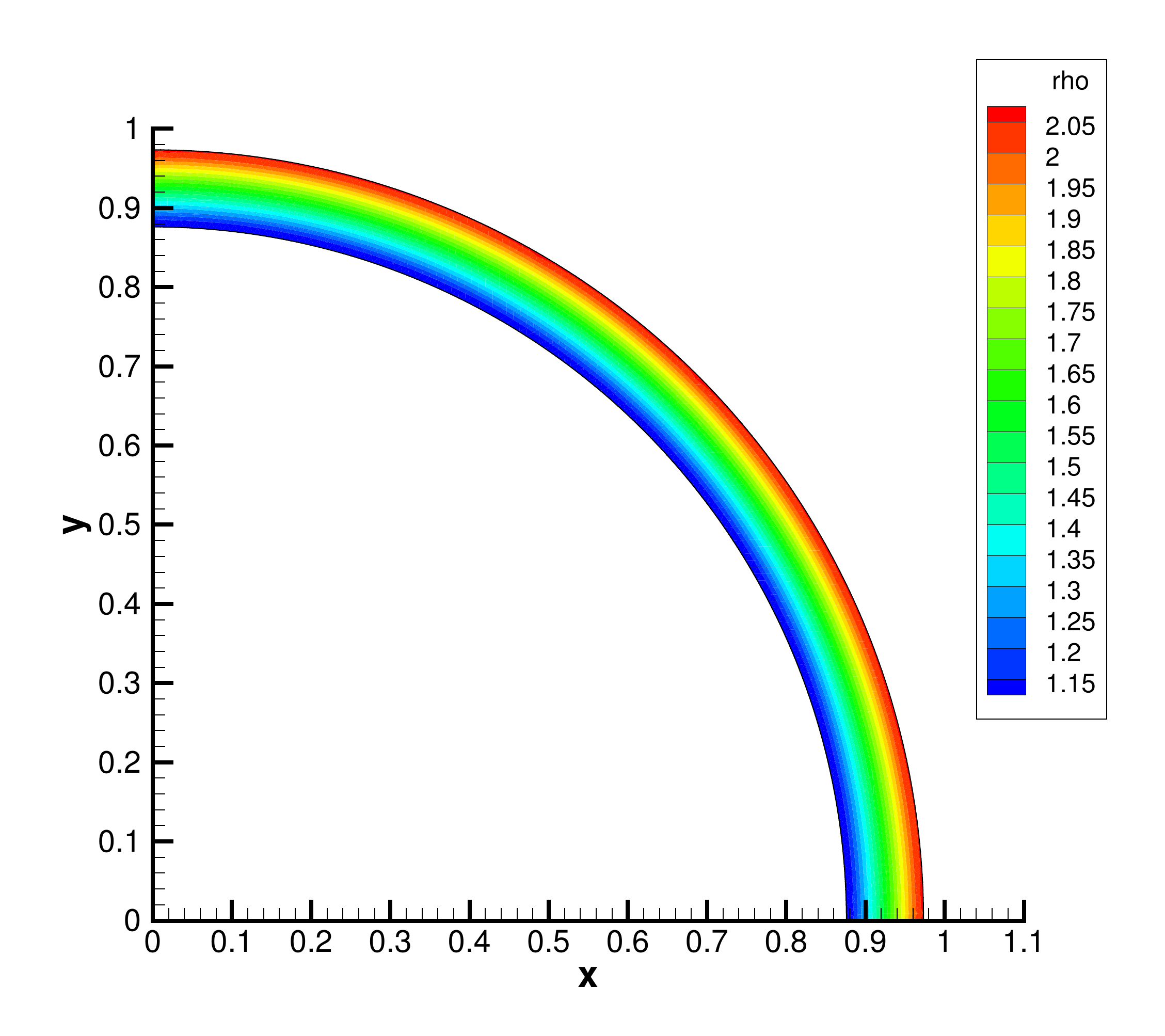} \\
\includegraphics[width=0.47\textwidth]{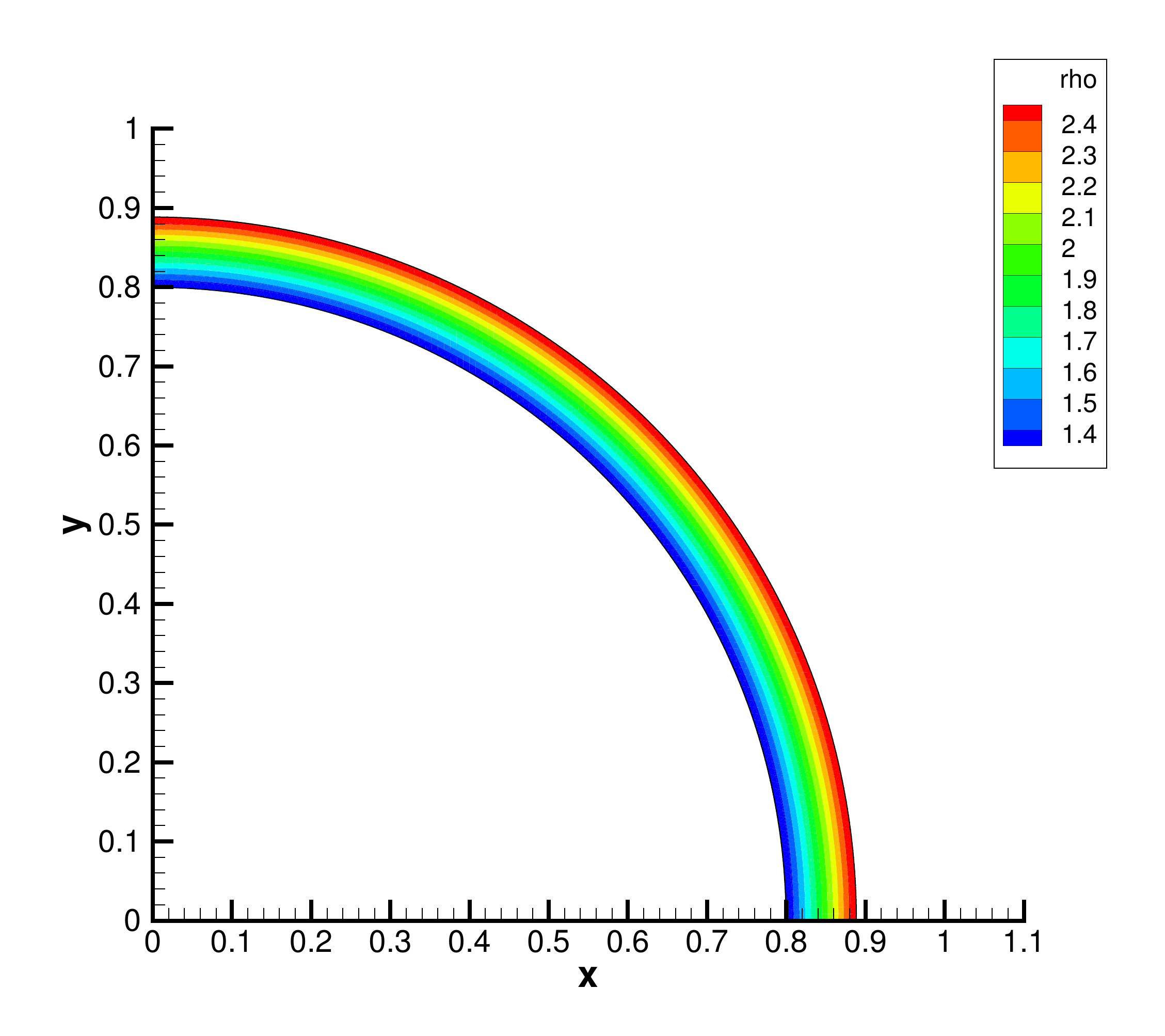}  &           
\includegraphics[width=0.47\textwidth]{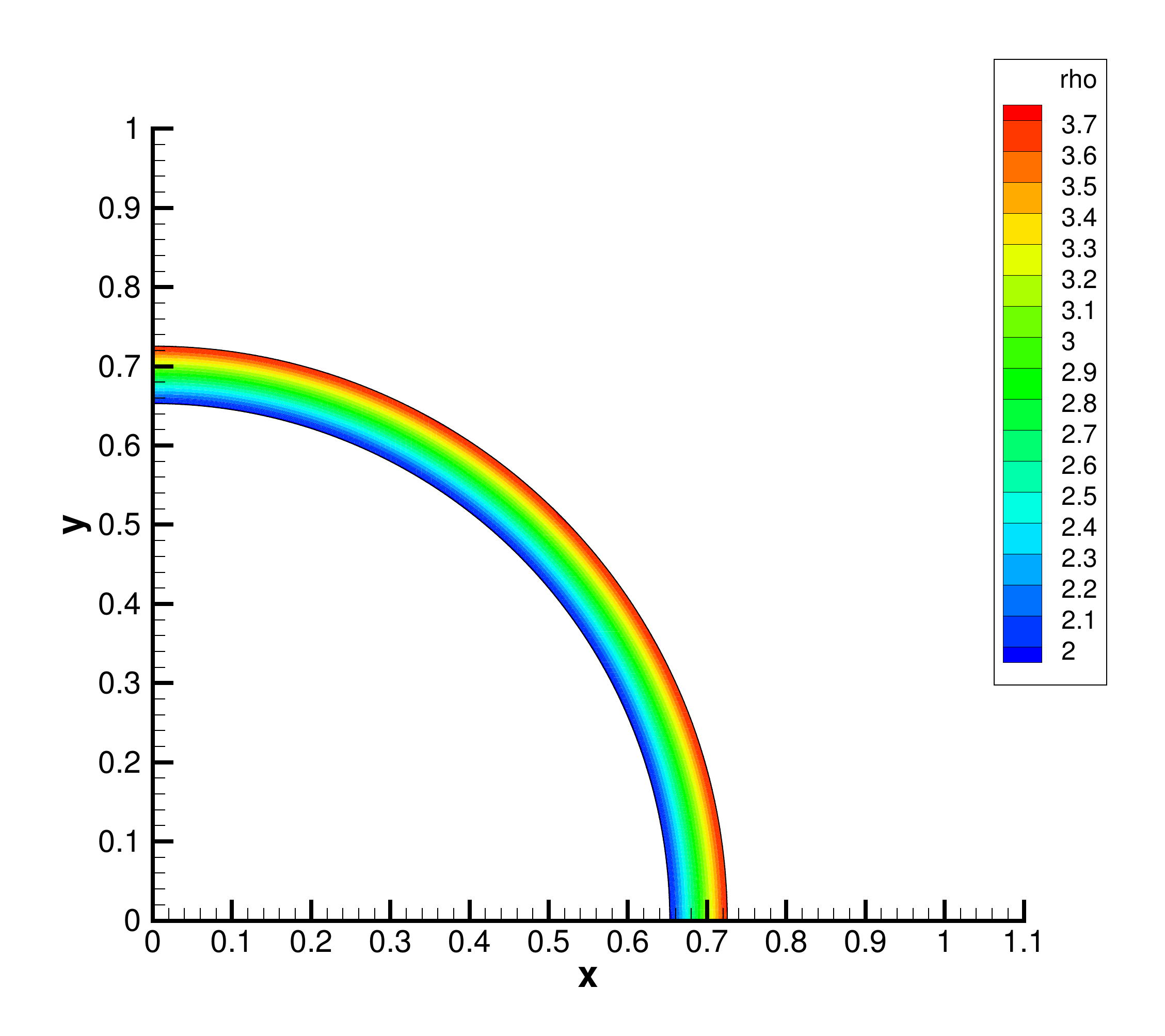} \\   
\includegraphics[width=0.47\textwidth]{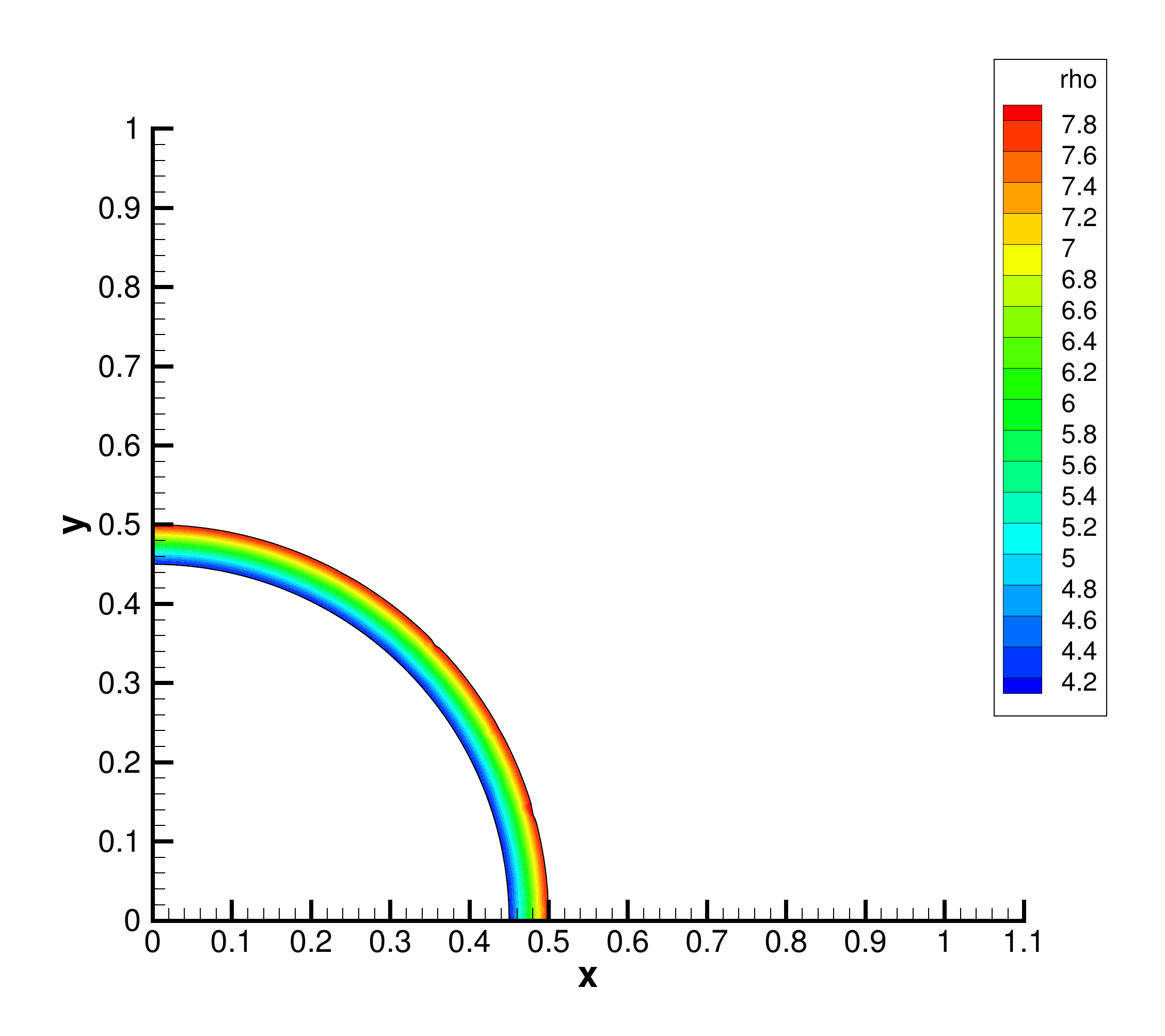}  &           
\includegraphics[width=0.47\textwidth]{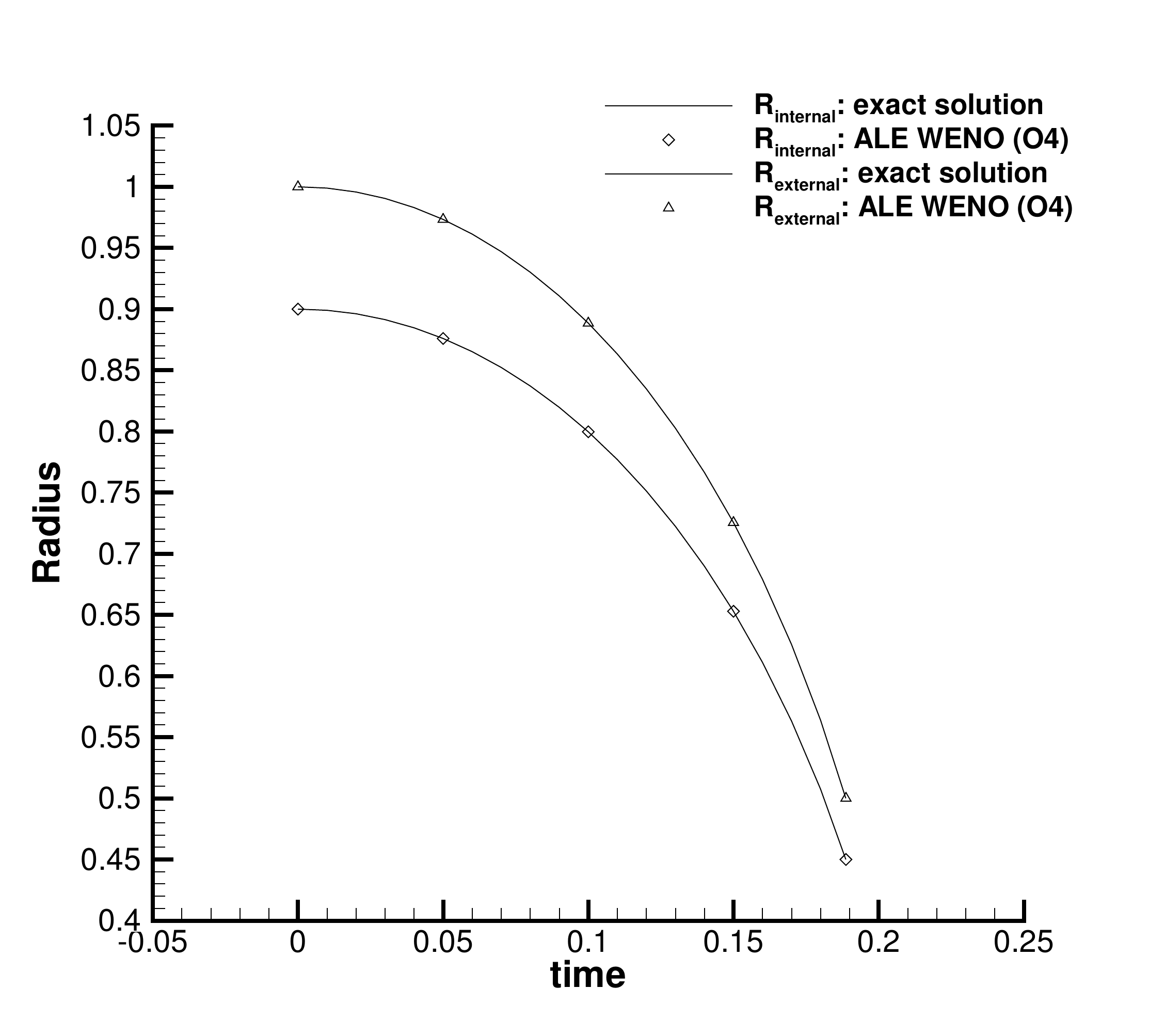} \\ 
\end{tabular} 
\caption{Density distribution for the Kidder problem at output times $t=0.00$, $t=0.05$, $t=0.10$, $t=0.15$ and $t=t_{f}$ (from top left to bottom left). Evolution of the internal and external radius of the shell and comparison between analytical and numerical solution (bottom right).} 
\label{fig.Kidder}
\end{center}
\end{figure}

In \cite{Kidder1976} Kidder proposed this test problem, which consists in an isentropic compression of a shell filled with an ideal gas. A self-similar analytical solution is available and can be used to check whether the numerical scheme generates spurious entropy during the isentropic compression, or not. The computational domain is a portion of a shell initially bounded by $r_i(t) \leq r \leq r_e(t)$, where $r$ denotes the general radial coordinate while $r_i(t),r_e(t)$ represent the time-dependent internal and external radius, respectively. The exact solution for a fluid particle initially located at radius $r$ is expressed as a function of the radius and the homothety rate $h(t)$, 
\begin{equation}
  R(r,t) = h(t)r, \qquad h(t) = \sqrt{1-\frac{t^2}{\tau^2}},
\label{eqKidderEx}
\end{equation} 
where $\tau$ denotes the focalisation time and is computed as 
\begin{equation}
\tau = \sqrt{\frac{\gamma-1}{2}\frac{(r_{e,0}^2-r_{i,0}^2)}{c_{e,0}^2-c_{i,0}^2}},
\end{equation}
with $c_{i,e}=\sqrt{\gamma\frac{p_{i,e}}{\rho_{i,e}}}$ the sound speeds at the inner and outer boundary, respectively. The initial density distribution $\rho_0$ is given by 
\begin{equation}
\rho_0 = \rho(r,0) = \left(\frac{r_{e,0}^2-r^2}{r_{e,0}^2-r_{i,0}^2}\rho_{i,0}^{\gamma-1}+\frac{r^2-r_{i,0}^2}{r_{e,0}^2-r_{e,0}^2}\rho_{e,0}^{\gamma-1}\right)^{\frac{1}{\gamma-1}} 
\end{equation}
where $r_i(0)=r_{i,0}=0.9$ and $r_e(0)=r_{e,0}=1.0$ are the initial values for the internal and external radius, respectively, while $\rho_{i,0}=1$ and $\rho_{e,0}=2$ give the initial values of density defined at the internal and at the external frontier of the shell, respectively. The ratio of specific heats is taken to be $\gamma=2$ and the initial velocity field is set to zero, i.e. $u=v=0$. We assume a uniform initial entropy, i.e. $s_0= \frac{p_0}{\rho_0^\gamma} = 1$, hence the initial pressure distribution is expressed as $p_0(r) = s_0\rho_0(r)^\gamma$. Sliding wall boundary conditions are set on the lateral faces of the shell, whereas the internal and the external frontier are assigned with a space-time dependent state, which is computed according to the exact analytical solution $R(r,t)$ (see \cite{Kidder1976} for details). As done in \cite{Despres2009,Maire2009}, the final time is taken to be $t_f=\frac{\sqrt{3}}{2}\tau$, so that the compression rate is $h(t_f)=0.5$ and the exact location of the shell is delimited by $0.45 \leq r \leq 0.5$. Figure \ref{fig.Kidder} shows the numerical results obtained with a fourth order version of the ALE WENO scheme together with the multidimensional HLLC flux on a computational grid with a characteristic mesh size of $h=1/100$. The CFL number used was $\textnormal{CFL}=0.95$. The  evolution of the density distribution has been plotted as well as the time-dependent location of the internal and the external frontier. Furthermore Table \ref{tab:radiusKidder} reports the absolute error $|err|$ of the frontier positions, which is defined as the difference between the analytical and the numerical location of the internal and external radius at the final time.

\begin{table}[!htbp]
	\begin{center}
		\begin{tabular}{|c|c|c|c|}
		\hline
		    		  			& $R_{ex}$ 		& $R_{num}$  & $|err|$     \\
		\hline
		\textit{Internal radius}	& 0.45000000 	& 0.45000031  & 0.31E-06 \\
		\hline
		\textit{External radius}	& 0.50000000 	& 0.50000613  & 6.13E-06 \\
		\hline
		\end{tabular}
	\end{center}
	\caption{Absolute error for the internal and external radius location between exact $R_{ex}$ and numerical $R_{num}$ solution.}
	\label{tab:radiusKidder}
\end{table}

\subsection{The Saltzman Problem} 
\label{sec.Saltzman}
\begin{figure}[htbp]
	\centering
		\includegraphics[width=0.75\textwidth]{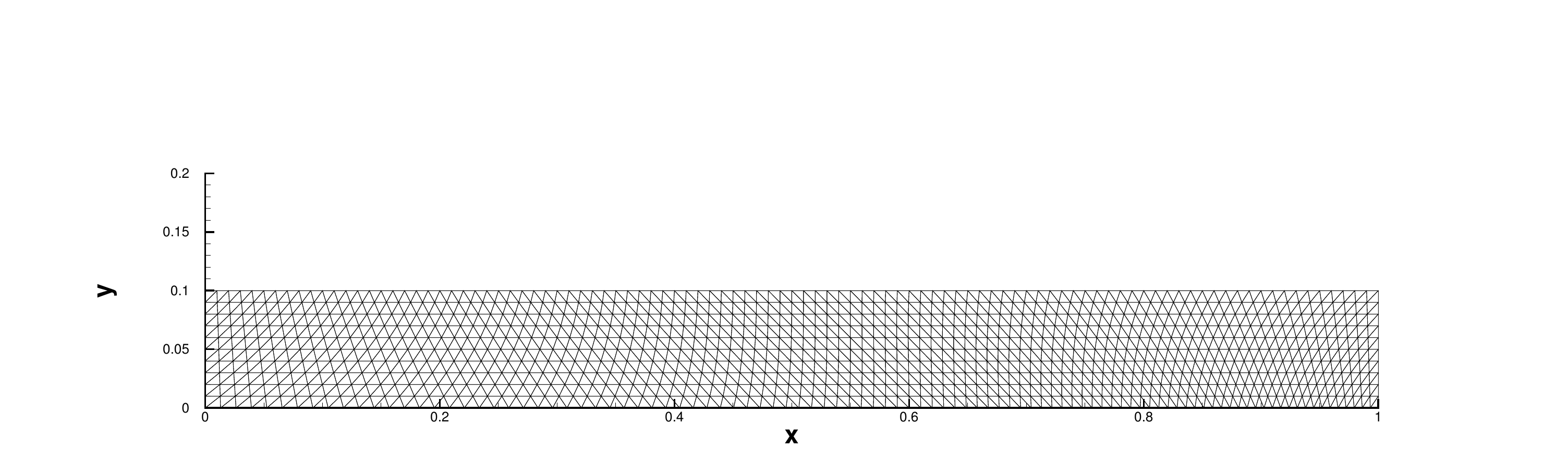}
		\caption{Initial mesh configuration for the Saltzman problem with a total number of elements of $N_E=2 \cdot 100\times10=2000$.}
	\label{fig:Saltz_initialGrid}
\end{figure}

The Saltzman test problem was first proposed by Dukowicz et al. in \cite{SaltzmanOrg} and involves a strong one-dimensional shock wave driven by a piston that is pushing and compressing a gas contained in a closed channel. The initial rectangular domain is $\Omega(0)=[0;1]\times[0;0.1]$ and the computational mesh is composed of $2 \cdot 100 \times 10$ right-angled triangular elements, as depicted in Figure \ref{fig:Saltz_initialGrid}. The piston is moving with velocity $\mathbf{v}_p = (1,0)$ and initially the fluid is at rest with an internal energy of $e_0=10^{-4}$ and a density of $\rho_0=1$. Therefore the initial vector of conserved variables is $\Q_{0} = \left( \rho_0, u_0, v_0, \rho E \right) = \left( 1, 0, 0, 10^{-4} \right)$. According to \cite{chengshu2}, the ratio of specific heats is taken to be $\gamma = \frac{5}{3}$ and the final time is set to $t_f=0.6$. We impose moving slip wall boundary condition on the piston and fixed slip wall boundaries on the remaining sides of the domain.  

\begin{figure}[!htbp]
\begin{center}
\begin{tabular}{cc} 
\includegraphics[width=0.47\textwidth]{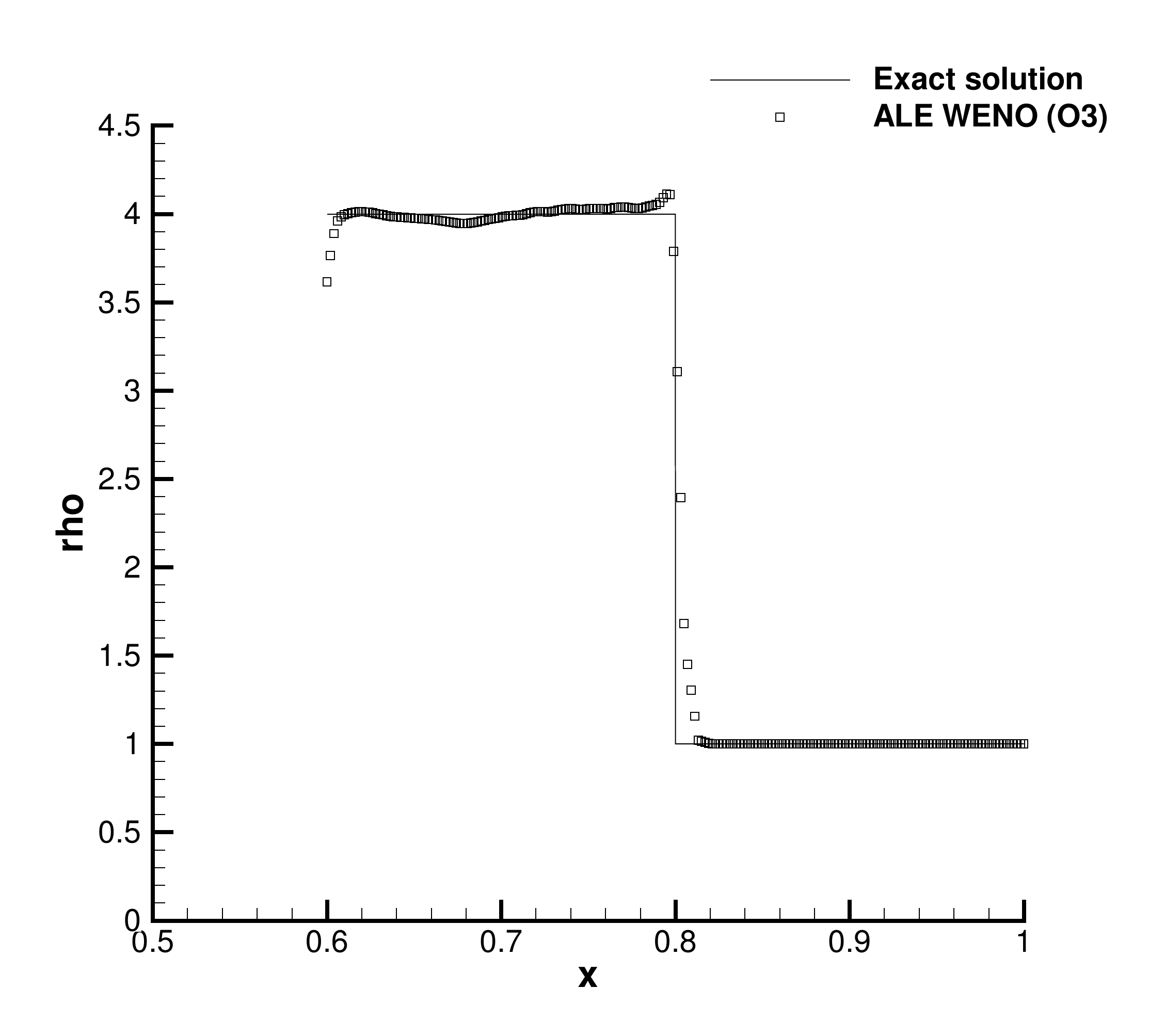}  &           
\includegraphics[width=0.47\textwidth]{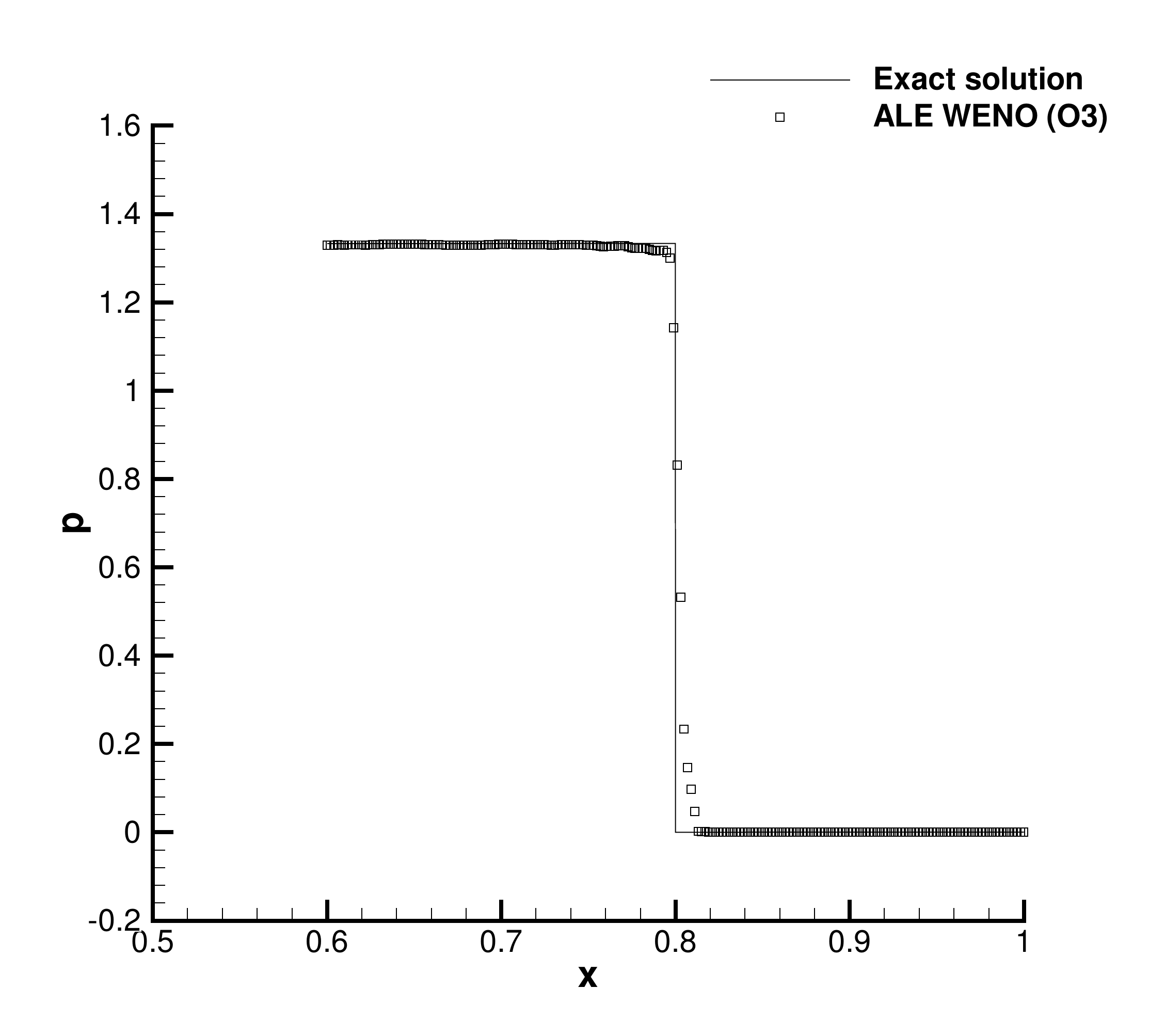} \\
\includegraphics[width=0.47\textwidth]{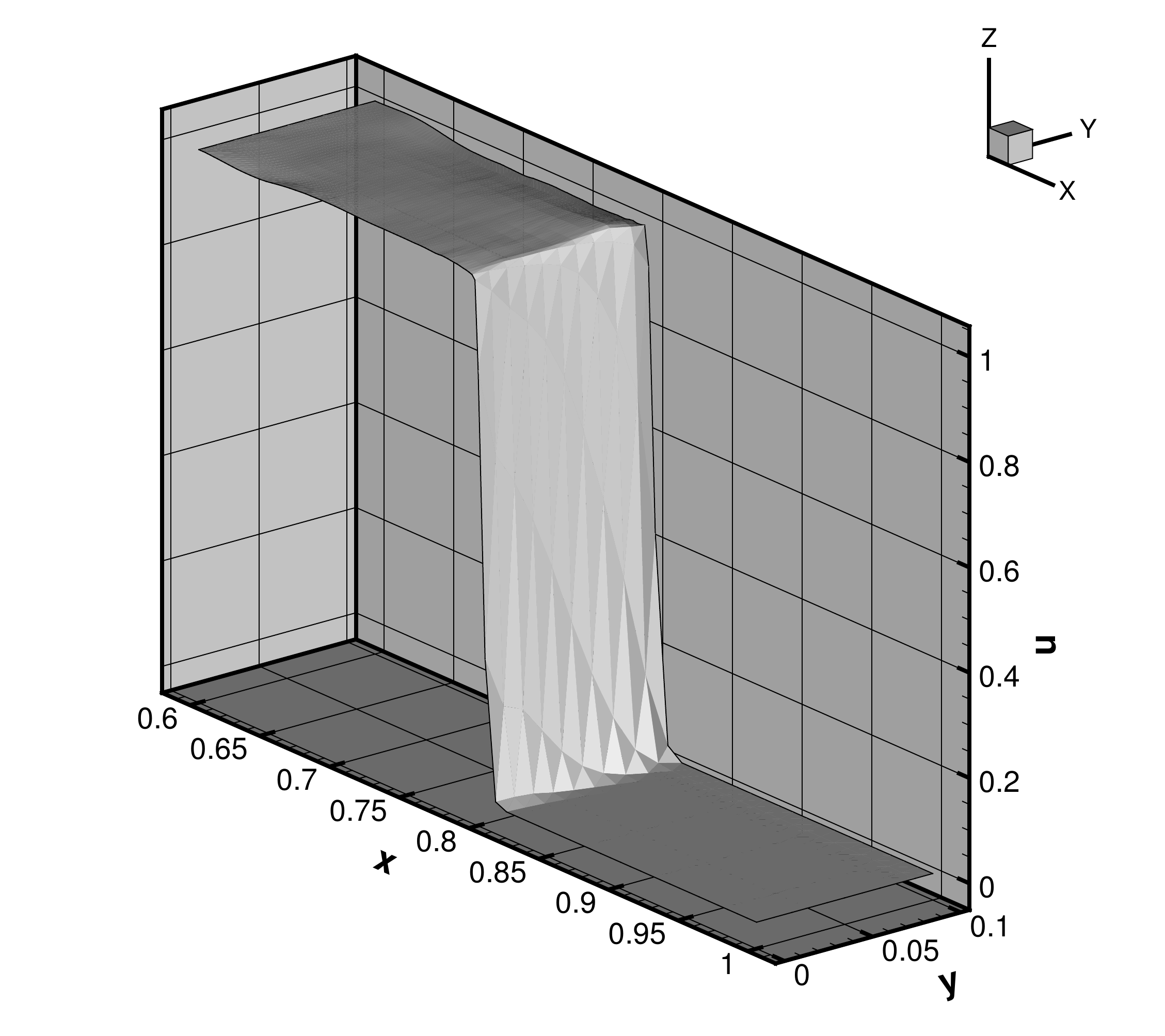}  &           
\includegraphics[width=0.47\textwidth]{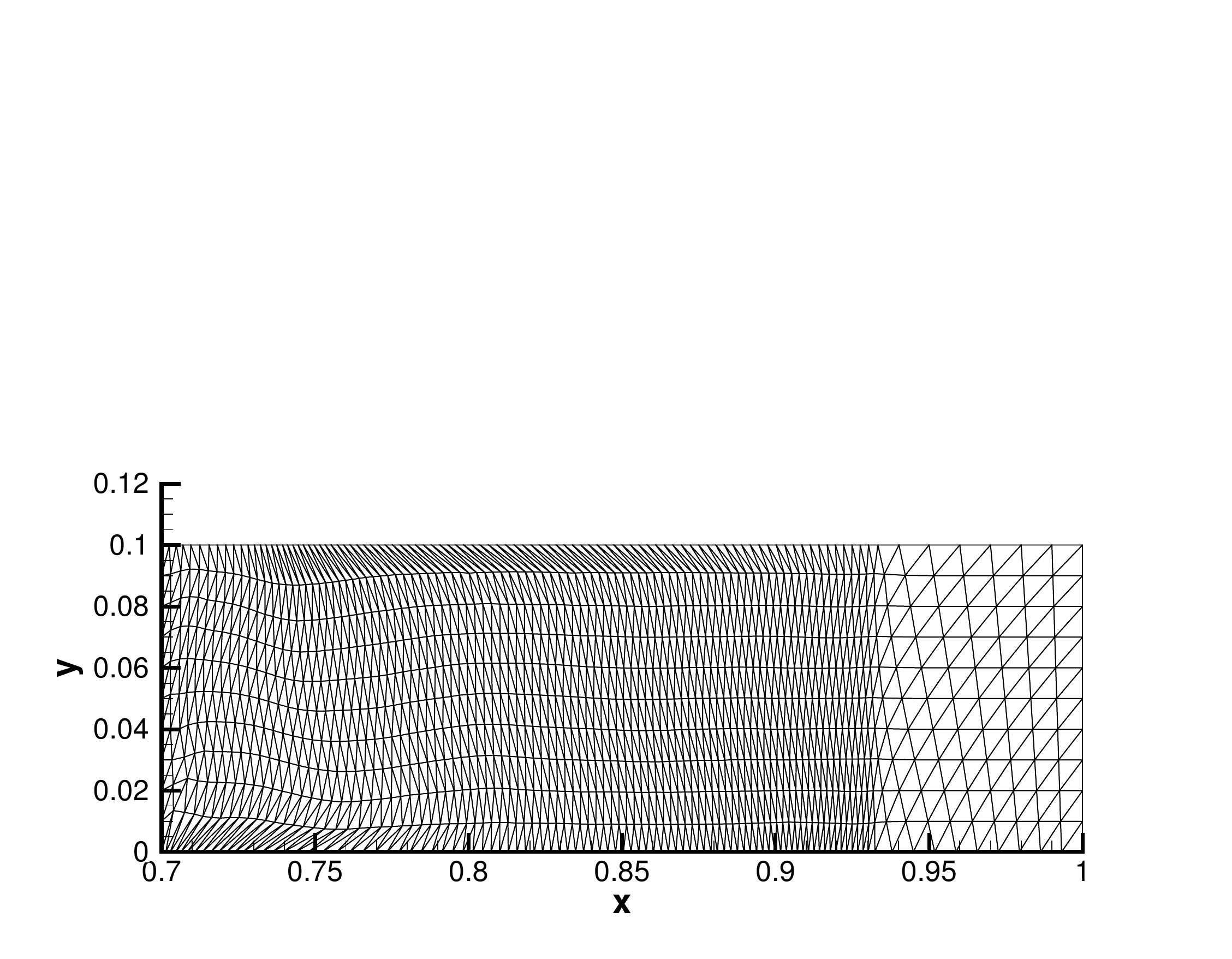} \\    
\end{tabular} 
\caption{Top: comparison between numerical and analytical solution for density and pressure for the Saltzman problem at time $t=0.6$. Bottom: velocity distribution ($t=0.6$) and mesh configuration at time $t=0.7$.} 
\label{fig.Saltzman}
\end{center}
\end{figure}

The exact solution $\Q_{ex}$ can be computed by solving a one-dimensional Riemann problem, \cite{BoscheriDumbserLag,ToroBook}, and reads:
\begin{equation}
  \Q_{ex}(\x,t_f) = \left\{ \begin{array}{ccc} \left( 4, 1, 0, 2.5     \right) & \textnormal{ if } & x \leq x_f, \\
                                               \left( 1, 0, 0, 10^{-4} \right) & \textnormal{ if } & x > x_f,        
                      \end{array}  \right. 
\end{equation}
where $x_f=0.8$ is the final shock location at time $t_f$. Since the piston is strongly compressing the fluid at the initial times of the simulation, particular care has to be taken in order to respect the geometric CFL condition of those elements that lie near the piston. That is why we could not run this challenging test problem adopting a Courant number higher than $0.7$, as done for the previous test cases. We use the third order version of the ADER-WENO ALE scheme and the multidimensional HLL flux to obtain the results shown in Figure \ref{fig.Saltzman}, where a good agreement with the exact solution can be noticed. 

\subsection{The Sedov Problem} 
\label{sec.Sedov}
\begin{figure}[!htbp]
\begin{center}
\begin{tabular}{cc} 
\includegraphics[width=0.47\textwidth]{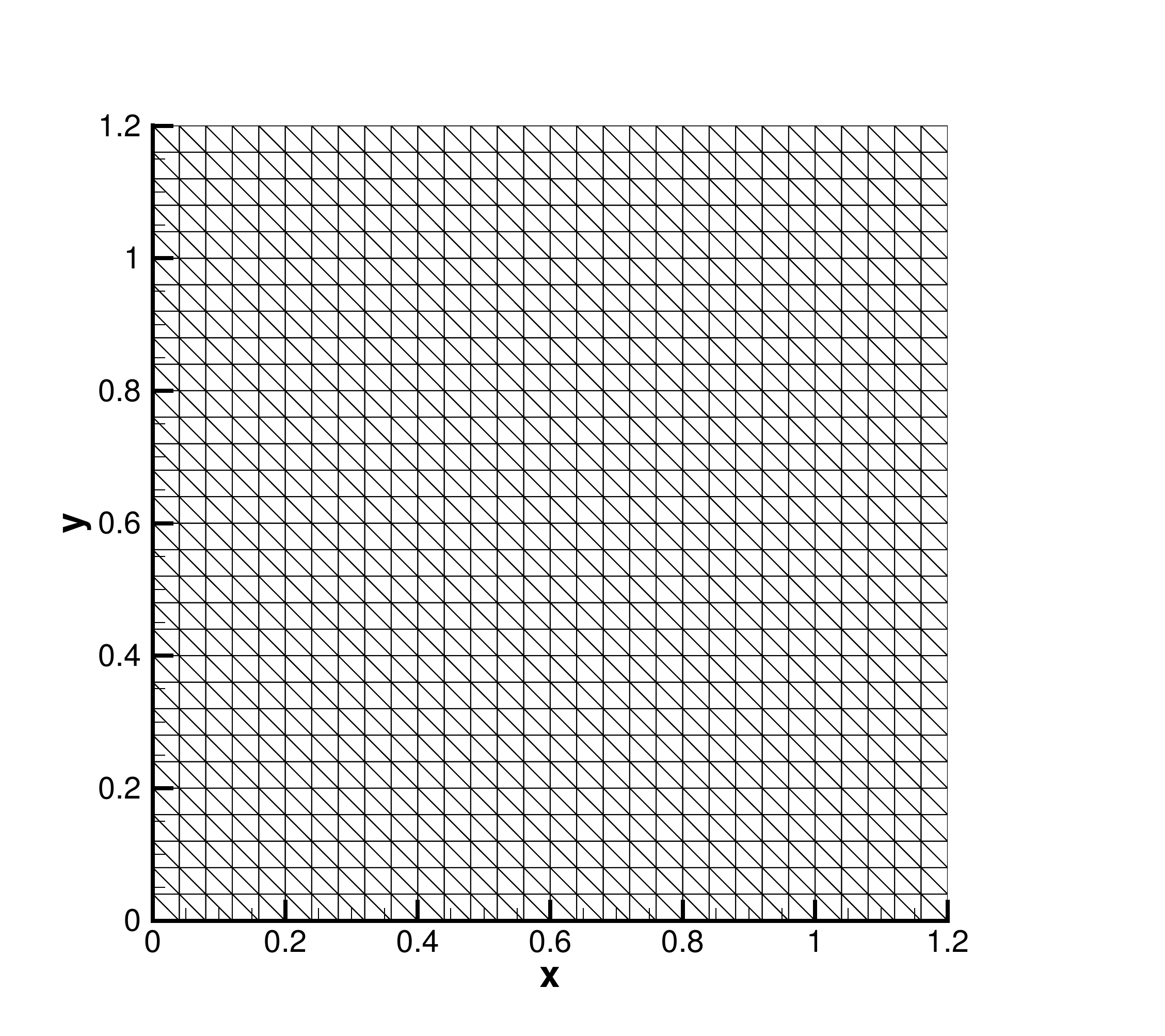}  &           
\includegraphics[width=0.47\textwidth]{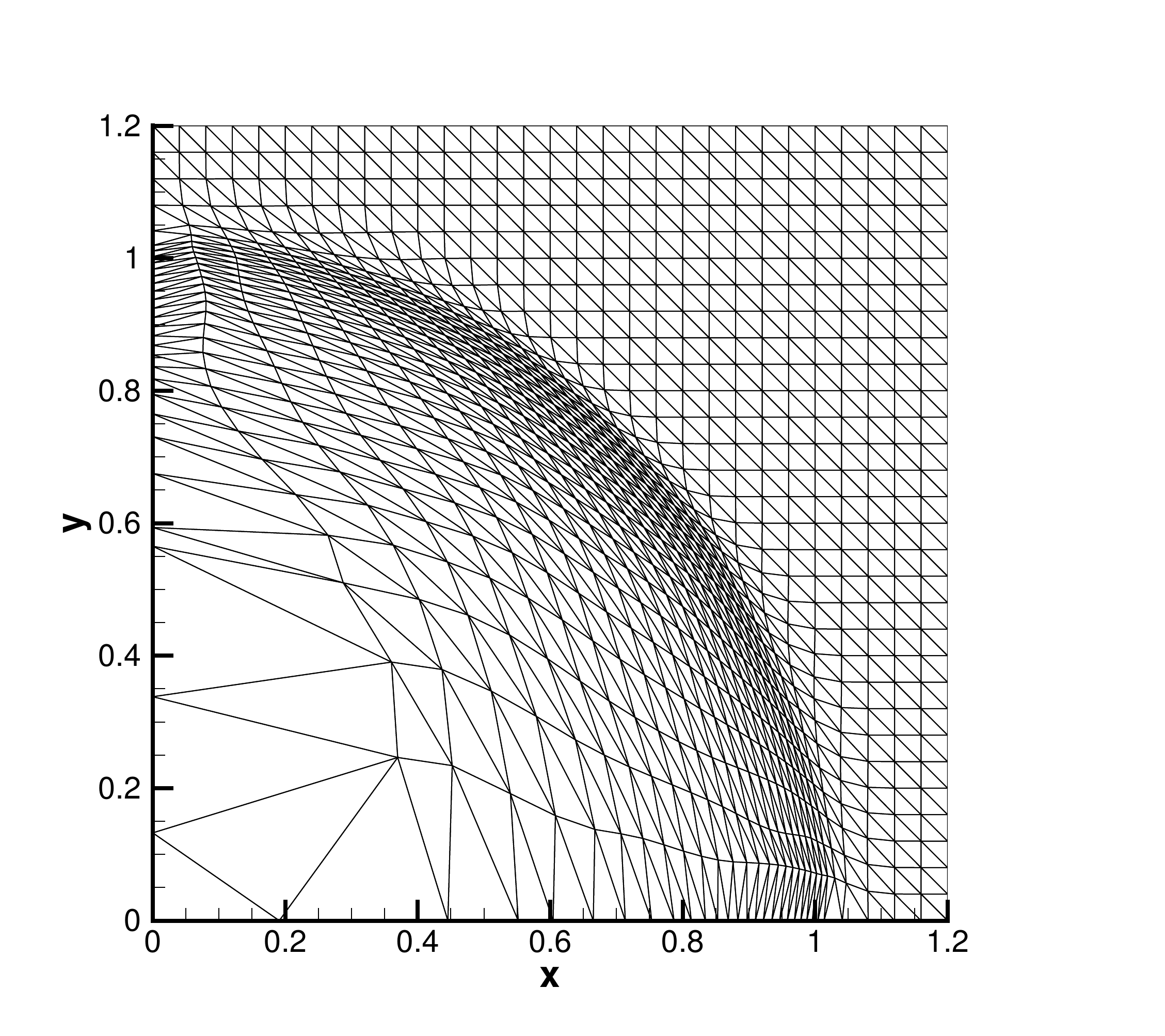} \\
\includegraphics[width=0.47\textwidth]{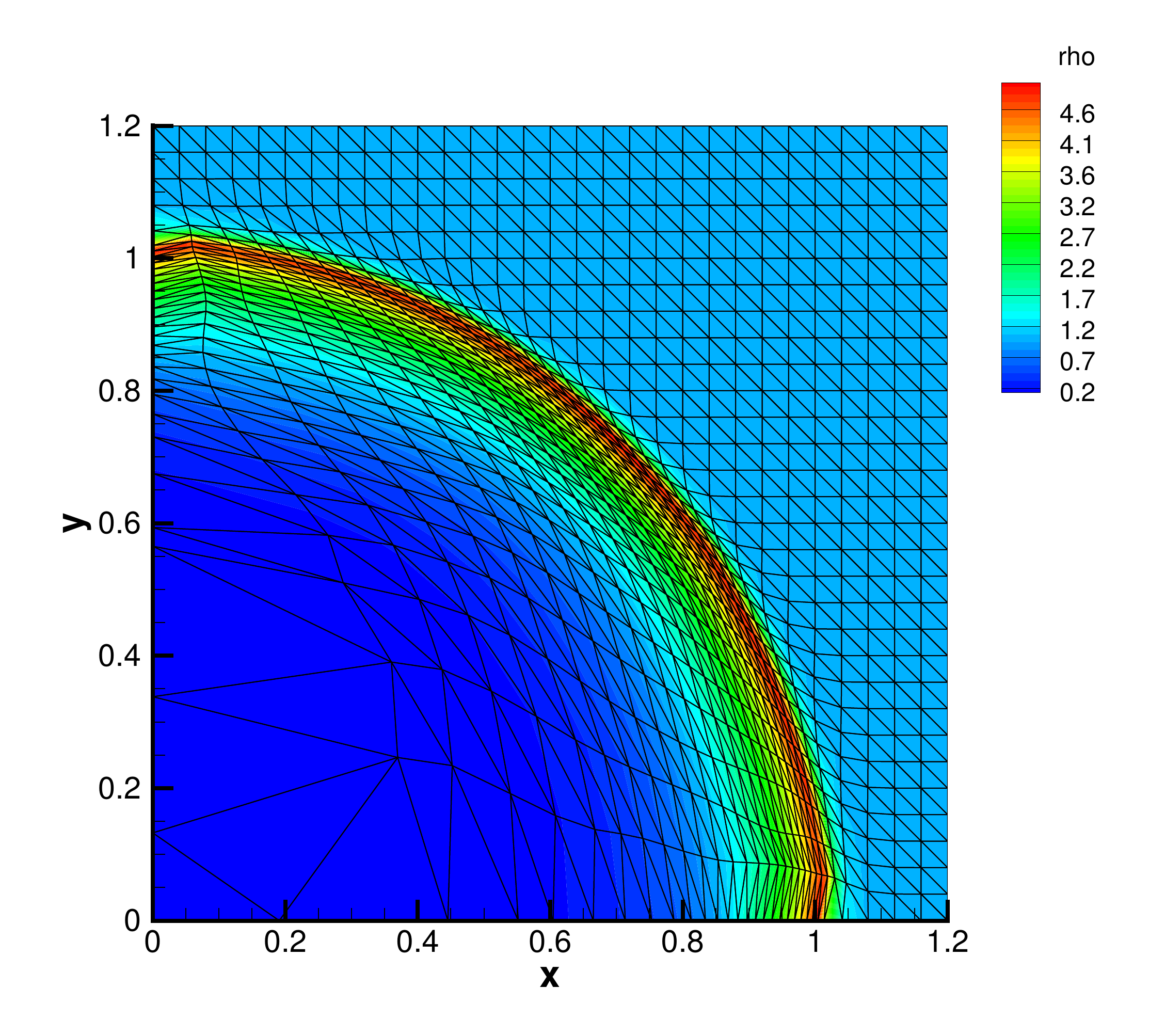}  &           
\includegraphics[width=0.47\textwidth]{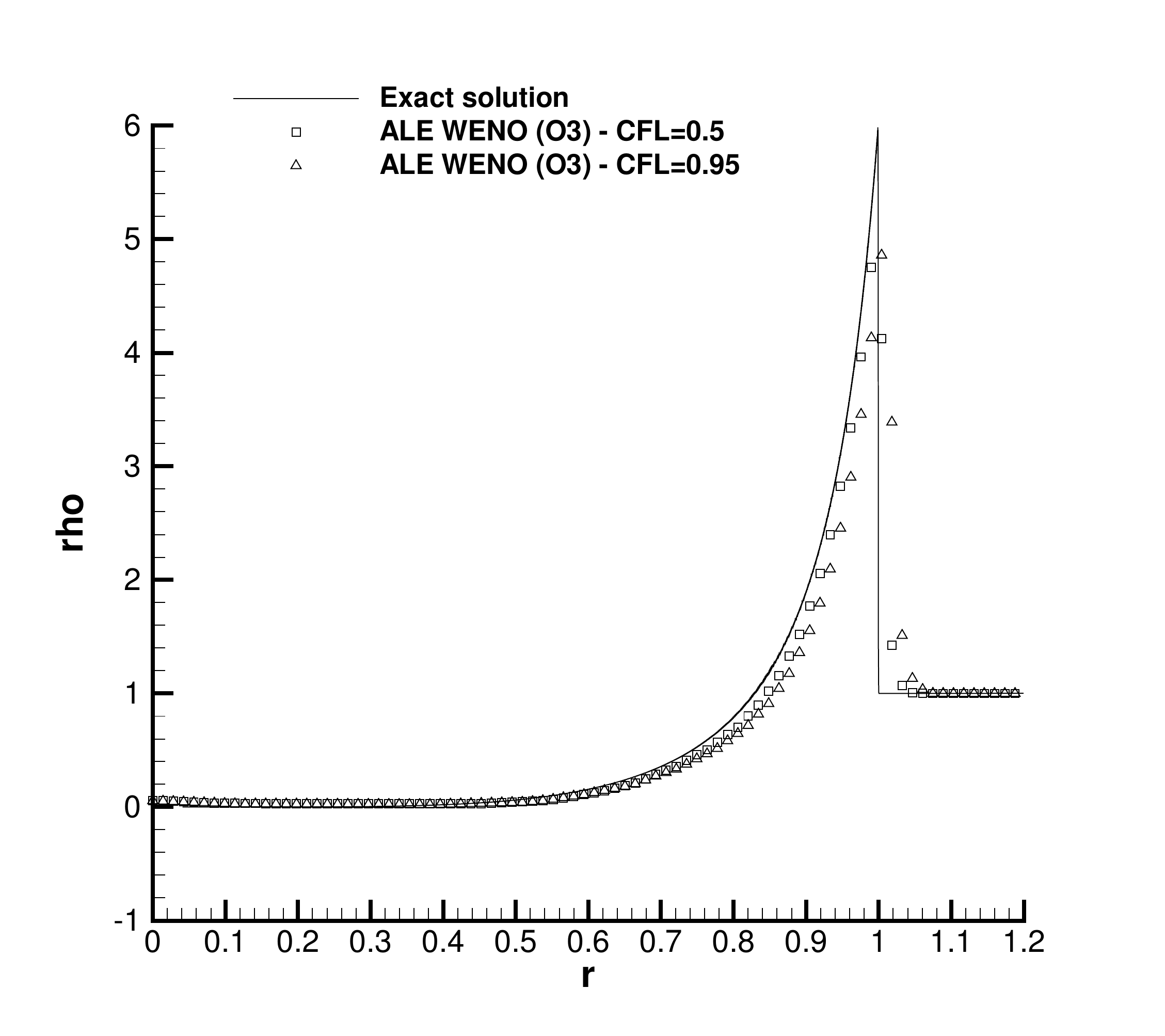} \\    
\end{tabular} 
\caption{Top: initial and final mesh configuration for the Sedov problem. Bottom: density distribution at the final time $t_f=0.6$ and comparison between the exact solution (solid line) and two different third order accurate numerical solution obtained with $\textnormal{CFL}=0.5$ and $\textnormal{CFL}=0.95$.} 
\label{fig.Sedov2D}
\end{center}
\end{figure}

A point-symmetric explosion with the generation of a blast wave describes the Sedov problem. It is a very challenging test case for which Kamm et al. \cite{SedovExact} proposed an exact solution with cylindrical symmetry which depends on self-similarity arguments. According to \cite{SedovExact}, the gas has a unity initial density $\rho_0=1$ and a quasi-zero initial pressure $p_0=10^{-6}$, which is imposed to the whole computational domain except at the origin $O=(0,0)$, where the pressure is set to
\begin{equation}
p_{or} = (\gamma-1)\rho_0 \frac{\epsilon_0}{V_{or}},
\label{eqn.p0.sedov}
\end{equation}  
with  $\epsilon_0=0.244816$ denoting the total amount of released energy, $V_{or}$ representing the volume of the cell $T_{or}$ located at the origin and $\gamma=1.4$ being the ratio of specific heats. The initial computational domain is a square $\Omega(0)=[0;1.2]\times[0;1.2]$ and the initial mesh is composed by $(30\times30)$ square elements, each of those has been split into two right-angled triangles. Since this test case was first proposed for Cartesian grids, the volume $V_{or}$ of the origin cell is here taken to be the volume of the two triangles which compose the square element located at the origin of the domain. The exact position of the cylindrical shock wave at the final time of the simulation $t_f=1$ is at radius $r=\sqrt{x^2+y^2}=1$. We impose sliding wall boundary conditions on each side of the domain. Figure \ref{fig.Sedov2D} shows the numerical results obtained with a third order ADER-WENO ALE scheme using the multidimensional HLL flux. 
The mesh is highly distorted and compressed by the shock wave, but the numerical solution agrees well with the exact solution, as depicted in Figure 
\ref{fig.Sedov2D}. The rezoning step described in Section \ref{sec.meshmotion} was necessary in order to reduce the mesh deformation and to avoid tangled elements.
 
\subsection{The Noh Problem} 
\label{sec.Noh}
\begin{figure}[!htbp]
\begin{center}
\begin{tabular}{cc} 
\includegraphics[width=0.47\textwidth]{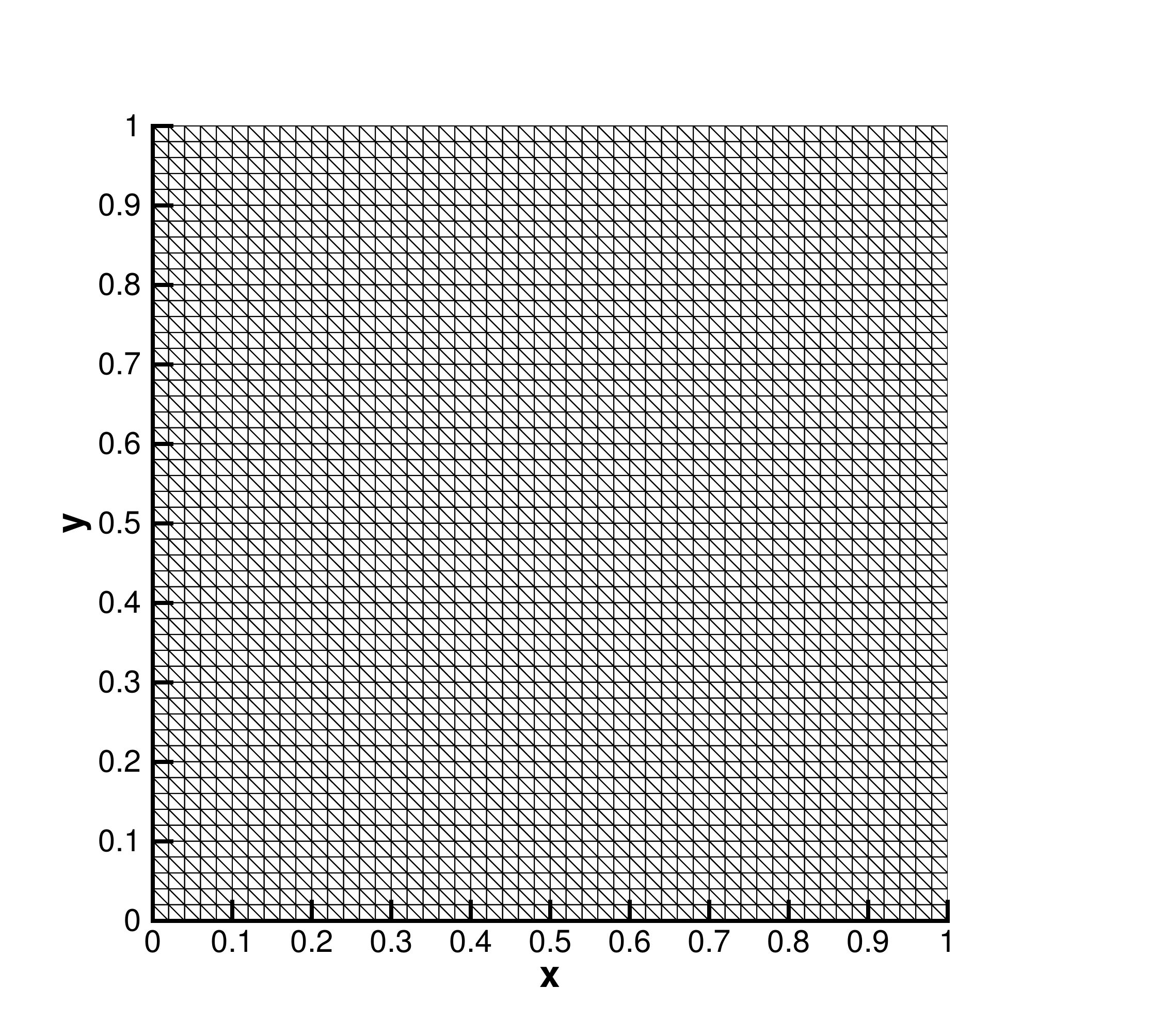}  &           
\includegraphics[width=0.47\textwidth]{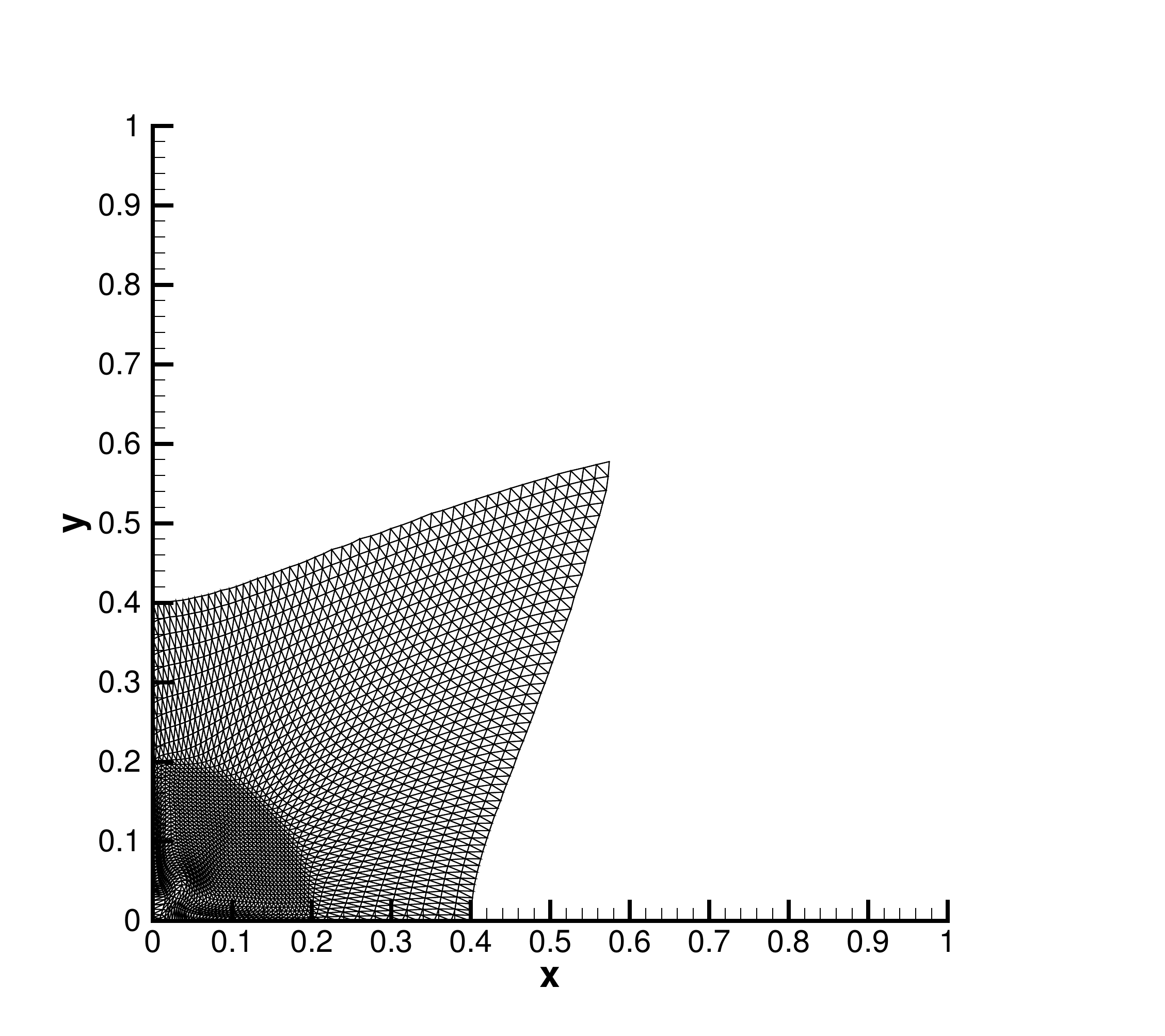} \\
\includegraphics[width=0.47\textwidth]{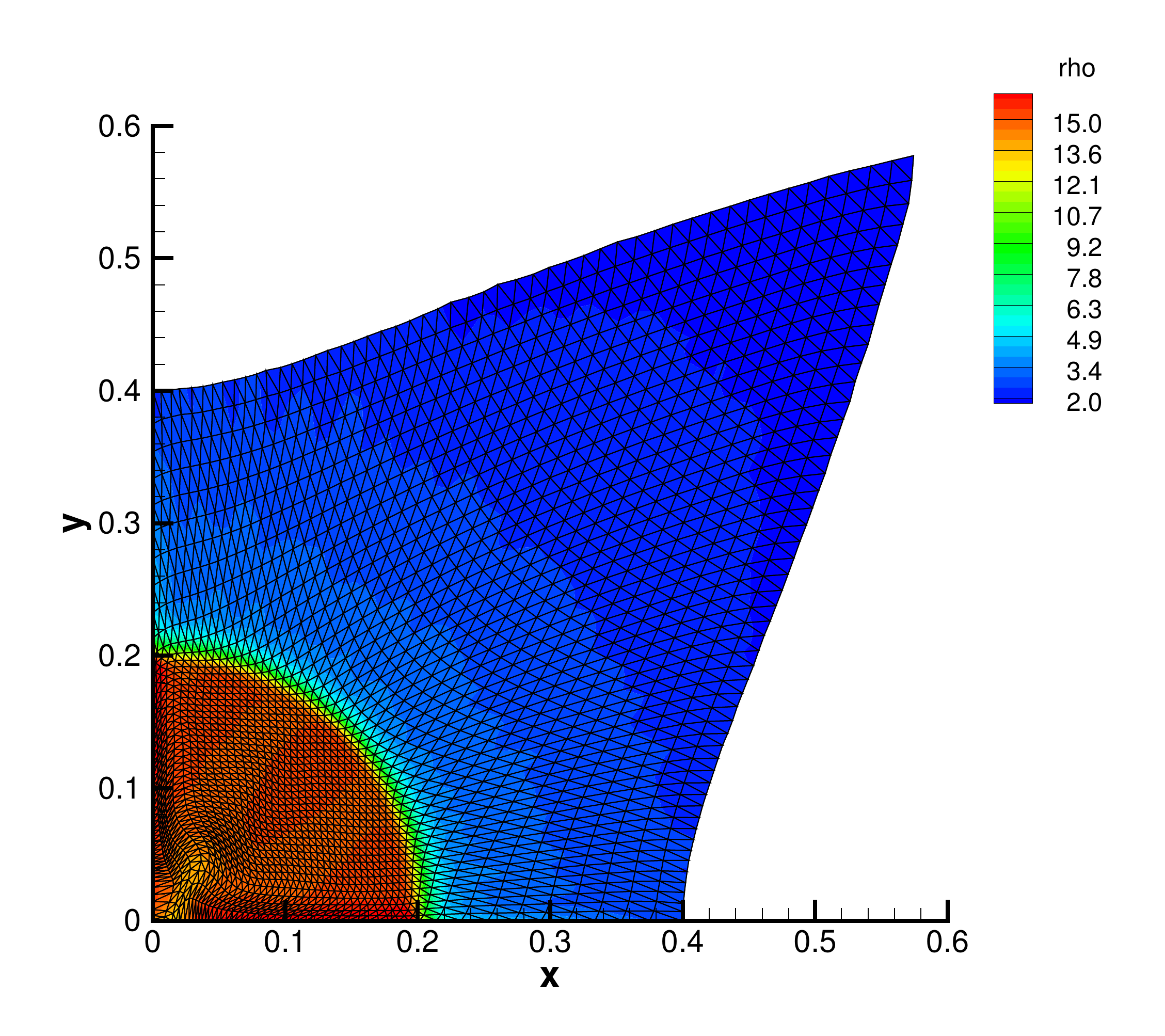}  &           
\includegraphics[width=0.47\textwidth]{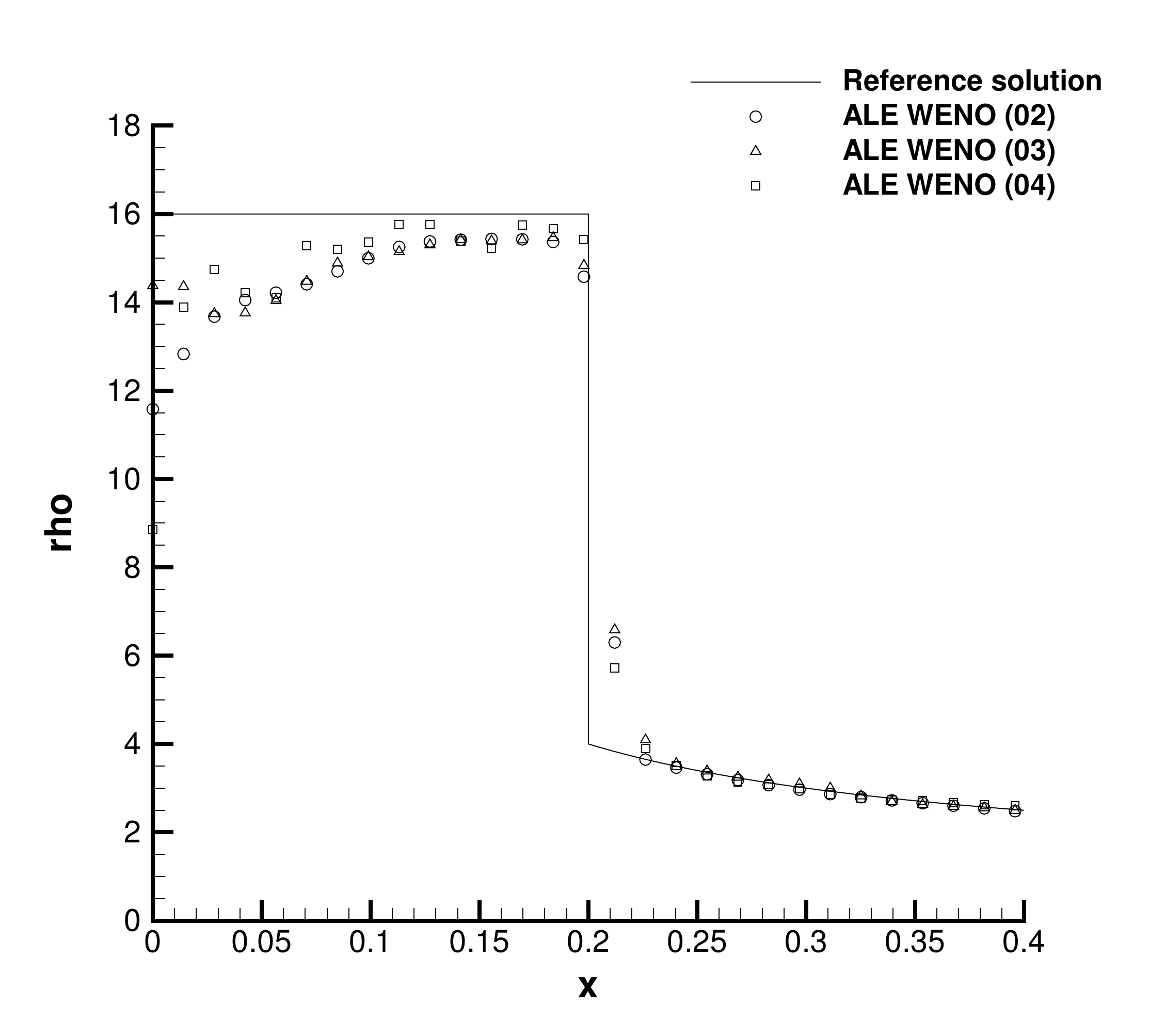} \\    
\end{tabular} 
\caption{Top: mesh configuration for the Noh problem at the initial time $t=0$ and at the final time $t_f=0.6$. Bottom: fourth order accurate density distribution at the final time and comparison between the exact solution (solid line) and three different high order accurate numerical results, i.e. $2^{nd}$, $3^{rd}$ and $4^{th}$ order ALE ADER-WENO finite volume schemes using the multi-dimensional HLL Riemann solver with $\textnormal{CFL}=0.9$.} 
\label{fig.Noh2D}
\end{center}
\end{figure}

In \cite{Noh} Noh introduced this test case, that involves a strong outward traveling shock wave produced by the compression of a zero pressure gas.  Initially the computational domain is square shaped, $\Omega(0)=[0;1.0]\times[0;1.0]$. The domain is discretized with a total number of elements of $N_E=5000$, obtained by splitting into triangles $50\times50$ square elements, as depicted in Figure \ref{fig.Noh2D}. A gas is initially assigned a unity density $\rho_0=1$ and an initial unit inward velocity $\mathbf{v}_in=(u,v)$, whose components are given by
\begin{equation}
u = -\frac{x}{r}, \qquad v = -\frac{y}{r}, 
\end{equation}  
where $r=\sqrt{x^2+y^2}$ is the general radial position. We set moving boundaries on the top and on the right boundaries, while no-slip wall boundary conditions have been imposed on the remaining sides. The ratio of specific heats is set to $\gamma=\frac{5}{3}$ and the initial pressure is $p=10^{-6}$ everywhere. According to \cite{Noh,Maire2009,Maire2009b}, we set a final time of $t_f=0.6$, hence the exact solution is given by an outward traveling shock wave located at radius $R=0.2$. The maximum density value occurs on the plateau behind the shock wave and it reaches the value of $\rho_f=16$, while the velocity of the shock wave is $v_{sh}=\frac{1}{3}$ along the radial direction. This is a well-known and very difficult test case, since the elements are highly deformed and distorted due to the very  strong shock wave, therefore we use the rezoning algorithm presented in Section \ref{sec.meshmotion} to recover a better mesh quality. Figure \ref{fig.Noh2D} shows  the initial and the final mesh configuration and a comparison between the exact solution and three high order accurate numerical results obtained with the ALE ADER-WENO finite volume schemes based on genuinely multi-dimensional HLL Riemann solvers presented in this paper. A Courant number of $\textnormal{CFL}=0.9$ has been used for all the numerical simulations and one can notice that the quality of the solution becomes the better as the order of accuracy of the scheme increases. 

\subsection{Numerical Convergence Study for the ideal MHD equations} 
\label{sec.conv.Rates-MHD}
\begin{figure}[!htbp]
\begin{center}
\begin{tabular}{cc} 
\includegraphics[width=0.47\textwidth]{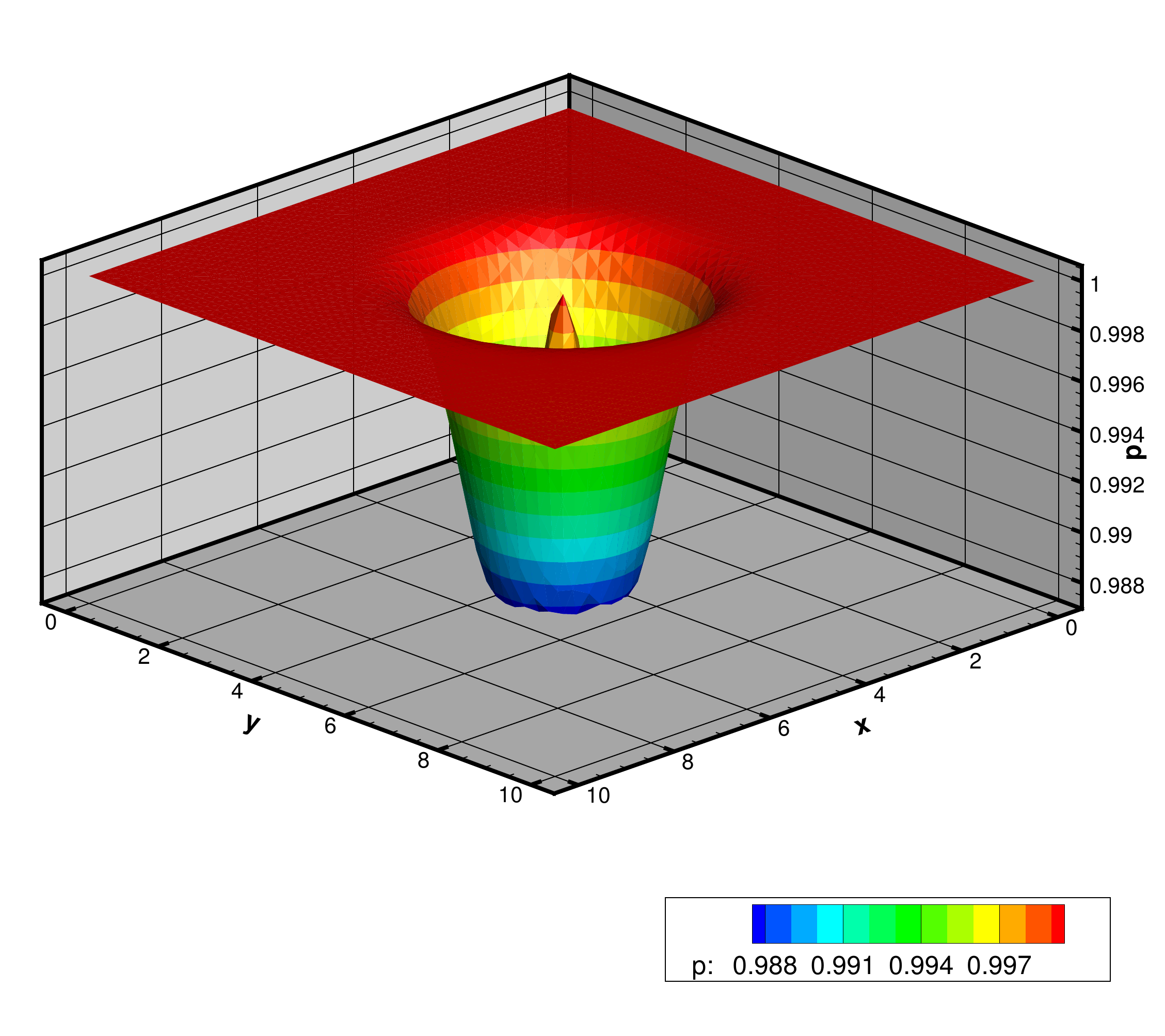}  &           
\includegraphics[width=0.47\textwidth]{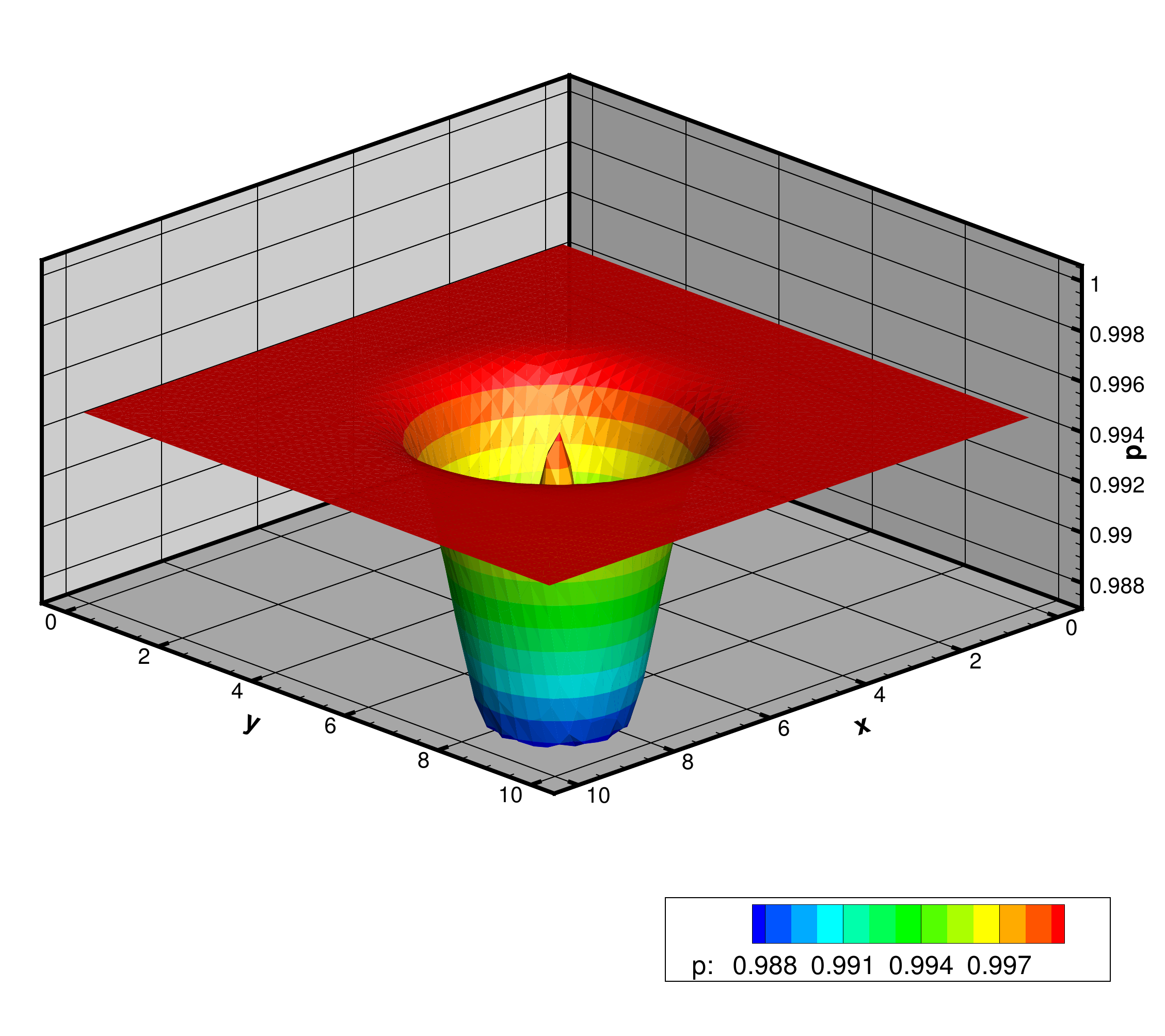} \\
\includegraphics[width=0.47\textwidth]{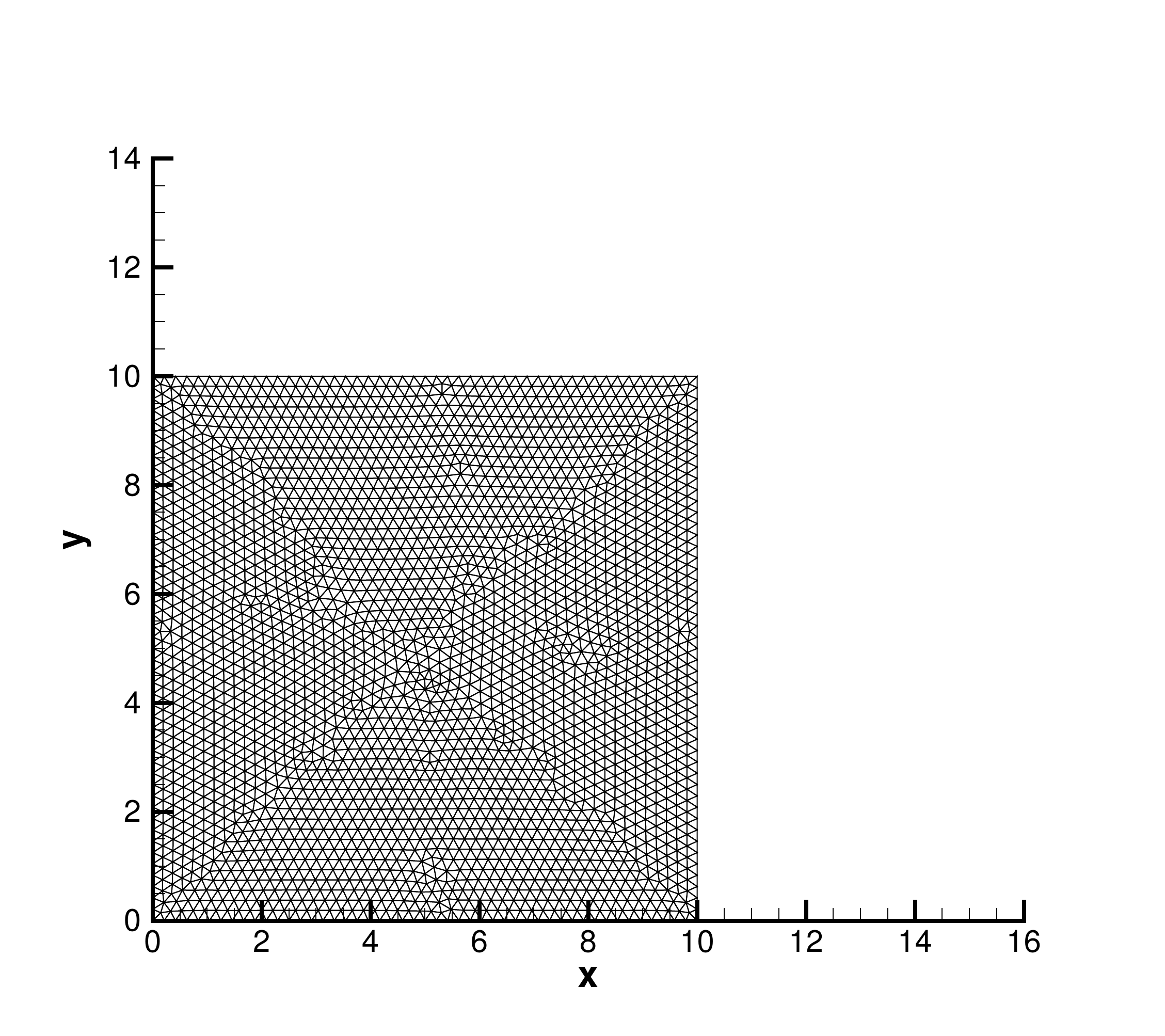}  &           
\includegraphics[width=0.47\textwidth]{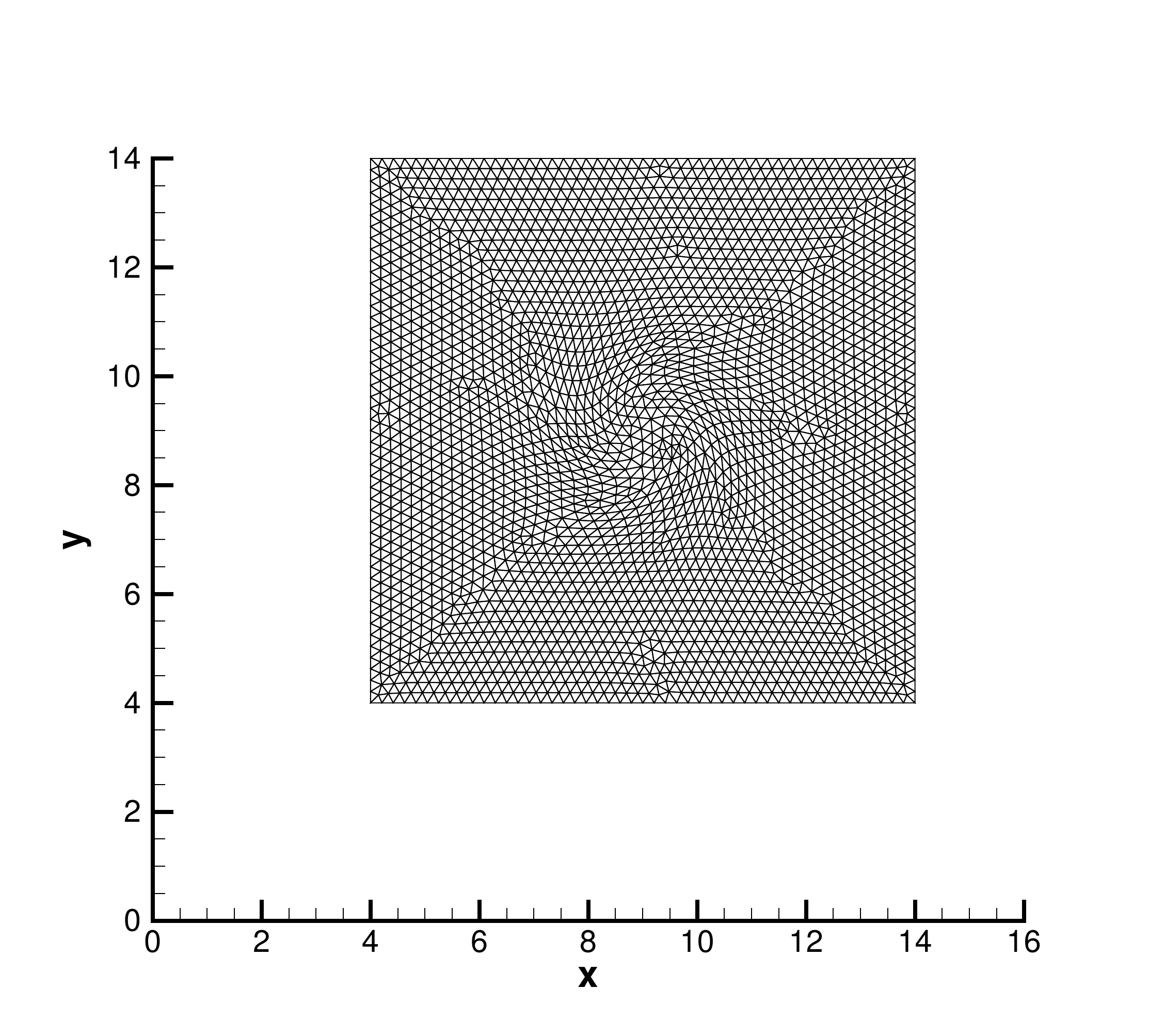} \\  
\end{tabular} 
\caption{Top: pressure distribution for the ideal MHD vortex problem at time $t=0.0$ and $t=4.0$. Bottom: mesh configuration at time $t=0.0$ and $t=4.0$. } 
\label{fig.MHDVortex}
\end{center}
\end{figure}

We use the convected smooth vortex test problem proposed by Balsara et al. \cite{Balsara2004} in order to carry out the numerical convergence studies for the ideal classical MHD equations. This test case is defined on a square computational domain $\Omega(0)=[0;10]\times[0;10]$ with periodic boundaries everywhere. As for the hydrodynamic isentropic vortex presented in Section \ref{sec.conv.Rates-Eul}, the initial condition is given by a linear superposition of a constant flow and some fluctuations in the velocity and magnetic fields, which read
\begin{equation}
\label{MHDVortIC}
(\rho, u, v, p, B_x, B_y, \Psi) = (1+\delta \rho, 1+\delta u, 1+\delta v, 1+\delta p, 1+\delta B_x, 1+\delta B_y, 0),
\end{equation} 
with the following perturbations:
\begin{eqnarray}
\left[\begin{array}{c} \delta u \\ \delta v \\ \delta p \\ \delta B_x \\ \delta B_y \end{array}\right] &=& \left[\begin{array}{c} \frac{\epsilon}{2\pi}e^{\frac{1}{2}(1-r^2)}(5-y) \\ \frac{\epsilon}{2\pi}e^{\frac{1}{2}(1-r^2)}(x-5) \\    \frac{1}{8\pi} \left(\frac{\mu}{2\pi}\right)^2(1-r^2) e^{(1-r^2)}-\frac{1}{2}\left(\frac{\epsilon}{2\pi}\right)^2 e^{(1-r^2)} \\ \frac{\mu}{2\pi}e^{\frac{1}{2}(1-r^2)}(5-y)
\\ \frac{\mu}{2\pi}e^{\frac{1}{2}(1-r^2)}(x-5)
\end{array}\right].  
\end{eqnarray}
According to \cite{Balsara2004}, we set the parameters $\epsilon=1$ and $\mu=\sqrt{4\pi}$ as well as the ratio of specific heats $\gamma=\frac{5}{3}$. The speed for the divergence cleaning is taken to be $c_h=2$ and the velocity $\v_c=(1,1)$ convects the vortex. The fluid motion of the vortex would lead to high element distortions and deformations, as clearly depicted in Figure \ref{fig.MHDVortex}, therefore the final time of the computation is $t_f=1.0$, because we do not want the rezoning step to be used for the convergence rate studies. 

The exact solution is given by the initial condition shifted in space by a factor $\mathbf{s}=(s_x,s_y)=\mathbf{v}\cdot t_f$. We perform the vortex problem on four successively refined meshes from first up to fourth order of accuracy and for each simulation we compute the error in $L_2$ norm, given by Eqn. 
\eqref{eqnL2error}. The multidimensional HLLC Riemann solver for the MHD equations has been used, see \cite{BalsaraMultiDRS}.

\begin{table}  
\caption{Numerical convergence results for the ideal MHD equations using the Lagrangian one-step WENO finite volume schemes with genuinely multidimensional HLL Riemann solvers presented in this article. The error norms refer to the variable $\rho$ (density) at time $t=1.0$ for first up to fourth order version of the scheme.} 
\begin{center} 
\begin{small}
\renewcommand{\arraystretch}{1.0}
\begin{tabular}{cccccc} 
\hline
  $h(\Omega(t_f))$ & $\epsilon_{L_2}$ & $\mathcal{O}(L_2)$ & $h(\Omega,t_f)$ & $\epsilon_{L_2}$ & $\mathcal{O}(L_2)$ \\
\hline
  \multicolumn{3}{c}{$\mathcal{O}1$} & \multicolumn{3}{c}{$\mathcal{O}2$} \\
\hline
3.26E-01 & 5.4032E-03 & -   & 3.25E-01  & 1.2393E-02 & -     \\ 
2.36E-01 & 4.7048E-03 & 0.4 & 2.46E-01  & 9.5840E-03 & 0.9   \\ 
1.63E-01 & 4.0697E-03 & 0.4 & 1.63E-01  & 5.7617E-03 & 1.2   \\ 
1.28E-01 & 3.5298E-03 & 0.6 & 1.28E-01  & 3.5875E-03 & 2.0   \\ 
\hline 
  \multicolumn{3}{c}{$\mathcal{O}3$} & \multicolumn{3}{c}{$\mathcal{O}4$}   \\
\hline
6.75E-01 & 1.6836E-02 & -   & 6.73E-01  & 1.9276E-02 & -      \\ 
3.25E-01 & 3.3009E-03 & 2.2 & 3.26E-01  & 1.0209E-03 & 4.1   \\ 
2.47E-01 & 1.0170E-03 & 4.3 & 2.47E-01  & 2.6494E-04 & 4.9   \\ 
1.63E-01 & 2.9097E-04 & 3.0 & 1.63E-01  & 5.1003E-05 & 4.0   \\ 
\hline 
\end{tabular}
\end{small}
\end{center}
\label{tab.convMHD}
\end{table}

\subsection{The MHD Rotor Problem} 
\label{sec.MHDRotor}
\begin{figure}[!htbp]
\begin{center}
\begin{tabular}{cc} 
\includegraphics[width=0.47\textwidth]{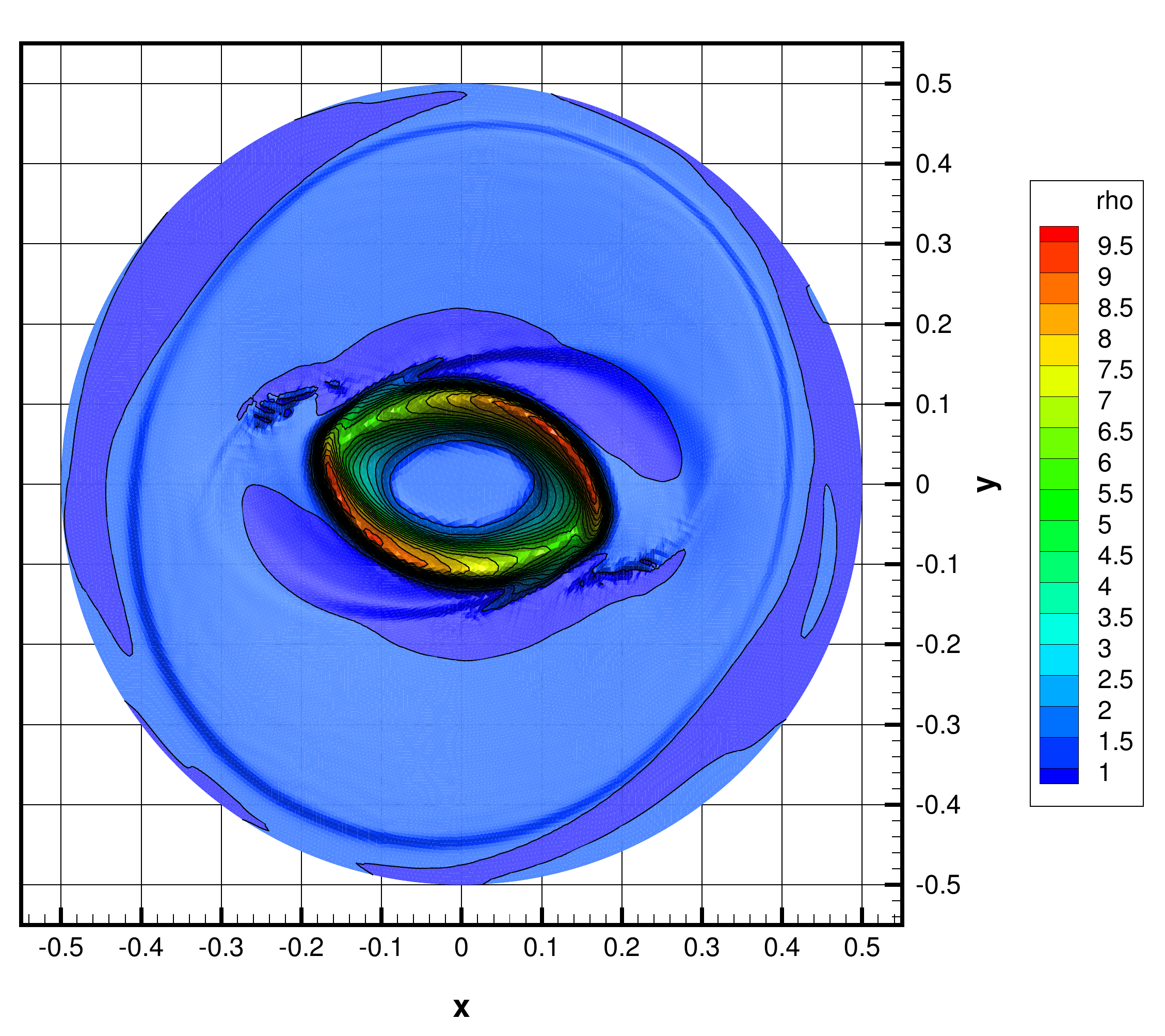}  &           
\includegraphics[width=0.47\textwidth]{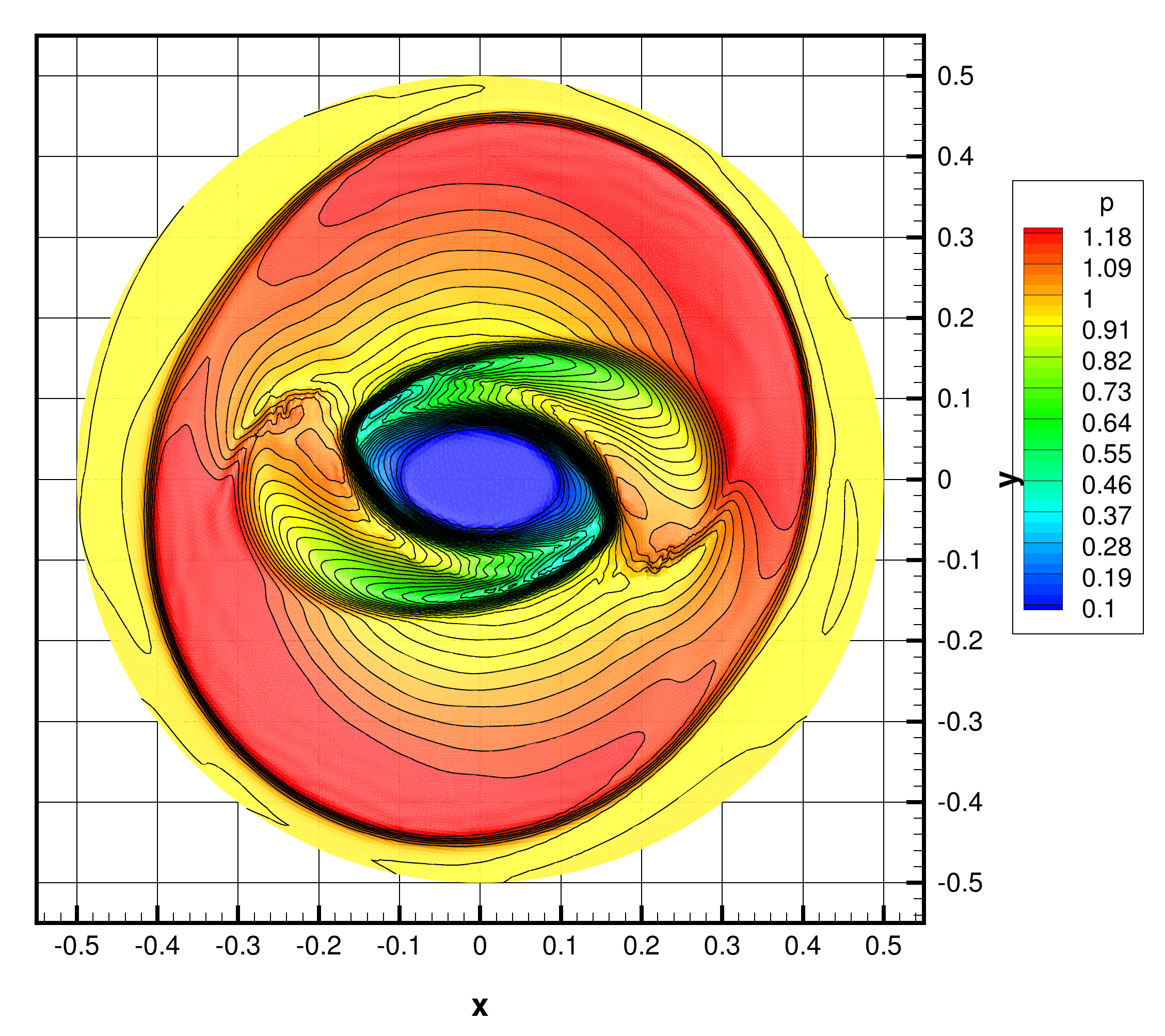} \\  
\end{tabular} 
\caption{Third order numerical results for the ideal MHD rotor problem: density and pressure at time $t=0.25$.} 
\label{fig.MHDRotor}
\end{center}
\end{figure}

In \cite{BalsaraSpicer1999} Balsara and Spicer solve the ideal MHD rotor problem that consists in a high density fluid that is rotating around the center of a circular computational domain of radius $R=0.5$, i.e. $\Omega(0)=\left\{ \mathbf{x} : \left\| \mathbf{x} \right\| \leq R \right\}$. The high density region is  delimited by a circle of radius $R_i=0.1$, splitting the domain in the internal region, which is filled by the moving fluid, and the external region, characterized by a low density fluid at rest. According to \cite{BalsaraSpicer1999}, at $r=1$ the toroidal velocity is $v_t=(\omega \cdot R)=1$, since the angular velocity $\omega$ of the rotor is assumed to be constant. The initial pressure $p=1$ is constant in the whole domain as well as the initial magnetic field $\mathbf{B}=(2.5,0,0)^T$, while the initial density distribution is $\rho=10$ for $0 \leq r \leq R_i$ and $\rho=1$ elsewhere. Furthermore, we use a linear taper to smear out the initial discontinuity occurring at radius $R_i=0.1$ involving both velocity and density. The taper is applied for $0.1 \leq r \leq 0.13$, so that at radius $r=0.1$ and $r=0.13$ density and velocity match exactly the values of the inner and the outer region, respectively. Any further detail regarding the taper can be found in \cite{BalsaraSpicer1999}. The divergence cleaning speed is taken to be $c_h=2$ and the ratio of specific heats is $\gamma=1.4$. We impose transmissive boundary conditions at the external boundary and we use a third order ALE WENO scheme with the the multi-dimensional HLL flux for the MHD equations \cite{BalsaraMultiDRS} on a computational grid with a characteristic mesh size of $h=1/200$. Figure \ref{fig.MHDRotor} shows the numerical results at the final time $t_f=0.25$ obtained with 
$\textnormal{CFL}=0.95$. The Alfv\'en waves produced by the rotor tend to diminish the angular momentum of the rotor as the simulation goes on and a good agreement with the solution presented in \cite{BalsaraSpicer1999} can be noticed.

\subsection{The MHD Blast Wave Problem} 
\label{sec.MHDBlast}
\begin{figure}[!htbp]
\begin{center}
\begin{tabular}{cc} 
\includegraphics[width=0.47\textwidth]{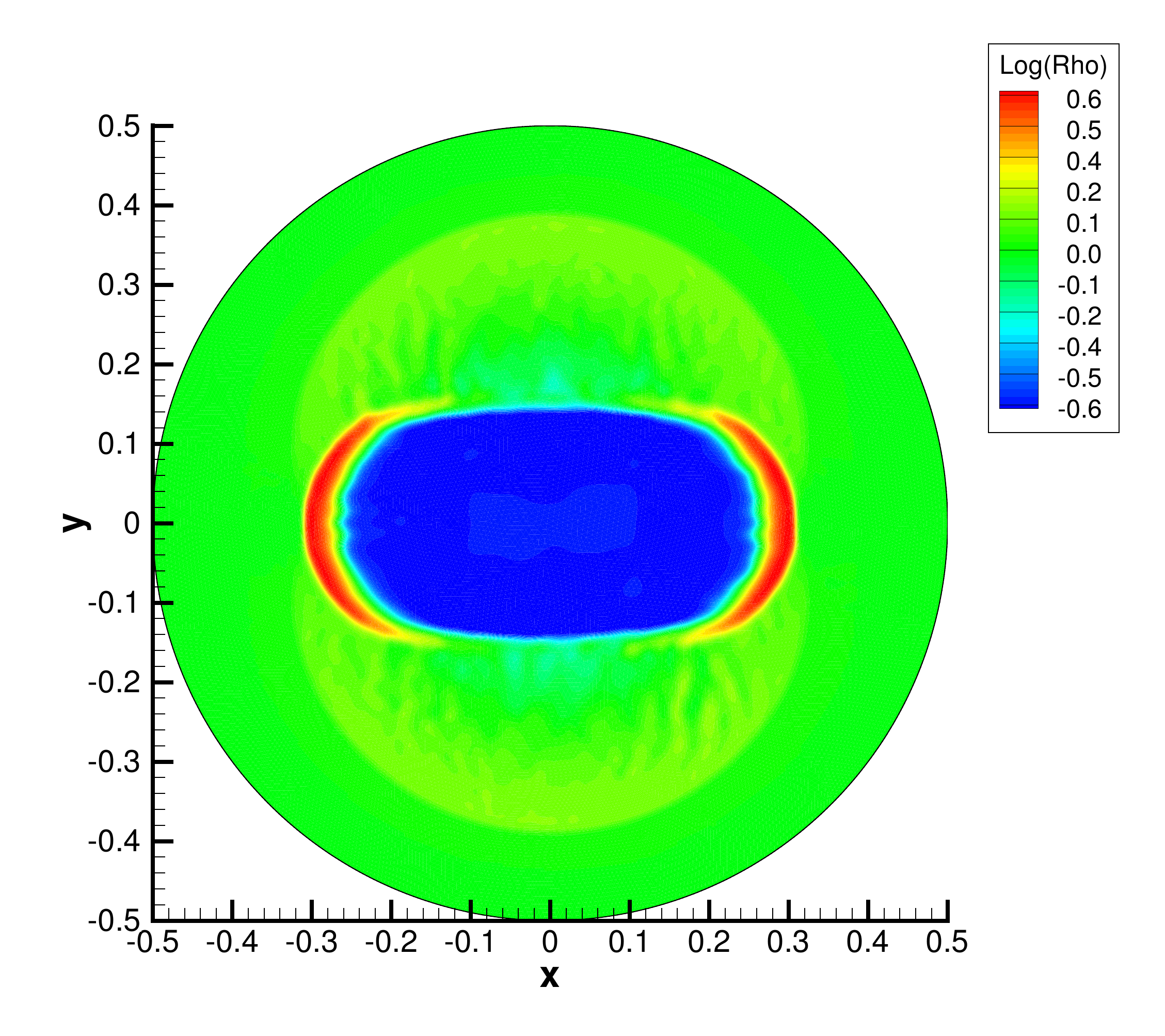}  &           
\includegraphics[width=0.47\textwidth]{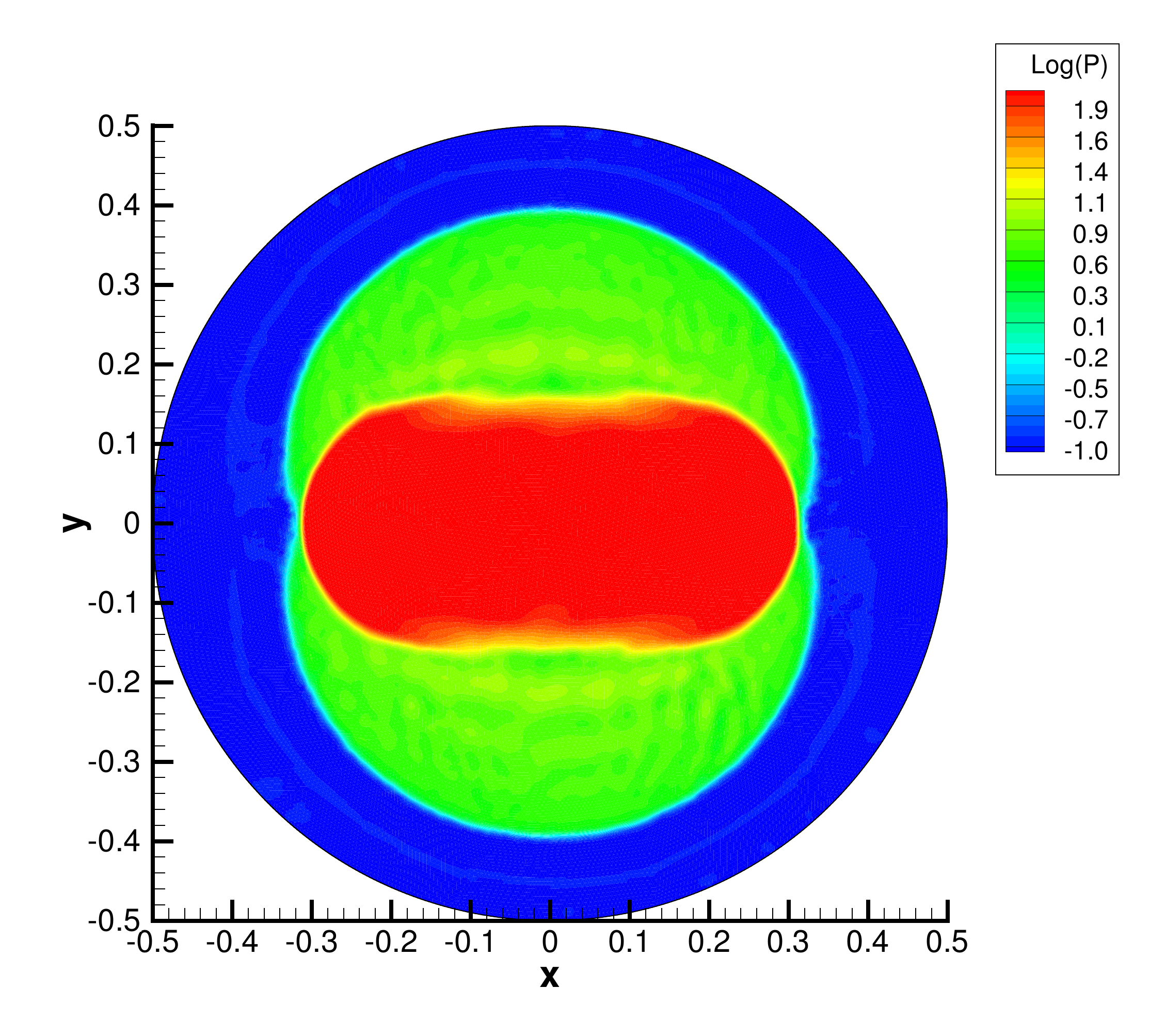} \\   
\end{tabular} 
\caption{Numerical results for the MHD blast wave problem at time $t=0.01$. Left: logarithm (base 10) of the density. Right: logarithm (base 10) of the pressure.} 
\label{fig.Blast}
\end{center}
\end{figure}

The MHD blast wave problem is characterized by a strong circular fast magnetosonic shock wave traveling from the center towards the boundary of a circular computational domain $\Omega(0)$ of radius $R_0=1.0$ with the generic radial position defined as usual by $r=\sqrt{x^2+y^2}$. The initially constant magnetic field in the $x-$direction limits the expansion of the shock wave along the $y-$direction, hence stretching and propagating the wave towards the $x-$oriented boundaries. This well known test case was first proposed in \cite{BalsaraSpicer1999b} and the initial condition $\U(\x,0)=(\rho,u,v,p,B_x,B_y,\Psi)$ is given by two different states 
\begin{equation}
  \U(\x,0) = \left\{ \begin{array}{ccc} \U_i & \textnormal{ if } & r \leq R, \\ 
                                        \U_o & \textnormal{ if } & r > R,        
                      \end{array}  \right. 
\end{equation}
where $\Q_i$ denotes the \textit{inner} state, which is bounded by a circle of radius $R=0.1$, and $\Q_e$ is the \textit{outer} state that covers the remaining part of the domain. The initial states in primitive variables read 
\begin{equation}
\U_i = \left( 1.0, 0.0, 0.0, 0.1, 70, 0.0, 0.0 \right), \qquad \U_e = \left( 1.0, 0.0, 0.0, 1000, 70, 0.0, 0.0 \right).
\label{eq.BlastIC}
\end{equation}
The final time of the simulation is taken to be $t_f=10^{-2}$ and the ratio of specific heats is $\gamma=1.4$. We run the simulation on the same mesh used for the ideal MHD rotor problem, with a characteristic mesh size of $h=1/200$ and transmissive boundary conditions everywhere. Since strong deformations occur mostly where the shock front is traveling, we use the rezoning algorithm presented in Section \ref{sec.meshmotion} to avoid mesh tangling or element overlapping. The numerical results depicted in Figure \ref{fig.Blast} have been computed using the third order accurate version of the ALE ADER-WENO finite volume scheme with the multidimensional HLL flux for the MHD equations \cite{BalsaraMultiDRS} and $\textnormal{CFL}=0.95$. We show the logarithm of density and pressure, and the solution looks very similar to the results given in \cite{Balsara2004}.

\section{Conclusions}
\label{sec.concl}
We have presented the first high-order unstructured ADER-WENO ALE finite volume schemes based on genuinely multidimensional HLL Riemann solvers 
\cite{BalsaraMultiDRS}. We applied the new algorithm to both hydrodynamics and magnetohydrodynamics and convergence studies up to fourth order of accuracy 
in space and time have been carried out. The multidimensionality in the Riemann solver allowed the numerical method to run with a less severe CFL condition 
on the timestep, namely taking $\textnormal{CFL}=0.95$ in two space dimensions, instead of setting it to the usual limit of $\textnormal{CFL} \leq 0.5$ 
typical for unsplit Godunov schemes in two space dimensions using conventional one-dimensional Riemann solvers. For those test cases that involve strong 
shock or shear waves, a rezoning algorithm has been used in order to recover a better mesh quality. 

A possible research field in the future would be the extension of the genuinely multidimensional HLL Riemann solvers to three-dimensional tetrahedral meshes, 
first in the  Eulerian framework and later also for Lagrangian schemes on moving meshes.  
The multidimensional evolution and integration of the gradients and second derivatives could also be included in the multidimensional Riemann solver, in such 
a way that a quadrature free approach for the flux evaluation can be developed, which might result in further improvements concerning computational efficiency.  

\section*{Acknowledgments}
MD and WB have been financed by the European Research Council (ERC) under the European Union's Seventh Framework 
Programme (FP7/2007-2013) with the research project \textit{STiMulUs}, ERC Grant agreement no. 278267. 

DSB acknowledges support via NSF grants NSF-AST-1009091 and NSF-ACI-1307369. DSB also acknowledges support via NASA grants from the Fermi program 
as well as NASA-NNX 12A088G. 

The authors acknowledge PRACE for awarding access to the SuperMUC supercomputer of the Leibniz Rechenzentrum (LRZ) based in Munich, Germany. 

\bibliography{LagrangeMDRS}

\begin{thebibliography}{100}

\bibitem{AbgrallmultiD1}
R.~Abgrall.
\newblock {Approximation du probl\`eme de Riemann vraiment multidimensionnel
  des \'equations d'Euler par une m\'ethode de type Roe, I: La
  lin\'earisation}.
\newblock {\em C.R. Acad. Sci. Ser. I}, 319:499 -- 504, 1994.

\bibitem{AbgrallmultiD2}
R.~Abgrall.
\newblock {Approximation du probl\`eme de Riemann vraiment multidimensionnel
  des \'equations d'Euler par une m\'ethode de type Roe, II: Solution du
  probl\`eme de Riemann approch\'e}.
\newblock {\em C.R. Acad. Sci. Ser. I}, 319:625 -- 629, 1994.

\bibitem{AbgrallmultiD3}
R.~Abgrall.
\newblock On essentially non-oscillatory schemes on unstructured meshes :
  Analysis and implementation.
\newblock {\em Journal of Computational Physics}, 114:45 -- 58, 1994.

\bibitem{AboiyarIske}
T.~Aboiyar, E.H. Georgoulis, and A.~Iske.
\newblock {Adaptive ADER Methods Using Kernel-Based Polyharmonic Spline WENO
  Reconstruction}.
\newblock {\em SIAM Journal on Scientific Computing}, 32:3251--3277, 2010.

\bibitem{BalsaraRMHD}
D.~Balsara.
\newblock Total variation diminishing scheme for relativistic
  magnetohydrodynamics.
\newblock {\em The Astrophysical Journal Supplement Series}, 132:83--101, 2001.

\bibitem{Balsara2004}
D.~Balsara.
\newblock Second-order accurate schemes for magnetohydrodynamics with
  divergence-free reconstruction.
\newblock {\em The Astrophysical Journal Supplement Series}, 151:149--184,
  2004.

\bibitem{BalsaraSpicer1999b}
D.~Balsara and D.~Spicer.
\newblock Maintaining pressure positivity in magneto\-hydrodynamic simulations.
\newblock {\em Journal of Computational Physics}, 148:133–--148, 1999.

\bibitem{BalsaraSpicer1999}
D.~Balsara and D.~Spicer.
\newblock A staggered mesh algorithm using high order godunov fluxes to ensure
  solenoidal magnetic fields in magnetohydrodynamic simulations.
\newblock {\em Journal of Computational Physics}, 149:270--292, 1999.

\bibitem{BalsaraAMRdivfree2001}
D.S. Balsara.
\newblock Divergence-free adaptive mesh refinement for magnetohydrodynamics.
\newblock {\em Journal of Computational Physics}, 174:614--648, 2001.

\bibitem{balsarahlle2d}
D.S. Balsara.
\newblock {Multidimensional HLLE Riemann solver: Application to Euler and
  magnetohydrodynamic flows}.
\newblock {\em Journal of Computational Physics}, 229:1970--1993, 2010.

\bibitem{balsarahllc2d}
D.S. Balsara.
\newblock {A two-dimensional HLLC Riemann solver for conservation laws:
  Application to Euler and magnetohydrodynamic flows}.
\newblock {\em Journal of Computational Physics}, 231:7476--7503, 2012.

\bibitem{BalsaraPositivity2012}
D.S. Balsara.
\newblock Self-adjusting, positivity preserving high order schemes for
  hydrodynamics and magnetohydrodynamics.
\newblock {\em Journal of Computational Physics}, 231:7504--7517, 2012.

\bibitem{BalsaraMultiDRS}
D.S. Balsara, M.~Dumbser, and R.~Abgrall.
\newblock {Multidimensional HLLC Riemann Solver for Unstructured Meshes}.
\newblock {\em Journal of Computational Physics}.
\newblock submitted to.

\bibitem{BalsaraRumpf}
D.S. Balsara, T.~Rumpf, M.~Dumbser, and C.D. Munz.
\newblock Efficient, high accuracy {ADER-WENO} schemes for hydrodynamics and
  divergence-free magnetohydrodynamics.
\newblock {\em Journal of Computational Physics}, 228:2480--2516, 2009.

\bibitem{BalsaraShu2000}
D.S. Balsara and C.W. Shu.
\newblock Monotonicity preserving weighted essentially non-oscillatory schemes
  with increasingly high order of accuracy.
\newblock {\em Journal of Computational Physics}, 160:405--452, 2000.

\bibitem{StencilRec1990}
{T.J.} Barth and {P.O.} Frederickson.
\newblock Higher order solution of the euler equations on unstructured grids
  using quadratic reconstruction.
\newblock {\em 28th Aerospace Sciences Meeting}, pages AIAA paper no. 90--0013,
  January 1990.

\bibitem{Artzi}
M.~Ben-Artzi and J.~Falcovitz.
\newblock A second-order godunov-type scheme for compressible fluid dynamics.
\newblock {\em Journal of Computational Physics}, 55:1--32, 1984.

\bibitem{ShashkovRemap2}
M.~Berndt, J.~Breil, S.~Galera, M.~Kucharik, P.H. Maire, and M.J. Shashkov.
\newblock {Two-step hybrid conservative remapping for multimaterial arbitrary
  Lagrangian-Eulerian methods}.
\newblock {\em Journal of Computational Physics}, 230:6664--6687, 2011.

\bibitem{WAF}
S.J. Billett and E.F. Toro.
\newblock On waf-type schemes for multidimensional hyperbolic conservation
  laws.
\newblock {\em Journal of Computational Physics}, 130:1 -- 24, 1997.

\bibitem{ShashkovRemap1}
P.~Bochev, D.~Ridzal, and M.J. Shashkov.
\newblock Fast optimization-based conservative remap of scalar fields through
  aggregate mass transfer.
\newblock {\em Journal of Computational Physics}, 246:37--57, 2013.

\bibitem{FCT2}
D.L. Book, J.P. Boris, and K.~Hain.
\newblock {Flux-corrected transport II: Generalizations of the method }.
\newblock {\em Journal of Computational Physics}, 18:248--283, 1975.

\bibitem{FCT1}
J.P. Boris and D.L. Book.
\newblock {Flux-corrected transport. I. SHASTA, a fluid transport algorithm
  that works}.
\newblock {\em Journal of Computational Physics}, 11:38--69, 1973.

\bibitem{BoscheriDumbserLag}
W.~Boscheri and M.~Dumbser.
\newblock {Arbitrary--Lagrangian--Eulerian One--Step WENO Finite Volume Schemes
  on Unstructured Triangular Meshes}.
\newblock {\em Communications in Computational Physics}, 14:1174--1206, 2013.

\bibitem{LagNS}
W.~Boscheri, M.~Dumbser, and D.S. Balsara.
\newblock {High Order Lagrangian ADER-WENO Schemes on Unstructured Meshes --
  Application of Several Node Solvers to Hydrodynamics and
  Magnetohydrodynamics}.
\newblock {\em International Journal of Numerical Methods in Fluids}.
\newblock submitted to.

\bibitem{Raviart.GRP.2}
A.~Bourgeade, P.~LeFloch, and P.A. Raviart.
\newblock An asymptotic expansion for the solution of the generalized {Riemann}
  problem. {Part II}: application to the gas dynamics equations.
\newblock {\em Annales de l'institut Henri Poincar\'e (C) Analyse non
  lin\'eaire}, 6:437--480, 1989.

\bibitem{ShashkovMultiMat3}
J.~Breil, T.~Harribey, P.H. Maire, and M.J. Shashkov.
\newblock {A multi-material ReALE method with MOF interface reconstruction}.
\newblock {\em Computers and Fluids}, 83:115--125, 2013.

\bibitem{Despres2009}
G.~Carr\'e, S.~Del Pino, B.~Despr\'es, and E.~Labourasse.
\newblock A cell-centered lagrangian hydrodynamics scheme on general
  unstructured meshes in arbitrary dimension.
\newblock {\em Journal of Computational Physics}, 228:5160--5183, 2009.

\bibitem{Feistauer4}
J.~Cesenek, M.~Feistauer, J.~Horacek, V.~Kucera, and J.~Prokopova.
\newblock Simulation of compressible viscous flow in time-dependent domains.
\newblock {\em Applied Mathematics and Computation}, 219:7139--7150, 2013.

\bibitem{chengshu1}
J.~Cheng and C.W. Shu.
\newblock {A high order ENO conservative Lagrangian type scheme for the
  compressible Euler equations}.
\newblock {\em Journal of Computational Physics}, 227:1567--1596, 2007.

\bibitem{chengshu3}
J.~Cheng and C.W. Shu.
\newblock {A cell-centered Lagrangian scheme with the preservation of symmetry
  and conservation properties for compressible fluid flows in two-dimensional
  cylindrical geometry}.
\newblock {\em Journal of Computational Physics}, 229:7191--7206, 2010.

\bibitem{chengshu4}
J.~Cheng and C.W. Shu.
\newblock {Improvement on spherical symmetry in two-dimensional cylindrical
  coordinates for a class of control volume Lagrangian schemes}.
\newblock {\em Communications in Computational Physics}, 11:1144--1168, 2012.

\bibitem{MOOD}
S.~Clain, S.~Diot, and R.~Loub\`ere.
\newblock {A high--order finite volume method for systems of conservation laws
  -- Multi--dimensional Optimal Order Detection (MOOD)}.
\newblock {\em Journal of Computational Physics}, 230:4028--4050, 2011.

\bibitem{CBS-book}
B.~Cockburn, G.~E. Karniadakis, and {C.W.} Shu.
\newblock {\em Discontinuous {Galerkin} Methods}.
\newblock Lecture Notes in Computational Science and Engineering. Springer,
  2000.

\bibitem{ColellaRS}
P.~Colella.
\newblock A direct eulerian muscl scheme for gas dynamics.
\newblock {\em SIAM J. Sci. Statist. Comput.}, 6:104 -- 117, 1985.

\bibitem{Colella1990}
P.~Colella.
\newblock Multidimensional upwind methods for hyperbolic conservation laws.
\newblock {\em Journal of Computational Physics}, 87:171 -- 200, 1990.

\bibitem{Dedneretal}
A.~Dedner, F.~Kemm, D.~Kr\"oner, C.-D. Munz, T.~Schnitzer, and M.~Wesenberg.
\newblock Hyperbolic divergence cleaning for the {MHD} equations.
\newblock {\em Journal of Computational Physics}, 175:645--673, 2002.

\bibitem{DepresMazeran2003}
B.~Despr\'es and C.~Mazeran.
\newblock Symmetrization of lagrangian gas dynamic in dimension two and
  multimdimensional solvers.
\newblock {\em C.R. Mecanique}, 331:475--480, 2003.

\bibitem{Despres2005}
B.~Despr\'es and C.~Mazeran.
\newblock Lagrangian gas dynamics in two-dimensions and lagrangian systems.
\newblock {\em Archive for Rational Mechanics and Analysis}, 178:327--372,
  2005.

\bibitem{Feistauer1}
L.~Dubcova, M.~Feistauer, J.~Horacek, and P.~Svacek.
\newblock {Numerical simulation of interaction between turbulent flow and a
  vibrating airfoil}.
\newblock {\em Computing and Visualization in Science}, 12:207--225, 2009.

\bibitem{Dubiner}
M.~Dubiner.
\newblock Spectral methods on triangles and other domains.
\newblock {\em Journal of Scientific Computing}, 6:345--390, 1991.

\bibitem{SaltzmanOrg}
{J.K.} Dukovicz and B.~Meltz.
\newblock Vorticity errors in multidimensional lagrangian codes.
\newblock {\em Journal of Computational Physics}, 99:115 -- 134, 1992.

\bibitem{DukowiczRS}
J.K. Dukowicz.
\newblock A general non-iterative riemann solver for godunov's method.
\newblock {\em Journal of Computational Physics}, 61:119 -- 137, 1984.

\bibitem{Dumbser20088209}
M.~Dumbser, {D.S.} Balsara, E.F. Toro, and {C.-D.} Munz.
\newblock A unified framework for the construction of one-step finite volume
  and discontinuous galerkin schemes on unstructured meshes.
\newblock {\em Journal of Computational Physics}, 227:8209 -- 8253, 2008.

\bibitem{DumbserBoscheriLagNC}
M.~Dumbser and W.~Boscheri.
\newblock High-order unstructured lagrangian one-step weno finite volume
  schemes for non-conservative hyperbolic systems: Applications to compressible
  multi-phase flows.
\newblock {\em Computers and Fluids}, 86:405 -- 432, 2013.

\bibitem{DumbserEnauxToro}
M.~Dumbser, C.~Enaux, and E.F. Toro.
\newblock Finite volume schemes of very high order of accuracy for stiff
  hyperbolic balance laws.
\newblock {\em Journal of Computational Physics}, 227:3971--4001, 2008.

\bibitem{DumbserKaeser06b}
M.~Dumbser and M.~K\"aser.
\newblock Arbitrary high order non-oscillatory finite volume schemes on
  unstructured meshes for linear hyperbolic systems.
\newblock {\em Journal of Computational Physics}, 221:693--723, 2007.

\bibitem{Dumbser2007204}
M.~Dumbser, M.~K\"aser, {V.A.} Titarev, and {E.F.} Toro.
\newblock Quadrature-free non-oscillatory finite volume schemes on unstructured
  meshes for nonlinear hyperbolic systems.
\newblock {\em Journal of Computational Physics}, 226:204 -- 243, 2007.

\bibitem{OsherUniversal}
M.~Dumbser and E.~F. Toro.
\newblock On universal {Osher}--type schemes for general nonlinear hyperbolic
  conservation laws.
\newblock {\em Communications in Computational Physics}, 10:635--671, 2011.

\bibitem{OsherNC}
M.~Dumbser and E.~F. Toro.
\newblock A simple extension of the {Osher} {Riemann} solver to
  non-conservative hyperbolic systems.
\newblock {\em Journal of Scientific Computing}, 48:70--88, 2011.

\bibitem{Lagrange1D}
M.~Dumbser, A.~Uuriintsetseg, and O.~Zanotti.
\newblock {On Arbitrary--Lagrangian--Eulerian One--Step WENO Schemes for Stiff
  Hyperbolic Balance Laws}.
\newblock {\em Communications in Computational Physics}, 14:301--327, 2013.

\bibitem{EinfeldtHLL}
B.~Einfeldt.
\newblock On godunov-type methods for gas dynamics.
\newblock {\em SIAM J. Numer. Anal.}, 25:294 -- 318, 1988.

\bibitem{munz91}
B.~Einfeldt, C.~D. Munz, P.~L. Roe, and B.~Sj\"ogreen.
\newblock On godunov-type methods near low densities.
\newblock {\em Journal of Computational Physics}, 92:273--295, 1991.

\bibitem{Feistauer3}
M.~Feistauer, J.~Horacek, M.~Ruzicka, and P.~Svacek.
\newblock Numerical analysis of flow-induced nonlinear vibrations of an airfoil
  with three degrees of freedom.
\newblock {\em Computers and Fluids}, 49:110--127, 2011.

\bibitem{Feistauer2}
M.~Feistauer, V.~Kucera, J.~Prokopova, and J.~Horacek.
\newblock {The ALE discontinuous Galerkin method for the simulatio of air flow
  through pulsating human vocal folds}.
\newblock {\em AIP Conference Proceedings}, 1281:83--86, 2010.

\bibitem{Fey1}
M.~Fey.
\newblock Multidimensional upwinding 1. the method of transport for solving the
  euler equations.
\newblock {\em Journal of Computational Physics}, 143:159 -- 180, 1998.

\bibitem{Fey2}
M.~Fey.
\newblock Multidimensional upwinding 2. decomposition of the euler equation
  into advection equation.
\newblock {\em Journal of Computational Physics}, 143:159 -- 199, 1998.

\bibitem{Raviart.GRP.1}
P.~Le Floch and P.A. Raviart.
\newblock An asymptotic expansion for the solution of the generalized {Riemann}
  problem. {Part I}: General theory.
\newblock {\em Annales de l'institut Henri Poincar\'e (C) Analyse non
  lin\'eaire}, 5:179--207, 1988.

\bibitem{ShashkovMultiMat1}
M.M. Francois, M.J. Shashkov, T.O. Masser, and E.D. Dendy.
\newblock {A comparative study of multimaterial Lagrangian and Eulerian methods
  with pressure relaxation.}
\newblock {\em Computers and Fluids}, 83:126--136, 2013.

\bibitem{friedrich}
O.~Friedrich.
\newblock Weighted essentially non-oscillatory schemes for the interpolation of
  mean values on unstructured grids.
\newblock {\em Journal of Computational Physics}, 144:194--212, 1998.

\bibitem{MaireRezoning}
S.~Galera, P.H. Maire, and J.~Breil.
\newblock A two-dimensional unstructured cell-centered multi-material ale
  scheme using vof interface reconstruction.
\newblock {\em Journal of Computational Physics}, 229:5755--5787, 2010.

\bibitem{GodunovRS}
S.~K. Godunov.
\newblock Finite difference methods for the computation of discontinuous
  solutions of the equations of fluid dynamics.
\newblock {\em Mathematics of the USSR}, 47:271--306, 1959.

\bibitem{HartenHLL}
A.~Harten, P.~D. Lax, and B.~van Leer.
\newblock On upstream differencing and godunov-type schemes for hyperbolic
  conservation laws.
\newblock {\em SIAM Rev}, 25:289 -- 315, 1983.

\bibitem{HidalgoDumbser}
A.~Hidalgo and M.~Dumbser.
\newblock {ADER} schemes for nonlinear systems of stiff
  advection–diffusion–reaction equations.
\newblock {\em Journal of Scientific Computing}, 48:173--189, 2011.

\bibitem{Hirt1974}
C.~Hirt, A.~Amsden, and J.~Cook.
\newblock An arbitrary lagrangian–eulerian computing method for all flow
  speeds.
\newblock {\em Journal of Computational Physics}, 14:227–253, 1974.

\bibitem{HuShuVortex1999}
C.~Hu and {C.W.} Shu.
\newblock A high-order weno finite difference scheme for the equations of ideal
  magnetohydrodynamics.
\newblock {\em Journal of Computational Physics}, 150:561 -- 594, 1999.

\bibitem{JiangShu1996}
{G.S.} Jiang and {C.W.} Shu.
\newblock Efficient implementation of weighted eno schemes.
\newblock {\em Journal of Computational Physics}, pages 202 -- 228, 1996.

\bibitem{SedovExact}
{J.R.} Kamm and {F.X.} Timmes.
\newblock On efficient generation of numerically robust sedov solutions.
\newblock {\em Technical Report LA-UR-07-2849,LANL}, 2007.

\bibitem{orth-basis}
G.~E. Karniadakis and S.~J. Sherwin.
\newblock {\em Spectral/hp Element Methods in CFD}.
\newblock Oxford University Press, 1999.

\bibitem{KaeserIske2005}
M.~K\"aser and A.~Iske.
\newblock Ader schemes on adaptive triangular meshes for scalar conservation
  laws.
\newblock {\em Journal of Computational Physics}, 205:486 -- 508, 2005.

\bibitem{Kidder1976}
{R.E.} Kidder.
\newblock Laser-driven compression of hollow shells: power requirements and
  stability limitations.
\newblock {\em Nucl. Fus.}, 1:3 -- 14, 1976.

\bibitem{KnuppRezoning}
{P.M.} Knupp.
\newblock Achieving finite element mesh quality via optimization of the
  jacobian matrix norm and associated quantities. part ii -- a framework for
  volume mesh optimization and the condition number of the jacobian matrix.
\newblock {\em Int. J. Numer. Meth. Engng.}, 48:1165 -- 1185, 2000.

\bibitem{ShashkovRemap5}
M.~Kucharik, J.~Breil, S.~Galera, P.H. Maire, M.~Berndt, and M.J. Shashkov.
\newblock {Hybrid remap for multi-material ALE}.
\newblock {\em Computers and Fluids}, 46:293--297, 2011.

\bibitem{ShashkovRemap3}
M.~Kucharik and M.J. Shashkov.
\newblock {One-step hybrid remapping algorithm for multi-material arbitrary
  Lagrangian-Eulerian methods}.
\newblock {\em Journal of Computational Physics}, 231:2851--2864, 2012.

\bibitem{ShashkovRemap4}
R.~Liska, M.J. Shashkov~P. V\'achal, and B.~Wendroff.
\newblock {Synchronized flux corrected remapping for ALE methods}.
\newblock {\em Computers and Fluids}, 46:312--317, 2011.

\bibitem{chengshu2}
W.~Liu, J.~Cheng, and C.W. Shu.
\newblock {High order conservative Lagrangian schemes with Lax–Wendroff type
  time discretization for the compressible Euler equations}.
\newblock {\em Journal of Computational Physics}, 228:8872--8891, 2009.

\bibitem{LoubereMaireShashkov}
R.~Loub\`ere, P.H. Maire, and M.J. Shashkov.
\newblock {ReALE: A Reconnection Arbitrary-Lagrangian-Eulerian method in
  cylindrical geometry}.
\newblock {\em Computers and Fluids}, 46:59--69, 2011.

\bibitem{EvolutionGalerkin2}
M.~Lukacova-Medvidova, K.W. Morton, and G.~Warnecke.
\newblock {Finite volume evolution Galerkin methods for Euler equations of gas
  dynamics}.
\newblock {\em International Journal of Numerical Methods in Fluids},
  40:425--434, 2002.

\bibitem{EvolutionGalerkin1}
M.~Lukacova-Medvidova, K.W. Morton, and G.~Warnecke.
\newblock {Finite volume evolution Galerkin methods for hyperbolic systems}.
\newblock {\em SIAM Journal on Scientific Computing}, 26:1--30, 2005.

\bibitem{Maire2009}
{P.H.} Maire.
\newblock A high-order cell-centered lagrangian scheme for two-dimensional
  compressible fluid flows on unstructured meshes.
\newblock {\em Journal of Computational Physics}, 228:2391--2425, 2009.

\bibitem{Maire2011}
P.H. Maire.
\newblock A high-order one-step sub-cell force-based discretization for
  cell-centered lagrangian hydrodynamics on polygonal grids.
\newblock {\em Computers and Fluids}, 46(1):341--347, 2011.

\bibitem{Maire2010}
P.H. Maire.
\newblock A unified sub-cell force-based discretization for cell-centered
  lagrangian hydrodynamics on polygonal grids.
\newblock {\em International Journal of Numerical Methods in Fluids},
  65:1281--1294, 2011.

\bibitem{Maire2009b}
{P.H.} Maire and B.~Nkonga.
\newblock {Multi-scale Godunov-type method for cell-centered discrete
  Lagrangian hydrodynamics}.
\newblock {\em Journal of Computational Physics}, 228:799--821, 2009.

\bibitem{munz94}
C.D. Munz.
\newblock {On Godunov--type schemes for Lagrangian gas dynamics}.
\newblock {\em SIAM Journal on Numerical Analysis}, 31:17--42, 1994.

\bibitem{Noh}
{W.F.} Noh.
\newblock Errors for calculations of strong shocks using artificial viscosity
  and an artificial heat flux.
\newblock {\em Journal of Computational Physics}, 72:78 -- 120, 1987.

\bibitem{Olliver2002}
C.~Olliver-Gooch and M.~Van Altena.
\newblock A high-order-accurate unstructured mesh finite-volume scheme for the
  advection–diffusion equation.
\newblock {\em Journal of Computational Physics}, 181:729 -- 752, 2002.

\bibitem{scovazzi1}
A.~L\'opez Ortega and G.~Scovazzi.
\newblock {A geometrically--conservative, synchronized, flux--corrected remap
  for arbitrary Lagrangian--Eulerian computations with nodal finite elements}.
\newblock {\em Journal of Computational Physics}, 230:6709--6741, 2011.

\bibitem{OsherRS}
S.~Osher and F.~Solomon.
\newblock Upwind difference schemes for hyperbolic systems of conservation
  laws.
\newblock {\em Mathematics of Computation}, 38:339 -- 374, 1982.

\bibitem{Peery2000}
J.S. Peery and D.E. Carroll.
\newblock Multi-material ale methods in unstructured grids,.
\newblock {\em Computer Methods in Applied Mechanics and Engineering},
  187:591--619, 2000.

\bibitem{Roe}
{P.L.} Roe.
\newblock Approximate {Riemann} solvers, parameter vectors, and difference
  schemes.
\newblock {\em Journal of Computational Physics}, 43:357--372, 1981.

\bibitem{Rumsey1993}
C.B. Rumsey, B.~van Leer, and P.L. Roe.
\newblock A multidimensional flux function with application to the euler and
  navier-stokes equations.
\newblock {\em Journal of Computational Physics}, 105:306 -- 323, 1993.

\bibitem{RusanovRS}
V.~V. Rusanov.
\newblock Calculation of interaction of non-steady shock waves with obstacles.
\newblock {\em J. Comput. Math. Phys. USSR}, 1:267 -- 305, 1961.

\bibitem{ShashkovCellCentered}
S.K. Sambasivan, M.J. Shashkov, and D.E. Burton.
\newblock {A finite volume cell-centered Lagrangian hydrodynamics approach for
  solids in general unstructured grids}.
\newblock {\em International Journal of Numerical Methods in Fluids},
  72:770--810, 2013.

\bibitem{ShashkovMultiMat4}
S.K. Sambasivan, M.J. Shashkov, and D.E. Burton.
\newblock {Exploration of new limiter schemes for stress tensors in Lagrangian
  and ALE hydrocodes}.
\newblock {\em Computers and Fluids}, 83:98--114, 2013.

\bibitem{scovazzi2}
G.~Scovazzi.
\newblock {Lagrangian shock hydrodynamics on tetrahedral meshes: A stable and
  accurate variational multiscale approach}.
\newblock {\em Journal of Computational Physics}, 231:8029--8069, 2012.

\bibitem{Smith1999}
{R.W.} Smith.
\newblock {AUSM(ALE)}: a geometrically conservative arbitrary
  lagrangian--eulerian flux splitting scheme.
\newblock {\em Journal of Computational Physics}, 150:268–286, 1999.

\bibitem{stroud}
{A.H.} Stroud.
\newblock {\em Approximate Calculation of Multiple Integrals}.
\newblock Prentice-Hall Inc., Englewood Cliffs, New Jersey, 1971.

\bibitem{toro3}
{V.A.} Titarev and {E.F.} Toro.
\newblock {ADER}: Arbitrary high order {Godunov} approach.
\newblock {\em Journal of Scientific Computing}, 17(1-4):609--618, December
  2002.

\bibitem{TitarevToro}
V.A. Titarev and E.F. Toro.
\newblock {ADER} schemes for three-dimensional nonlinear hyperbolic systems.
\newblock {\em Journal of Computational Physics}, 204:715--736, 2005.

\bibitem{MixedWENO2D}
V.A. Titarev, P.~Tsoutsanis, and D.~Drikakis.
\newblock {WENO schemes for mixed--element unstructured meshes}.
\newblock {\em Communications in Computational Physics}, 8:585--609, 2010.

\bibitem{Toro:2006a}
E.~F. Toro and V.~A. Titarev.
\newblock {Derivative Riemann solvers for systems of conservation laws and ADER
  methods}.
\newblock {\em Journal of Computational Physics}, 212(1):150--165, 2006.

\bibitem{ToroBook}
{E.F.} Toro.
\newblock {\em Riemann Solvers and Numerical Methods for Fluid Dynamics: a
  Practical Introduction.}
\newblock Springer, 2009.

\bibitem{ToroHLLC}
E.F. Toro, M.~Spruce, and W.~Speares.
\newblock Restoration of contact surface in the hll riemann solver.
\newblock {\em Shock Waves}, 4:25 -- 34, 1994.

\bibitem{MixedWENO3D}
P.~Tsoutsanis, V.A. Titarev, and D.~Drikakis.
\newblock {WENO schemes on arbitrary mixed-element unstructured meshes in three
  space dimensions}.
\newblock {\em Journal of Computational Physics}, 230:1585--1601, 2011.

\bibitem{vanLeerRS}
B.~van Leer.
\newblock Toward the ultimate conservative difference scheme. v. a second-order
  sequel to godunov's method.
\newblock {\em Journal of Computational Physics}, 32:101 -- 136, 1979.

\bibitem{ShashkovMultiMat2}
Y.V. Yanilkin, E.A. Goncharov, V.Y. Kolobyanin, V.V. Sadchikov, J.R. Kamm, M.J.
  Shashkov, and W.J. Rider.
\newblock Multi-material pressure relaxation methods for lagrangian
  hydrodynamics.
\newblock {\em Computers and Fluids}, 83:137--143, 2013.

\bibitem{ZhangShu3D}
{Y.T.} Zhang and {C.W.} Shu.
\newblock Third order {WENO} scheme on three dimensional tetrahedral meshes.
\newblock {\em Communications in Computational Physics}, 5:836--848, 2009.

\end{thebibliography}
\bibliographystyle{plain}

\end{document}